\documentclass[12pt]{article}
\usepackage{amsfonts, amsmath, amsthm}
\usepackage[all]{xy}
\usepackage{epic,eepic}

\newtheorem{theorem}{Theorem}[section]
\newtheorem{lemma}[theorem]{Lemma}

\newtheorem{cor}[theorem]{Corollary}
\newtheorem{prop}[theorem]{Proposition}

\def\FF{\mathbb{F}}

\def\QQ{\mathbb{Q}}
\def\RR{\mathbb{R}}
\def\ZZ{\mathbb{Z}}
\def\calC{\mathcal{C}}
\def\calO{\mathcal{O}}

\def\idp{\mathfrak{p}}

\def\alg{\mathrm{alg}}
\def\an{\mathrm{an}}
\def\con{\mathrm{con}}

\def\sep{\mathrm{sep}}
\def\perf{\mathrm{perf}}

\def\beq{\begin{equation}}
\def\eeq{\end{equation}}

\def\fp{\frac{1}{p}}

\def\GK{\Gamma^K}
\def\GL{\Gamma^L}
\def\Gperf{\Gamma^{\perf}}
\def\Gsep{\Gamma^{\sep}}

\def\Galg{\Gamma^{\alg}}

\def\Galgancon{\Galg_{\an, \con}}
\def\Galgcon{\Galg_{\con}}

\def\Gancon{\Gamma_{\an,\con}}
\def\Gcon{\Gamma_{\con}}

\def\GKancon{\GK_{\an,\con}}
\def\GKcon{\GK_{\con}}
\def\GLancon{\GL_{\an,\con}}
\def\GLcon{\GL_{\con}}
\def\Gperfcon{\Gperf_{\con}}
\def\Gsepcon{\Gsep_{\con}}

\def\be{\mathbf{e}}

\def\bv{\mathbf{v}}
\def\bw{\mathbf{w}}
\def\bx{\mathbf{x}}
\def\by{\mathbf{y}}

\DeclareMathOperator{\Frac}{Frac}
\DeclareMathOperator{\Gal}{Gal}

\DeclareMathOperator{\naive}{naive}

\DeclareMathOperator{\rank}{rank}
\DeclareMathOperator{\Span}{SatSpan}
\DeclareMathOperator{\Spec}{Spec}

\newcounter{fixmectr}

\begin{document}

\title{A $p$-adic local monodromy theorem}
\author{Kiran S. Kedlaya \\ Department of Mathematics \\ University
of California, Berkeley \\ Berkeley, CA 94720 \\
kedlaya@math.berkeley.edu}
\date{December 31, 2002}

\maketitle

\begin{abstract}
We produce a canonical filtration for locally free sheaves on an open
$p$-adic annulus equipped with a Frobenius structure.
Using this filtration, we deduce a conjecture of Crew on
$p$-adic differential equations, analogous to Grothendieck's local
monodromy theorem (also a consequence of results of 
Andr\'e and of Mebkhout). Namely, given a finite locally free sheaf
on an open $p$-adic annulus with
a connection and a compatible Frobenius structure, the
module admits a basis over a finite cover of the annulus
on which the connection acts via a nilpotent matrix.
\end{abstract}

\tableofcontents

\section{Introduction}

\subsection{Crew's conjecture on $p$-adic local monodromy}

The role of $p$-adic differential equations in algebraic geometry was
first pursued systematically by Dwork; the modern manifestation of
this role comes via the theory of isocrystals and $F$-isocrystals,
which over a field
of characteristic $p>0$ attempt to play the
part of local systems for the classical 
topology on complex varieties and lisse sheaves for the $l$-adic topology
when $l \neq p$. In order to get a usable theory, however, an additional
``overconvergence''
condition must be imposed, which has no analogue in either
the complex or $l$-adic cases.
For example, the cohomology of the affine line is infinite dimensional
if computed using convergent isocrystals, but has the expected dimension
if computed using overconvergent isocrystals.
This phenomenon was generalized by Monsky and Washnitzer \cite{bib:mw}
into a cohomology theory for smooth affine varieties, and then generalized
further by Berthelot to the theory of rigid cohomology, which has good
behavior for arbitrary varieties (see for example \cite{bib:ber2}).

Unfortunately, the use of overconvergent isocrystals to date has been
hampered by a gap in the
local theory of these objects; for example, it obstructed the proof
of finite dimensionality of Berthelot's rigid cohomology with
arbitrary coefficients (the case of constant coefficients was treated by
Berthelot in \cite{bib:ber1}). This gap can be described as a $p$-adic analogue
of Grothendieck's local monodromy theorem for $l$-adic cohomology.

The best conceivable analogue of Grothendieck's theorem would be that an
$F$-isocrystal becomes a successive extension of trivial isocrystals
after a finite \'etale base extension. Unfortunately, this assertion is not
correct; for example, it fails for the pushforward of the constant isocrystal
on a family of ordinary elliptic curves degenerating to a supersingular
elliptic curve (and for the Bessel isocrystal described in 
Section~\ref{subsec:bessel} over the affine line).

The correct analogue of the local monodromy theorem was formulated
conjecturally by Crew \cite[Section~10.1]{bib:crew2}, and reformulated
in a purely local form by Tsuzuki \cite[Theorem~5.2.1]{bib:tsu3};
we now introduce some terminology and notation needed to describe it.
(These definitions 
are reiterated more precisely in Chapter~\ref{sec:rings}.)
Let $k$ be a field of characteristic $p>0$,
and $\calO$ a finite totally 
ramified extension of a Cohen ring $C(k)$.
The Robba ring $\Gancon$ is defined
as the set of Laurent series over $\calO[\fp]$ which converge on some
open annulus with outer radius 1; its subring $\Gcon$ consists of series
which take integral values on some open annulus with outer radius 1,
and is a discrete valuation ring.
(See Chapter~\ref{sec:aux} to find out where the notations come from.)
We say a ring endomorphism $\sigma: \Gancon \to \Gancon$
is a \emph{Frobenius} for $\Gancon$ if it is a composition power
of a map preserving $\Gcon$ and reducing modulo
a uniformizer of $\Gcon$ to the $p$-th power map.
For example, one can choose
$t \in \Gcon$ whose reduction is a uniformizer in the ring
of Laurent series over $k$, then set $t^\sigma = t^q$.
Note that one cannot hope to define a Frobenius on the ring
of analytic functions on any fixed open annulus with outer radius 1, because
for $\eta$ close to 1, functions on the annulus of inner radius 1 pull
back under $\sigma$
to functions on the annulus of inner radius $\eta^{1/p}$. Instead, one
must work over an ``infinitely thin'' annulus
of radius 1.

Given a ring $R$ in which $p \neq 0$
and an endomorphism $\sigma: R \to R$,
we define a \emph{$\sigma$-module} $M$ as a finite locally free module
equipped with an $R$-linear map $F: M \otimes_{R,\sigma} R \to M$
that becomes an isomorphism over $R[\fp]$; the tensor product notation
indicates that $R$ is viewed as an $R$-module via $\sigma$.
For the rings considered in this paper, a finite locally free module
is automatically free; see Proposition~\ref{prop:free}.
Then $F$ can be viewed as an additive, $\sigma$-linear map $F: M \to M$
that acts on any basis of $M$ by a matrix invertible over $R[\fp]$.

We define
a \emph{$(\sigma, \nabla)$-module} as a $\sigma$-module plus
a connection $\nabla: M \to M \otimes \Omega^1_{R/\calO}$ (that is,
an additive map satisfying the Leibniz rule
$\nabla(c\bv) = c\nabla(\bv) + \bv \otimes dc$) 
that makes the following
diagram commute:
\[
\xymatrix{
M \ar^-{\nabla}[r] \ar^{F}[d] & M \otimes \Omega^1_{R/\calO} \ar^{F \otimes d\sigma}[d] \\
M  \ar^-{\nabla}[r] & M \otimes \Omega^1_{R/\calO}
}
\]
We say a $(\sigma, \nabla)$-module over $\Gancon$ is
\emph{quasi-unipotent} if, after tensoring $\Gancon$ over $\Gcon$ with
a finite extension of $\Gcon$, the module admits a filtration
by $(\sigma, \nabla)$-submodules such that each successive quotient 
admits a basis of elements in the kernel of $\nabla$. (If $k$ is perfect,
one may also insist that the extension of $\Gcon$ be residually separable.)
In these notations, Crew's conjecture is resolved by the following
theorem, which we will prove in a more precise form as Theorem~\ref{thm:main2}.
\begin{theorem}[Local monodromy theorem] \label{thm:monodromy}
Let $\sigma$ be any Frobenius for the Robba ring $\Gancon$. Then
every $(\sigma, \nabla)$-module over $\Gancon$ is 
quasi-unipotent.
\end{theorem}
Briefly put, a $p$-adic differential equation on an annulus with
a Frobenius structure has quasi-unipotent monodromy. It is worth noting
(though not needed in this paper) that for a given $\nabla$, whether
there exists a compatible $F$ does not depend on the choice of the Frobenius
map $\sigma$. This follows from the existence of change of Frobenius
functors \cite[Theorem~3.4.10]{bib:tsu3}.

The purpose of this paper is to establish some structural results
on modules over the Robba ring yielding a proof of Theorem~\ref{thm:monodromy}.
Note that Theorem~\ref{thm:monodromy} itself
has been established independently by
Andr\'e \cite{bib:and} and by Mebkhout \cite{bib:meb}. However,
as we describe in the next section, the methods
in this paper are essentially orthogonal to the methods of those
authors. In fact, the different approaches provide different auxiliary
information, various pieces of which may be of relevance in other contexts.

\subsection{Frobenius filtrations and Crew's conjecture}

Before outlining our approach to Crew's conjecture, we describe by
way of contrast the common features of the work of Andr\'e and Mebkhout.
Both authors build upon the results of
a series of papers by Christol and Mebkhout \cite{bib:cm1}, \cite{bib:cm2},
\cite{bib:cm3}, \cite{bib:cm4} concerning properties of
modules with connection over the Robba ring. Most notably, in \cite{bib:cm4}
is produced a canonical filtration (the ``weight filtration''), defined
whether or not the connection admits a Frobenius structure. Andr\'e
and Mebkhout show (in two different ways)
that when a Frobenius structure is present, the graded pieces of this
filtration can be shown to be quasi-unipotent.

The strategy in this paper is in a sense completely orthogonal
to the aforementioned approach. (For a more detailed comparison between the
various approaches to Crew's conjecture, see the November 2001
Seminaire Bourbaki
talk of Colmez \cite{bib:colmez}.) Instead of isolating the connection
data, we isolate the Frobenius structure and prove a structure theorem
for $\sigma$-modules over the Robba ring. This can be accomplished by
a ``big rings'' argument, where one first proves that $\sigma$-modules
can be trivialized over a large auxiliary ring, and then ``descends''
the construction back to the Robba ring. (Isolating Frobenius is in
a sense natural from the point of view of crystalline cohomology; for
example, this is the approach of Katz in \cite{bib:katz}.)

The model for our strategy of trivializing $\sigma$-modules 
over an auxiliary ring
is the Dieu\-donn\'e-Manin classification of $\sigma$-modules over a
complete discrete valuation ring $R$ of
mixed characteristic $(0,p)$
with algebraically closed residue field. (This classification is a
semilinear analogue of the diagonalization of matrices over an algebraically
closed field, except that here  there is no failure of semisimplicity.)
We give a quick statement here,
deferring the precise formulation to Section~\ref{subsec:dieumanin}.
For $\lambda \in \calO[\fp]$
and $d$ a positive integer, let $M_{\lambda, d}$
denote the $\sigma$-module of rank $d$ over $R[\fp]$ on which $F$ acts
by a basis $\bv_1, \dots, \bv_d$ as follows:
\begin{align*}
F\bv_1 &= \bv_2 \\
&\vdots \\
F\bv_{d-1} &= \bv_d \\
F\bv_d &= \lambda \bv_1.
\end{align*}
Define the \emph{slope} of $M_{\lambda,d}$ to be $v_p(\lambda)/d$.
Then the Dieudonn\'e-Manin classification states (in part) that over $R[\fp]$,
every $\sigma$-module is isomorphic to a direct sum
$\oplus_j M_{\lambda_j, d_j}$, and the slopes that occur do not
depend on the decomposition.

If $R$ is a
discrete valuation ring of mixed characteristic $(0,p)$, 
we may define the slopes of
a $\sigma$-module over $R[\fp]$ as the slopes in a Dieudonn\'e-Manin
decomposition over the maximal unramified extension of the completion
of $R$. 
However, this definition cannot be used immediately over $\Gancon$,
because that ring is not a discrete valuation ring. Instead, we
must first reduce to considering modules over $\Gcon$. Our main
theorem makes it possible to do so.
 Again, we give a quick formulation here and prove a more precise
result later as Theorem~\ref{thm:filt}.
(Note: the filtration in this theorem is similar to what
Tsuzuki \cite{bib:tsu3} calls a ``slope filtration for
Frobenius structures''.)
\begin{theorem} \label{thm:intromain}
Let $M$ be a $\sigma$-module over $\Gancon$.
Then there is a canonical filtration
$0 = M_0 \subset M_1 \subset \cdots \subset M_l = M$
of $M$ by saturated $\sigma$-submodules such that:
\begin{enumerate}
\item[(a)]
each quotient $M_i/M_{i-1}$ is isomorphic over $\Gancon$
to a $\sigma$-module $N_i$ defined over $\Gcon[\fp]$;
\item[(b)]
the slopes of $N_i$ are all equal to some rational number $s_i$;
\item[(c)] $s_1 < \cdots < s_{l}$.
\end{enumerate}
\end{theorem}

The relevance of this theorem to Crew's conjecture is that $(\sigma,
\nabla)$-modules over $\Gcon[\fp]$ with a single slope can be shown to
be quasi-unipotent using a result of Tsuzuki \cite{bib:tsu1}.
The essential case is that of a unit-root $(\sigma, \nabla)$-module 
over $\Gcon$, in which all slopes are 0.
Tsuzuki showed that such modules become isomorphic to a direct sum
of trivial $(\sigma, \nabla)$-modules after
a finite base extension, and even gave precise information about
what extension is needed. This makes it possible to deduce
the local monodromy theorem from Theorem~\ref{thm:intromain}.

\subsection{Applications}

We now describe some consequences of the results of this paper,
starting with some applications via
Theorem~\ref{thm:monodromy}. One set of consequences occurs in
the study of Berthelot's rigid cohomology
(a sort of ``grand unified theory''
of $p$-adic cohomologies). For example, Theorem~\ref{thm:monodromy}
can be used to establish finite dimensionality of rigid cohomology with 
coefficients in an overconvergent $F$-isocrystal; see
\cite{bib:crew2} for the case of a curve and \cite{bib:mefinite} for the
general case. It can also be used to generalize the results
of Deligne's ``Weil II'' to overconvergent $F$-isocrystals; 
this is carried out in \cite{bib:meweil}, building on
work of Crew \cite{bib:crew1}, \cite{bib:crew2}.
In addition, it can be used to treat certain types of ``descent'',
such as Tsuzuki's full faithfulness conjecture \cite{bib:tsu4},
which asserts
that convergent morphisms between overconvergent
$F$-isocrystals are themselves overconvergent; see \cite{bib:me7}.

Another application of Theorem~\ref{thm:monodromy} has been found by
Berger  \cite{bib:berger}, who
has exposed a close relationship between $F$-isocrystals and $p$-adic
Galois representations. In particular, he shows that
Fontaine's ``conjecture de monodromie
$p$-adique'' for $p$-adic Galois representations (that every
de~Rham representation is potentially semistable) follows from 
Theorem~\ref{thm:monodromy}.

Further applications of Theorem~\ref{thm:intromain} exist that do not
directly pass through Theorem~\ref{thm:monodromy}.
For example, Andr\'e (work in progress) has formulated a $q$-analogue of Crew's
conjecture, in which the single differential equation is replaced by
a formal deformation. He has established this analogue using
Theorem~\ref{thm:filt} plus a
$q$-analogue of Tsuzuki's unit-root theorem (Proposition~\ref{prop:tsuzuki}),
and has deduced a finiteness theorem for
rigid cohomology of $q$-$F$-isocrystals. (It should also be possible to 
obtain these results
using a $q$-analogue of the Christol-Mebkhout theorem, and indeed
Andr\'e and di~Vizio have made progress in this direction;
however, at the
time of this writing, some technical details had not yet been worked
out.)

We also plan to establish, in a subsequent paper, a
conjecture of Shiho \cite[Conjecture~3.1.8]{bib:shiho},
on extending overconvergent $F$-isocrystals to log-$F$-isocrystals after
a generically \'etale base change. This result appears to require
a more sophisticated analogue of Theorem~\ref{thm:filt}, in which the
``one-dimensional'' Robba ring is replaced by a ``higher-dimensional''
analogue. (One might suspect that this
conjecture should follow from Theorem~\ref{thm:monodromy} and some
clever geometric arguments, but the situation appears to be more subtle.)
Berthelot (private communication) has suggested that
a suitable result in this direction
may help in constructing Grothendieck's six operations in the category
of arithmetic $\mathcal{D}$-modules, which would provide a
$p$-adic analogue of the constructible sheaves in \'etale cohomology.

\subsection{Structure of the paper}
\label{subsec:outlineproof}

We now outline the strategy of the proof of Theorem~\ref{thm:intromain},
and in the process
describe the structure of the paper. We note in passing that
some of the material appears
in the author's doctoral dissertation \cite{bib:methesis}, written
under Johan de~Jong, and/or in a sequence
of unpublished
preprints \cite{bib:me3}, \cite{bib:me4}, \cite{bib:me5}, \cite{bib:me6}.
However, the present document avoids
any logical dependence on unpublished results.

In Chapter~\ref{sec:rings}, we recall some of the basic rings of
the theory of $p$-adic differential equations; they include the
Robba ring, its integral subring and the completion of the latter
(denoted the ``Amice ring'' in some sources). In Chapter~\ref{sec:aux},
we construct some less familiar rings by
augmenting the classical constructions.
These augmentations
are inspired by (and in some cases identical to)
the auxiliary rings used by de Jong \cite{bib:dej1} in his extension
to equal characteristic of Tate's theorem \cite{bib:tate}
on $p$-divisible groups
over mixed characteristic discrete valuation rings. (They also resemble
the ``big rings'' in Fontaine's theory of $p$-adic Galois representations.)
 In particular,
a key augmentation, denoted $\Galgancon$, is a sort of
``maximal unramified extension'' of the Robba ring, and a great
effort is devoted to showing that it shares the B\'ezout property
with the Robba ring; that is, every finitely generated
ideal in $\Galgancon$ is
principal. (This chapter is somewhat technical; we suggest that
the reader skip it on first reading, and refer back to it as needed.)

With these augmented rings in hand, in Chapter~\ref{sec:special}
we show that every $\sigma$-module over the Robba ring can be
equipped with a canonical filtration over $\Galgancon$; this amounts
to an ``overconvergent'' analogue of the Dieudonn\'e-Manin classification.
From this filtration we read off a sequence of slopes, which in
case we started with a quasi-unipotent $(\sigma, \nabla)$-module
agree with the slopes of Frobenius on a nilpotent basis; the Newton
polygon with these slopes is called the \emph{special Newton polygon}.

By contrast, in Chapter~\ref{sec:generic}, we associate to a
$(\sigma, \nabla)$-module over $\Gcon$ the Frobenius slopes produced
by the Dieudonn\'e-Manin classification. The Newton polygon with
these slopes is called the \emph{generic Newton polygon}. Following
\cite{bib:dej1}, we construct
some canonical filtrations associated with the generic Newton
polygon. This chapter is logically independent of Chapter~\ref{sec:special}
except at its conclusion, when the two notions of Newton polygon are
compared. In particular, we show that the special Newton polygon lies
above the generic Newton polygon with the same endpoint, and obtain
additional structural consequences in case the Newton polygons coincide.

Finally, in Chapter~\ref{sec:conclude}, we take the ``generic''
and ``special'' filtrations, both defined over large auxiliary rings,
and descend them back to the Robba ring itself.
The basic strategy here is to separate positive and negative powers
of the series parameter, using the auxiliary filtrations to guide
the process. Starting with a $\sigma$-module over the Robba ring, the
process yields a $\sigma$-module over $\Gcon$ whose generic and
special Newton polygons coincide. The structural consequences mentioned
above yield Theorem~\ref{thm:intromain}; by
applying Tsuzuki's theorem on unit-root $(\sigma, \nabla)$-modules
(Proposition~\ref{prop:tsuzuki}), we deduce a precise form
of Theorem~\ref{thm:monodromy}.

\subsection{An example: the Bessel isocrystal}
\label{subsec:bessel}

To clarify the remarks of the previous section, we include a classical
example to  illustrate the different structures we have described,
especially the generic and special Newton polygons.
Our example is the Bessel isocrystal, first studied
by Dwork \cite{bib:dwork}; our description is a summary of
the discussion of Tsuzuki \cite[Example~6.2.6]{bib:tsu3}
(but see also Andr\'e \cite{bib:and1}).

Let $p$ be an odd prime, and put $\calO = \ZZ_p[\pi]$,
where $\pi$ is a $(p-1)$-st root of $-p$.
Choose $\eta < 1$, and let $R$ be the ring of Laurent series over $\calO$
convergent for $|t| > \eta$ in the variable.
Let $\sigma$ be the Frobenius lift on $\calO$ such that $t^\sigma
= t^p$. Then for suitable $\eta$, there exists a $(\sigma,
\nabla)$-module $M$ of rank two over $R$
admitting a basis $\bv_1, \bv_2$ such that
\begin{align*}
F\bv_1 &= A_{11} \bv_1 + A_{12} \bv_2 \\
F\bv_2 &= A_{21} \bv_1 + A_{22} \bv_2 \\
\nabla \bv_1 &= t^{-2} \pi^{2} \bv_2 \otimes dt \\
\nabla \bv_2 &= t^{-1} \bv_1 \otimes dt.
\end{align*}
Moreover, the matrix $A$ satisfies
\[
\det(A) = p \qquad \mbox{and} \qquad A \equiv \begin{pmatrix}
1 & 0 \\ 0 & 0 \end{pmatrix} \pmod{p}.
\]
It follows that the two generic Newton slopes are nonnegative (because
the entries of $A$ are integral), their sum is 1 (by the determinant
equation), and the smaller of the two is zero (by the congruence).
Thus the generic Newton slopes are 0 and 1.

On the other hand, $M$ becomes unipotent after
adjoining a square root of $t$ to the Robba ring. More precisely,
if $y = (t/4)^{1/2}$, define
\[
f_{\pm} = 1 + \sum_{n=1}^\infty (\pm 1)^n \frac{(1 \times 3 \times
\cdots \times (2n-1))^2}{(8 \pi)^n n!} y^n
\]
and set
\[
\bw_{\pm} = f_{\pm} \be_1 + \left( y \frac{df_{\pm}}{dy} + \left(
\frac{1}{2} \mp \pi y^{-1} \right) f_{\pm} \right) \be_2.
\]
Then
\[
\nabla \bw_{\pm} = \left( \frac{-1}{2} \pm \pi y^{-1} \right) \bw_{\pm}
\otimes \frac{dy}{y}.
\]
Using the compatibility between the Frobenius and connection structures,
we deduce that
\[
F\bw_{\pm} = \alpha_{\pm} y^{-(p-1)/2} \exp(\pm \pi(y^{-1} - y^{-\sigma}))
\]
for some $\alpha_+, \alpha_- \in \calO[\fp]$ 
with $\alpha_+ \alpha_- = 2^{1-p} p$.
By the invariance of Frobenius under the automorphism $y \to -y$ of
$\Gancon[y]$, we deduce
that $\alpha_+$ and $\alpha_-$ have the same valuation.

It follows from this discussion that
$M$ is unipotent over $\Gancon[y]$, and 
the two slopes of the special Newton polygon are equal, necessarily
to $1/2$ since their sum is 1. In particular, the special Newton polygon
lies above the generic Newton polygon and has the same endpoint, but
the two polygons are not equal in this case.

\subsection*{Acknowledgments}
The author was supported by a Clay Mathematics Institute Liftoffs
grant and a National Science Foundation Postdoctoral Fellowship.
Thanks to the organizers of the Algorithmic Number Theory program
at MSRI, the Arizona Winter School in Tucson, and
the Dwork Trimester in Padua for their
hospitality, and to Laurent Berger, Pierre Colmez,
Johan de~Jong and the referee for helpful suggestions.

\section{A few rings}
\label{sec:rings}

In this chapter, we set some notations and conventions,
and define some of the basic rings
used in the local
study of $p$-adic differential equations. We also review
the basic properties of rings in which every finitely generated ideal
is principal (B\'ezout rings), and introduce
$\sigma$-modules and $(\sigma, \nabla)$-modules.

\subsection{Notations and conventions}

Recall that for every field $K$ of characteristic $p>0$, there exists
a complete discrete valuation ring with
fraction field of characteristic 0, maximal ideal generated by $p$, and
residue field isomorphic to $K$,
and that this ring is unique up to noncanonical
isomorphism. Such a ring is called a \emph{Cohen ring} for $K$;
see \cite{bib:bou} for the basic properties of such rings. If $K$ is perfect,
the Cohen ring is unique up to canonical isomorphism, and coincides with
the ring $W(K)$ of Witt vectors over $K$. (Note in passing: for $K$ perfect,
we use brackets to denote Teichm\"uller lifts into $W(K)$.)

Let $k$ be a field of characteristic $p>0$, and $C(k)$ a Cohen ring
for $k$. 
Let $\calO$ be a finite totally ramified extension of $C(k)$,
let $\pi$ be a uniformizer of $\calO$, and
fix once and for all a ring endomorphism $\sigma_0$
on $\calO$ lifting the absolute Frobenius $x \mapsto x^p$ on $k$.
Let $q = p^f$ be a power of $p$ and put $\sigma = \sigma_0^f$.
(In principle, one could dispense with $\sigma_0$ and simply take
$\sigma$ to be any ring endomorphism lifting the $q$-power Frobenius.
We will eschew this additional level of generality so that we can
invoke the results of \cite{bib:tsu1}.)
Let $v_p$ denote the valuation on $\calO[\fp]$ normalized so that
$v_p(p) = 1$, and let $|\cdot|$ denote the norm on $\calO[\fp]$
given by $|x|= p^{-v_p(x)}$.

Let $\calO_0$ denote the fixed ring of $\calO$ under $\sigma$. 
If $k$ is algebraically
closed, then the equation $u^\sigma = (\pi^\sigma/\pi)u$ in $u$ has
a nonzero solution modulo $\pi$, and so by a variant of Hensel's lemma 
(see Proposition~\ref{prop:hensel}) has a nonzero solution in $\calO$.
Then $(\pi/u)$ is a uniformizer of
$\calO$ contained in $\calO_0$, and hence
$\calO_0$ has the same value group as $\calO$. That being the case,
we can and will take $\pi \in \calO_0$ in case $k$ is algebraically closed.

We wish to alert the reader to several notational conventions in force
throughout the paper. The first of these is ``exponent consolidation''.
The expression $(x^{-1})^\sigma$,
for $x$ a ring element or matrix and $\sigma$ a ring endomorphism,
will often be abbreviated $x^{-\sigma}$.
This is not to be confused with $x^{\sigma^{-1}}$; the former is the image
under $\sigma$ of the multiplicative inverse of $x$, the latter is the preimage
of $x$ under $\sigma$ (if it exists).
Similarly, if $A$ is a matrix, then
$A^T$ will denote the transpose of $A$, and the expression $(A^{-1})^T$
will be abbreviated $A^{-T}$.

We will use the summation notation
$\sum_{i=m}^n f(i)$ in some cases where $m > n$, in which case we mean
0 for $n=m-1$ and $- \sum_{i=n+1}^{m-1} f(i)$ otherwise. 
The point of this convention is that
$\sum_{i=m}^{n} f(i) = f(n) + \sum_{i=m}^{n-1} f(i)$ for all $n \in \ZZ$. 

We will perform a number of calculations involving matrices; these will
always be $n \times n$
matrices unless otherwise specified. Also, $I$ will denote the $n \times n$
identity matrix over any ring.

\subsection{Valued fields}

Let $k((t))$ denote the field of Laurent series over $k$.
Define a \emph{valued field} to be an algebraic extension $K$ of $k((t))$
for which there exist subextensions $k((t)) \subseteq L
\subseteq M \subseteq N \subseteq K$ such that:
\begin{enumerate}
\item[(a)] $L = k^{1/p^m}((t))$ for some $m \in \{0,1, \dots, \infty\}$;
\item[(b)] $M = k_M((t))$ for some separable algebraic extension $k_M/k^{1/p^m}$;
\item[(c)] $N = M^{1/p^n}$ for some $n \in \{0,1, \dots, \infty\}$;
\item[(d)] $K$ is a separable totally ramified algebraic extension of $N$.
\end{enumerate}
(Here $F^{1/p^\infty}$ means the perfection of the field $F$,
and $K/N$ totally ramified means that $K$ and $N$ have the same value group.)
Note that $n$ is uniquely determined by $K$: it is the largest integer $n$
such that $t$ has a $p^n$-th root in $K$. If $n < \infty$ (e.g., if 
$K/k((t))$ is finite), then $L,M,N$ are also determined by $K$:
$k_M^{1/p^n}$ must be the integral closure of $k$ in $K$, which
determines $k_M$, and $k^{1/p^m}$
must be the maximal purely inseparable subextension of $k_M/k$.

The following proposition shows that the definition of a valued
field is only restrictive if $k$ is imperfect. It also
guides the construction
of the rings $\Gamma^K$ in Section~\ref{subsec:cohen}.
\begin{prop} \label{prop:valued}
If $k$ is perfect, then any 
algebraic extension $K/k((t))$ is a valued field.
\end{prop}
\begin{proof}
Normalize the valuation $v$ on $k((t))$ so that $v(t) = 1$. 
Let $k_M$ be the integral closure of $k$ in $K$,
and define $L = k((t))$ and $M = k_M((t))$.
Then (a) holds for $m=0$ and (b) holds because $k$ is perfect.

Let $n$
be the largest nonnegative integer such that $t$ has a $p^n$-th
root in $K$, or $\infty$ if there is no largest integer.
Put
\[
N = \bigcup_{i=0}^\infty \left(K \cap M^{1/p^i} \right).
\]
Since $t^{1/p^i} \in K$ for all
$i \leq n$ and $k_M$ is perfect,
we have $M^{1/p^n} \subseteq N$. On the other hand,
suppose $x^{1/p^i} \in (K \cap M^{1/p^i}) \setminus (K \cap M^{1/p^{i-1}})$,
that is, $x \in M$ has a $p^i$-th root in $K$ but has no $p$-th root in
$M$. Then $v(x)$ is relatively prime to $p$, so we can find integers
$a$ and $b$ such that $y = x^a/t^{bp^i} \in M$ has a $p^i$-th root in
$K$ and $v(y)=1$. We can write every element of $M$ uniquely
as a power series in $y$, so every element of $M$ has a $p^i$-th root in
$K$. In particular, $t$ has a $p^i$-th root in $K$, so $i \leq n$.
We conclude that $N = M^{1/p^n}$, verifying (c).

If $y \in K^p \cap N$, then $y = z^p$ for some $z \in K$ and
$y^{p^i} \in M$ for some $i$. Then $z^{p^{i+1}} \in M$, so $z \in N$.
We thus have $K^p \cap N = N^p$, so $K/N$ is separable. To verify that
$K/N$ is totally ramified, let $K_0$ be any finite subextension of $K/k((t))$
and let $U$ be the maximal unramified subextension of $K_0/(K_0 \cap N)$. 
We now recall two basic facts from \cite{bib:serre} about finite
extensions of fields complete with respect to discrete valuations:
\begin{enumerate}
\item[1.] $K_0/U$ is totally ramified, because $K_0/(K_0 \cap N)$ and its residue
field extension are both separable.
\item[2.] There is a unique unramified extension of $K_0 \cap N$ 
yielding any specified separable residue field extension.
\end{enumerate}
Since $K_0 \cap N$ is a power series field, we can make unramified extensions
of $K_0 \cap N$ with any specified residue field extension by
extending the constant field $K_0 \cap k_M$. By the second
assertion above, $U/(K_0 \cap N)$ must then be a constant field extension.
However,
$k_M$ is integrally closed in $K$, so $U = K_0 \cap N$ and 
$K_0/(K_0 \cap N)$ is totally ramified
by the first assertion above.
Since $K$ is the union of its finite subextensions over $k((t))$, we conclude
that $K/N$ is totally ramified, verifying (d).
\end{proof}
The proposition fails for $k$ imperfect, as there are separable
extensions of $k((t))$ with inseparable residue field extensions.
For example, if $c$ has no $p$-th root in $k$, then
$K = k((t))[x]/(x^p-x-ct^{-p})$ is separable over $k((t))$
but induces an inseparable residue field extension. Thus $K$ cannot be a valued
field, as valued fields contain their residue field extensions.
What is true for any $k$ is that valued fields finite and
normal over $k((t))$ form a cofinal subset of the finite extensions
of $k((t))$.

We denote the perfect and algebraic closures of $k((t))$ by
$k((t))^{\perf}$ and $k((t))^{\alg}$; these are both valued fields.
We denote the separable closure of $k((t))$ by 
$k((t))^{\sep}$; this is a valued field only if $k$ is perfect, as we
saw above.

We say a valued field $K$ is \emph{nearly separable} if it is a separable
extension of $k^{1/p^i}((t))$ for some integer $i$. (That is, any
inseparability is concentrated in the constant field.) This definition allows
to approximate certain separability assertions for $k$ perfect in the case
of general $k$, where some separable extensions of $K$ fail to be valued 
fields. For example, 
\[
k^{1/p}((t))[x]/(x^p - x - c t^{-p}) =
k^{1/p}((t))[x]/((x - c^{1/p}t^{-1})^p - (x - c^{1/p}t^{-1}) - c^{1/p} t^{-1})
\]
is a nearly separable valued field. In general, given any separable extension
of $k((t))$, taking its compositum with $k^{1/p^i}((t))$ for sufficiently
large $i$ yields a a nearly separable valued field.

\subsection{The ``classical'' case $K = k((t))$}
\label{subsec:classical}

The definitions and results of Chapter~\ref{sec:aux} generalize
previously known definitions and results in the key case $K = k((t))$.
We treat this case first, both to allow readers familiar
with the prior constructions to get used to the notations of this paper,
and to provide a base on which to build additional rings in 
Chapter~\ref{sec:aux}.

For $K = k((t))$,
let $\Gamma^K$ be the ring of bidirectional power series
$\sum_{i \in \ZZ} x_i u^i$, with $x_i \in \calO$, such that
$|x_i| \to 0$ as $i \to -\infty$. Note that $\Gamma^K$ is a
discrete valuation ring complete under the $p$-adic
topology, whose residue field is isomorphic to $K$ via the map
$\sum x_i u^i \mapsto \sum \overline{x_i} t^i$ (using the bar to
denote reduction modulo $\pi$). In particular, if $\pi = p$, then $\Gamma^K$ is
a Cohen ring for $K$.

For $n$ in the value group of $\calO$,
we define the ``na\"\i ve partial valuations''
\[
v^{\naive}_n\left(\sum x_i u^i\right) = \min_{v_p(x_i) \leq n} \{i\},
\]
taking the maximum to be $+\infty$ if no such $i$ exist.
These partial valuations obey some basic rules:
\begin{align*}
v_n(x+y) &\geq \min\{v_n(x), v_n(y) \} \\
v_n(xy) &\geq \min_m \{v_m(x) + v_{n-m}(y)\}
\end{align*}
In both cases, equality always holds if the minimum is achieved exactly once.

Define the \emph{levelwise topology} on $\Gamma^K$ by declaring
the collection of sets
\[
\{x \in \Gamma^K: v^{\naive}_n(x) > c\},
\]
for each $c \in \QQ$ and each $n$ in the value group of $\calO$,
to be a neighborhood basis of 0.
The levelwise topology is finer
than the $\pi$-adic topology, and the Laurent polynomial ring
$\calO[u, u^{-1}]$ is dense in $\Gamma^K$ under the levelwise topology;
thus any levelwise continuous endomorphism of $\Gamma^K$ is determined
by the image of $t$.

The ring $\GKcon$ is isomorphic to the subring of $\Gamma^K$ consisting of
those series $\sum_{i \in \ZZ} x_i u^i$ satisfying the more stringent
convergence condition 
\[
\liminf_{i \to -\infty} \frac{v_p(x_i)}{-i} > 0.
\]
It is also a discrete valuation ring with residue field $K$, but is
not $\pi$-adically complete.

Using the na\"\i ve partial valuations, we can define actual valuations
on certain subrings of $\GKcon$. Let $\Gamma^K_{r, \naive}$ be the
set of $x = \sum x_i u^i$ in $\GKcon$ such that $\lim_{n \to \infty}
r v^{\naive}_n(x) + n = \infty$; the union of the subrings over all
$r$ is precisely $\GKcon$.
(Warning: the rings $\Gamma^K_{r, \naive}$ for individual $r$ are not
discrete valuation rings,
even though their union is.)
On this subring, we have the function
\[
w^{\naive}_r(x) = \min_n \{ r v^{\naive}_n(x) + n \}
= \min_i \{ r i + v_p(x_i) \}
\]
which can be seen to be a nonarchimedean valuation as follows. It is clear
that $w^{\naive}_r(x+y) \geq \min \{ w^{\naive}_r(x), w^{\naive}_r(y)\}$
from the inequality $v_n(x+y) \geq \min\{v_n(x),v_n(y)\}$.
As for multiplication, it is equally clear that
$w^{\naive}_r(xy) \geq w^{\naive}_r(x) + w^{\naive}_r(y)$; the subtle part
is showing equality. Choose $m$ and $n$ minimal so that
$w^{\naive}_r(x) = r v^{\naive}_m(x) + m$
and $w^{\naive}_r(y) = r v^{\naive}_n(y) + n$; then
\[
r v^{\naive}_{m+n}(xy) + m + n 
\geq \min_i \{ r v^{\naive}_i(x) + i + r v^{\naive}_{m+n-i}(y) + m+n-i\}.
\]
The minimum occurs only once, for $i = m$, so equality holds,
yielding $w^{\naive}_r(xy) = w^{\naive}_r(x) + w^{\naive}_r(y)$.

Since $w^{\naive}_r$ is a valuation, we have a corresponding norm
$|\cdot|^{\naive}_r$ given by $|x|^{\naive}_r
= p^{-w^{\naive}_r(x)}$.
This norm admits a geometric interpretation: the ring
$\Gamma^K_{r, \naive}[\fp]$ consists of power series which converge
and are bounded for $p^{-r} \leq |u| < 1$, where
$u$ runs over all finite extensions of $\calO[\fp]$.
Then $|\cdot|^{\naive}_r$ coincides with the supremum norm on the circle
$|u| = p^{-r}$.

Recall that $\sigma_0: \calO \to \calO$ is a lift of the $p$-power Frobenius
on $k$. We choose an extension of $\sigma_0$ to a levelwise continuous
endomorphism of $\Gamma^K$ that maps $\GKcon$ into itself, and which
lifts the $p$-power Frobenius on $k((t))$. In other words, choose $y
\in \GKcon$ congruent to $u^p$ modulo $\pi$, and define $\sigma_0$ by
\[
\sum_i a_i u^i \mapsto \sum_i a_i^{\sigma_0} y^i.
\]
Define $\sigma = \sigma_0^f$, where $f$ is again given by $q = p^f$.

Let $\GKancon$ be the ring of bidirectional power series
$\sum_i x_i u^i$, now with $x_i \in \calO[\fp]$, satisfying
\[
\liminf_{i \to -\infty} \frac{v_p(x_i)}{-i} > 0,
\qquad
\liminf_{i \to +\infty} \frac{v_p(x_i)}{i} \geq 0.
\]
In other words, for any series $\sum_i x_i u^i$ in $\GKancon$,
there exists $\eta > 0$ such that the series converges for $\eta
\leq |u| < 1$. This ring is
commonly known as the \emph{Robba ring}.
It contains $\GKcon[\fp]$, as the subring of functions which are
analytic and bounded on some annulus
$\eta \leq |u| < 1$, but neither contains nor is contained
in $\Gamma^K$.

We can view $\Gamma^K$ as the $\pi$-adic completion of $\GKcon$; our next
goal is to identify $\GKancon$ with a certain completion
of $\GKcon[\fp]$.
Let $\GK_{\an,r,\naive}$ be the ring of
series $x \in \GKancon$ such that $r v^{\naive}_n(x) + n \to \infty$
as $n \to \pm \infty$. Then $\GKancon$ is visibly the union of
the rings $\GK_{\an,r,\naive}$ over all $r>0$.
We equip $\GK_{\an,r,\naive}$ with the Fr\'echet topology
for the norms $|\cdot|^{\naive}_s$ for $0 < s \leq r$. 
These topologies
are compatible with the embeddings $\GK_{\an,r,\naive} \hookrightarrow
\GK_{\an,s,\naive}$ for $0 < s < r$ (that is, the topology on
$\GK_{\an,r,\naive}$ coincides with the subspace topology for the
embedding), so by taking the direct limit we
obtain a topology on $\GKancon$, which by abuse of language we
will also call the Fr\'echet topology. (A better name might be
the ``limit-of-Fr\'echet topology''.)
Note that $\GK_{r,\naive}[\fp]$ is dense in $\GK_{\an,r,\naive}$ for each
$r$, so $\GKcon[\fp]$ is dense in $\GKancon$.
\begin{prop} \label{prop:freccomp}
The ring $\GK_{\an,r,\naive}$ is complete (for the Fr\'echet topology).
\end{prop}
\begin{proof}
Let $\{x_i\}$ be a Cauchy sequence for the Fr\'echet topology,
consisting of elements of $\GK_{r, \naive}[\fp]$.
That means
that for $0 < s \leq r$ and any $c>0$, there exists $N$ such that
$w^{\naive}_s(x_i - x_j) > c$ for $i,j \geq N$. Write $x_i
= \sum_{l} x_{i,l} u^l$; then for each $l$, $\{x_{i,l}\}$ forms a
Cauchy sequence. More precisely,
for $i,j \geq N$, we have
\[
s l + v_p( x_{i,l} - x_{j,l}) > c.
\]
Since $\calO$ is complete, we can take the limit
$y_l$ of $\{x_{i,l}\}$, and it will satisfy $s l + v_p(x_{i,l} - y_l) > c$
for $i \geq N$.
Thus if we can show $y = \sum_l y_l u^l \in \GK_{\an,r,\naive}$, then
$\{x_i\}$ will converge to $y$ under $|\cdot|^{\naive}_s$ for each $s$.

Choose $s<r$ and $c>0$; we must show that $sl + v_p(y_l) \geq c$ for all
but finitely many $l$. There exists $N$ such that
$sl + v_p(x_{i,l} - y_l) \geq c$ for $i \geq N$. Choose
a single such $i$; then
\begin{align*}
sl + v_p(y_l) &\geq \min\{sl + v_p(x_{i,l} - y_l), sl + v_p(x_{i,l})\} \\
&\geq \min\{c, sl + v_p(x_{i,l})\}.
\end{align*}
Since $x_i \in \GK_{r, \naive}[\fp]$, 
$sl + v_p(x_{i,l}) \geq c$ for all but finitely many $l$.
For such $l$, we have $sl + v_p(y_l) \geq c$, as desired.
Thus $y \in \GK_{\an,r,\naive}$; as noted earlier,
$y$ is the limit of $\{x_i\}$ under each
$|\cdot|^{\naive}_s$, and so is the Fr\'echet limit.

We conclude that each Cauchy sequence with elements in
$\GK_{r,\naive}[\fp]$ has a limit in $\GK_{\an,r,\naive}$. Since
$\GK_{r,\naive}[\fp]$ is dense in $\GK_{\an,r,\naive}$ (one sequence
converging to $\sum_i x_i u^i$ is simply $\{\sum_{i \leq j} x_i 
u^i\}_{j=0}^\infty)$, $\GK_{\an,r,\naive}$ is complete for the
Fr\'echet topology, as desired.
\end{proof}

Unlike $\Gamma^K$ and $\GKcon$, $\GKancon$ is not a discrete valuation
ring. For one thing, $\pi$ is invertible in $\GKancon$. For another,
there are plenty of noninvertible elements of $\GKancon$, such as
\[
\prod_{i=1}^\infty \left( 1 - \frac{u^{p^i}}{p^i} \right).
\]
For a third, $\GKancon$ is not Noetherian; the ideal $(x_1, x_2, \dots)$,
where
\[
x_j = \prod_{i=j}^\infty \left( 1 - \frac{u^{p^i}}{p^i} \right),
\]
is not finitely generated.
However, as long as we restrict to ``finite''
objects, $\GKancon$ behaves well:
a theorem of Lazard \cite{bib:laz} (see also 
\cite[Proposition~4.6]{bib:crew2} and our own Section~\ref{subsec:anbezout2})
states that $\GKancon$
is a \emph{B\'ezout ring}, which is to say every finitely generated ideal is
principal. 

For $L$ a finite extension of $k((t))$, we have $L \cong k'((t'))$ for
some finite extension $k'$ of $k$ and some uniformizer $t'$, 
so one could define $\GL$, $\GLcon$, $\GLancon$ abstractly
as above. However, a better strategy will be to construct these in a
``relative'' fashion; the results will be the same as the abstract rings,
but the relative construction will give us more functoriality, and
will allow us to define $\GL, \GLcon, \GLancon$ even when $L$ is
an infinite algebraic extension of $k((t))$. We return to this approach
in Chapter~\ref{sec:aux}.

The rings defined above occur in numerous other contexts, so it is perhaps
not surprising that there are several sets of notation for them in
the literature.
One common set is
\[
\mathcal{E} = \Gamma^{k((t))}[{\textstyle \fp}], \qquad
\mathcal{E}^\dagger = \Gamma^{k((t))}_{\con}[{\textstyle \fp}], \qquad
\mathcal{R} = \Gamma^{k((t))}_{\an,\con}.
\]
The peculiar-looking notation we have set up will make it easier to deal
systematically with a number of additional rings we will be defining
in Chapter~\ref{sec:aux}.

\subsection{More on B\'ezout rings}

Since $\GKancon$ is a B\'ezout ring, as
are trivially all discrete valuation rings,
it will be useful to record some consequences of the B\'ezout property.

\begin{lemma} \label{lem:det1}
Let $R$ be a B\'ezout ring.
If $x_1, \dots, x_n \in R$ generate the unit ideal,
then there exists a matrix $A$ over $R$ with determinant
$1$ such that $A_{1i} = x_i$ for $i=1, \dots, n$.
\end{lemma}
\begin{proof}
We prove this by induction on $n$, the case $n=1$ being evident.
Let $d$ be a generator of $(x_1, \dots, x_{n-1})$. By the induction
hypothesis, we can find an $(n-1) \times (n-1)$ matrix $B$ of determinant
1 such that $B_{1i} = x_i/d$ for $i=1, \dots, n-1$; extend $B$ to an
$n \times n$ matrix by setting $B_{nn} = 1$ and $B_{in} = B_{ni} = 0$ for
$i=1, \dots, n-1$. Since $(d, x_n) = (x_1, \dots, x_n)$ is the unit ideal,
we can find $e, f \in R$ such that $de - f x_n = 1$. Define the matrix
\[
C = \begin{pmatrix}
d & 0 & \cdots & 0 & x_n \\
0 & 1 & \cdots & 0 & 0 \\
\vdots & & \ddots & & \vdots \\
0 & 0 & \cdots & 1 & 0 \\
f & 0 & \cdots & 0 & e
\end{pmatrix}; \qquad
\mbox{that is,} \qquad
C_{ij} = \begin{cases}
d & i=j=1 \\
1 & 2 \leq i = j \leq n-1 \\
e & i=j=n \\
x_n & i=1, j=n \\
f & i=n, j = 1 \\
0 & \mbox{otherwise.}
\end{cases}
\]
Then we may take $A = CB$.
\end{proof}

Given a finite free module $M$ over a domain $R$, we may regard $M$
as a subset of $M \otimes_R \Frac(R)$; given a subset $S$ of $M$, 
we define the \emph{saturated span} $\Span(S)$ of $S$
as the intersection of $M$ with the $\Frac(R)$-span
of $S$ within $M \otimes_R \Frac(R)$. Note that the following lemma does
not require any finiteness condition on $S$.
\begin{lemma} \label{lem:satspan}
  Let $M$ be a finite free module over a B\'ezout domain $R$. Then for
any subset $S$ of $M$, 
$\Span(S)$ is free and admits a basis that extends to a basis of $M$;
in particular, $\Span(S)$ is a direct summand of $M$.
\end{lemma}
\begin{proof}
We proceed by induction on the rank of $M$, the case of rank 0 being trivial.
Choose a basis $\be_1, \dots, \be_n$ of $M$. If $S$ is empty, there is
nothing to prove; otherwise, choose $\bv \in S$ and write $\bv = \sum_i c_i
\be_i$. Since $R$ is a B\'ezout ring, we can choose a generator $r$
of the ideal $(c_1, \dots, c_n)$. Put $\bw = \sum_i (c_i/r) \be_i$;
then $\bw \in \Span(S)$ since $r\bw = \bv$. By Lemma~\ref{lem:det1},
there exists an invertible matrix $A$ over $R$ with $A_{1i} = c_i/r$.
Put $\by_j = \sum_i A_{ji} \be_i$ for $j=2, \dots, n$; then
$\bw$ and the $\by_j$ form a basis of $M$ (because $A$ is invertible), so
$M/\Span(\bw)$ is free.
Thus the induction hypothesis applies to $M/\Span(\bw)$, where the
saturated span of the image of $S$ admits a basis $\bx_1, \dots, \bx_r$. 
Together
with $\bw$, any lifts of $\bx_1, \dots, \bx_r$ to $M$ form a basis of
$\Span(S)$ that extends to a basis of $M$, as desired.
\end{proof}

Note that the previous lemma immediately implies that every finite torsion-free
module over $R$ is free. (If $M$ is torsion-free and $\phi:
F \to M$ is a surjection
from a free module $F$, then $\ker(\phi)$ is saturated, so $M \cong
F/\ker(\phi)$ is free.) A similar argument yields the following vitally
important fact.
\begin{prop} \label{prop:free}
Let $R$ be a B\'ezout domain. Then every finite locally free module
over $R$ is free.
\end{prop}
\begin{proof}
Let $M$ be a finite locally free module over $R$. Since $\Spec R$ is connected,
the localizations of $M$ all have the same rank $r$. Choose a surjection
$\phi: F \to M$, where $F$ is a finite free $R$-module, and let $N
= \Span(\ker(\phi))$. Then we have a surjection $M \cong F/\ker(\phi)
\to F/N$, and $F/N$ is free. Tensoring $\phi$ with $\Frac(R)$, we obtain
a surjection $F \otimes_R \Frac(R) \to M \otimes_R \Frac(R)$ of vector
spaces of dimensions $n$ and $r$. Thus the kernel of this map
has rank $n-r$, which implies that $N$ has rank $n-r$
and $F/N$ is free of rank $r$.

Now localizing at each prime $\idp$ of
$R$, we obtain a surjection $M_\idp \to (F/N)_{\idp}$ of free modules
of the same rank. By a standard result, this map is in fact a bijection.
Thus $M \to F/N$ is locally bijective, hence is bijective, and $M$ is free
as desired.
\end{proof}

The following lemma is a weak form of Galois descent for B\'ezout rings;
its key value is that it does not require that the
ring extension be finite.
\begin{lemma} \label{lem:gal}
Let $R_1/R_2$ be an extension of B\'ezout domains and $G$
a group of automorphisms of $R_1$ over $R_2$, with fixed ring $R_2$.
Assume that
every $G$-stable finitely
generated ideal of $R_1$ contains a nonzero element of $R_2$.
Let $M_2$ be a finite free module over $R_2$ and $N_1$ a saturated
$G$-stable submodule of
$M_1 = M_2 \otimes_{R_2} R_1$ stable under $G$. Then $N_1$ is equal to 
$N_2 \otimes_{R_2} R_1$ for a saturated
submodule (necessarily unique) $N_2$ of $M_2$.
\end{lemma}
\begin{proof}
We induct on $n = \rank M_2$, the case $n=0$ being trivial.
Let $\be_1, \dots, \be_n$ be a basis of $M_2$, and let
$P_1$ be the intersection of $N_1$ with the span of $\be_2, \dots, \be_n$;
since $N_1$ is saturated, $P_1$ is a direct summand of 
$\Span(\be_2, \dots, \be_n)$ by Lemma~\ref{lem:satspan} and hence
also of $M_1$.
By the induction hypothesis, $P_1 = P_2 \otimes_{R_2} R_1$ for a
saturated submodule $P_2$ of $M_2$ (necessarily a direct summand by
Lemma~\ref{lem:satspan}).
If $N_1 = P_1$, we are done. Otherwise,
$N_1/P_1$ is a $G$-stable, finitely generated ideal of $R_1$
(since $N_1$ is finitely generated by Lemma~\ref{lem:satspan}),
so contains a nonzero element $c$ of $R_2$.
Pick $\bv \in N_1$ reducing to $c$; that is,
$\bv - c \be_1 \in \Span(\be_2, \dots, \be_n)$.

Pick generators $\bw_1, \dots, \bw_m$ of
$P_2$; since $P_2$ is a direct summand of $\Span(\be_2, \dots, \be_n)$,
we can choose $\bx_1, \dots, \bx_{n-m-1}$ in $M_2$ so that
$\be_1, \bw_1, \dots, 
\bw_m, \bx_1, \dots, \bx_{n-m-1}$ is a basis of $M$. In this basis,
we may write $\bv = c \be_1 + \sum_i d_i \bw_i + \sum_i f_i \bx_i$,
where $c$ is the element of $R_2$ chosen above.
Put $\by = \bv - \sum_i d_i \bw_i$. For any $\tau \in G$, we have
$\by^\tau = c \be_1 + \sum_i f_i^\tau \bx_i$, so on one hand,
$\by^\tau - \by$ is a linear combination of $\bx_1,
\dots, \bx_{n-m-1}$. On the other hand, $\by^\tau - \by$ belongs to
$N_1$ and so is a linear combination of $\bw_1, \dots, \bw_m$. This
forces $\by^\tau - \by = 0$ for all $\tau \in G$; since $G$ has fixed
ring $R_2$, we conclude $\by$ is defined over $R_2$. Thus we may
take $N_2 = \Span(\by, \bw_1, \dots, \bw_m)$.
\end{proof}
Note that the hypothesis that every $G$-stable finitely
generated ideal of $R_1$ contains a nonzero element of $R_2$
is always satisfied if $G$ is finite: for any nonzero $r$
in the ideal, $\prod_{\tau \in G} r^\tau$ is nonzero and $G$-stable,
so belongs to $R_2$.

\subsection{$\sigma$-modules and $(\sigma, \nabla)$-modules}
\label{subsec:sigmod}

The basic object in the local study of $p$-adic differential equations is
a module with connection and Frobenius structure. In our approach, we
separate these two structures and study the Frobenius structure closely
before linking it with the connection. To this end, in this section
we introduce $\sigma$-modules and $(\sigma, \nabla)$-modules,
and outline some basic facts of what might be dubbed ``semilinear algebra''.
These foundations, in part, date back to Katz
\cite{bib:katz} and were expanded by de~Jong \cite{bib:dej1}.

For $R$ an integral domain in which $p \neq 0$
and $\sigma$ a ring endomorphism of $R$,
we define a \emph{$\sigma$-module} over $R$ to be a finite locally free
$R$-module $M$ equipped with an $R$-linear map $F: M \otimes_{R, \sigma} R
\to M$ that becomes an isomorphism over $R[\fp]$; the tensor product
notation indicates that $R$ is viewed as an $R$-module via $\sigma$.
Note that we will only use this definition when $R$ is a B\'ezout ring,
in which case every finite locally
free $R$-module is actually free by Proposition~\ref{prop:free}.
Then to specify $F$, it is equivalent to specify an additive,
$\sigma$-linear map from $M$ to $M$ that acts on any basis of $M$ by
a matrix invertible over $R[\fp]$. We abuse notation and refer to this map
as $F$ as well; since we will only use the $\sigma$-linear map in what
follows (with one exception: in proving Proposition~\ref{prop:tsuzuki}), 
there should not be any confusion induced by this.

Now suppose $R$ is one of $\GK, \GK[\fp], \GKcon,
\GKcon[\fp]$ or $\GKancon$ for $K = k((t))$.
Let $\Omega^1_R$ be the free module over $R$ generated by a single
symbol $du$, and let $d: R \to \Omega^1_R$
be the $\calO$-linear derivation given by the
formula
\[
d\left(\sum_i x_i u^i\right) = \sum_i i x_i u^{i-1}\,du.
\]
We define a $(\sigma, \nabla)$-module over $R$
to be a $\sigma$-module
$M$ plus a connection $\nabla: M \to M \otimes_R \Omega^1_{R}$ 
(i.e., an additive map satisfying the Leibniz rule
$\nabla(c\bv) = c\nabla(\bv) + \bv \otimes dc$ for $c \in R$
and $\bv \in M$) that makes the
following 
diagram commute:
\[
\xymatrix{
M \ar^-{\nabla}[r] \ar^{F}[d] & M \otimes \Omega^1_{R} \ar^{F \otimes
d\sigma}[d] \\
M  \ar^-{\nabla}[r] & M \otimes \Omega^1_{R}
}
\]
Warning: this definition is not the correct one in general. For 
larger rings $R$, one must include the condition that $\nabla$ is
integrable. That is, writing $\nabla_1$ for the induced map
$M \otimes_R \Omega^1_R \to M \otimes_R \wedge^2 \Omega^1_R$, we have
$\nabla_1 \circ \nabla = 0$. This condition is superfluous in our
context because $\Omega^1_R$ has rank one, so $\nabla_1$ is automatically
zero.

A \emph{morphism} of $\sigma$-modules or $(\sigma, \nabla)$-modules is a
homomorphism of the underlying $R$-modules compatible with the additional
structure in the obvious fashion. An \emph{isomorphism} 
of $\sigma$-modules or $(\sigma, \nabla)$-modules is a
morphism admitting an inverse; an \emph{isogeny} is a morphism that becomes
an isomorphism over $R[\fp]$.

Direct sums, tensor products, exterior powers,
and subobjects are defined in
the obvious fashion, as are duals if $p^{-1} \in R$;
quotients also make sense provided that the quotient
$R$-module is locally free. In particular, if $M_1 \subseteq M_2$ is
an inclusion of $\sigma$-modules, the saturation of $M_1$ in $M_2$
is also a $\sigma$-submodule of $M_1$; if $M_1$ itself is saturated,
the quotient $M_2/M_1$ is locally free and hence is a $\sigma$-module.

Given $\lambda$ fixed by $\sigma$, we define the \emph{twist} of a
$\sigma$-module $M$ by $\lambda$ as the $\sigma$-module with the same
underlying module but whose Frobenius has been multiplied by $\lambda$.

We say a $\sigma$-module $M$ is \emph{standard} if it is isogenous to
a $\sigma$-module
with a basis $\bv_1, \dots, \bv_n$ such that $F\bv_i = \bv_{i+1}$
for $i=1, \dots, n-1$ and $F\bv_n = \lambda \bv_1$ for some
$\lambda \in R$ fixed by $\sigma$. (The fact that $\lambda$ is fixed
by $\sigma$ is inserted for convenience only.) If $M$ is actually a
$(\sigma, \nabla)$-module, we say $M$ is standard as a $(\sigma, 
\nabla)$-module
if the same condition holds with the additional restriction that
$\nabla \bv_i = 0$ for $i=1, \dots, n$ (i.e., the $\bv_i$ are ``horizontal
sections'' for the connection).
If $\bv$ is a nonzero element of $M$ such that 
$F\bv = \lambda \bv$ for some $\lambda$, we say
$\bv$ is an \emph{eigenvector} of $M$ of \emph{eigenvalue} $\lambda$
and \emph{slope} $v_p(\lambda)$.

Warning: elsewhere in the literature, the
slope may be normalized differently, namely as $v_p(\lambda)/v_p(q)$.
(Recall that $q = p^f$.)
Since we will hold $q$ fixed, this normalization will not affect our
results.

From Lemma~\ref{lem:gal}, we have the following descent lemma for
$\sigma$-modules. (The condition on $G$-stable ideals is satisfied
because $R_1/R_2$ is an unramified extension of discrete valuation rings.)
\begin{cor} \label{cor:galois}
Let $R_1/R_2$ be an unramified
extension of discrete valuation rings unramified over $\calO$,
and let $\sigma$ be
a ring endomorphism of $R_1$ carrying $R_2$ into itself. 
Let
$\Gal^\sigma(R_1/R_2)$ be the group of automorphisms of $R_1$ over $R_2$
commuting with $\sigma$; assume that this group
has fixed ring $R_2$. Let $M_2$ be a $\sigma$-module over $R_2$ and $N_1$
a saturated $\sigma$-submodule of $M_1 = M_2 \otimes_{R_2} R_1$ stable under
$\Gal^\sigma(R_1/R_2)$.
Then $N_1 = N_2 \otimes_{R_2} R_1$ for some $\sigma$-submodule $N_2$
of $M_2$.
\end{cor}

\section{A few more rings}
\label{sec:aux}

In this chapter, we define a number of additional
auxiliary rings used in our study of $\sigma$-modules.
Again, we advise the reader to skim this chapter on first reading
and return to it as needed.

\subsection{Cohen rings}
\label{subsec:cohen}

We proceed to generalizing the 
constructions of Section~\ref{subsec:classical}
to valued fields. This cannot be accomplished using Witt
vectors because $k((t))$ and its finite extensions are not perfect.
To get around this, we fix once and for all a levelwise continuous
Frobenius lift $\sigma_0$ on
$\Gamma^{k((t))}$ carrying $\Gamma^{k((t))}_{\con}$ into itself; 
all of our constructions will
be made relative to the choice of $\sigma_0$.

Recall that a valued field $K$ is defined to be an algebraic extension
of $k((t))$ admitting subextensions
$k((t)) \subseteq L \subseteq M \subseteq N \subseteq K$ such that:
\begin{enumerate}
\item[(a)] $L = k^{1/p^m}((t))$ for some $m \in \{0,1, \dots, \infty\}$;
\item[(b)] $M = k_M((t))$ for some separable algebraic extension $k_M/k^{1/p^m}$;
\item[(c)] $N = M^{1/p^n}$ for some $n \in \{0,1, \dots, \infty\}$;
\item[(d)] $K$ is a separable totally ramified algebraic extension of $N$.
\end{enumerate}
We will associate to each valued field $K$ a complete discrete valuation ring $\Gamma^K$
unramified over $\calO$, equipped with a Frobenius lift $\sigma_0$ extending
the definition of $\sigma_0$ on $\Gamma^{k((t))}$. This assignment will be functorial
in $K$.

Let $\calC$ be the category of 
complete discrete valuation rings unramified over $\calO$,
in which morphisms are unramified morphisms of rings (i.e., morphisms
which induce isomorphisms of the value groups). If $R_0, R_1 \in \calC$
have residue fields $k_0, k_1$ and a homomorphism
$\phi: k_0 \to k_1$ is given, we say
the morphism $f: R_0 \to R_1$ is \emph{compatible} (with $\phi$) if the diagram
\[
\xymatrix{
R_0 \ar^f[r] \ar[d] & R_1 \ar[d] \\
k_0 \ar^\phi[r] & k_1
}
\]
commutes.

\begin{lemma} \label{lem:sep1}
Let $k_1/k_0$ be a finite separable extension of fields,
and take $R_0 \in \calC$ with residue field $k_0$.
Then there exists $R_1 \in \calC$ with residue field $k_1$
and a compatible morphism $R_0 \to R_1$.
\end{lemma}
\begin{proof}
By the primitive element theorem, there exists a monic separable polynomial
$\overline{P}(x)$ over $k_0$ and an isomorphism
$k_1 \cong k_0[x]/(\overline{P}(x))$. Choose a monic polynomial
$P(x)$ over $R_0$ lifting $\overline{P}(x)$ and set
$R_1 = R_0[x]/(P(x))$. Then the inclusion $R_0 \to R_0[x]$
induces the desired morphism $R_0 \to R_1$.
\end{proof}

\begin{lemma} \label{lem:sep2}
Let $k_0 \to k_1 \to k_2$ be homomorphisms of fields, with 
$k_1/k_0$ finite separable. For $i=0,1,2$, take $R_i \in \calC$
with residue field $k_i$. Let $f: R_0 \to R_1$
and $g: R_0 \to R_2$ be compatible morphisms. Then there exists
a unique compatible morphism $h: R_1 \to R_2$
such that $g = h \circ f$.
\end{lemma}
\begin{proof}
As in the previous proof, choose a monic separable polynomial
$\overline{P}(x)$ over $k_0$ and an isomorphism
$k_1 \cong k_0[x]/(\overline{P}(x))$. 
Let $y$ be the image of $x + (\overline{P}(x))$ in $k_1$,
and let $z$ be the image of $y$ in $k_2$.

Choose a monic polynomial
$P(x)$ over $R_0$ lifting $\overline{P}(x)$, and view $R_0$ as a subring
of $R_1$ and $R_2$ via $f$ and $g$, respectively. By Hensel's lemma,
there exist unique roots $\alpha$ and $\beta$ of $P(x)$ in
$R_1$ and $R_2$ reducing to $y$ and $z$, respectively, so
$h$ must satisfy $h(\alpha) = \beta$ if it exists.
Then $R_0[x]/(P(x)) \cong R_1$ by the map sending $x$ to $\alpha$
and $R_0[x]/(P(x)) \hookrightarrow R_2$ by the map sending $x$ to $\beta$,
so there exists a unique $h:
R_1 \to R_2$ such that $h(\alpha) = \beta$,
and this gives the desired morphism.
\end{proof}
\begin{cor} \label{cor:sep2}
If $k_1/k_0$ is finite Galois, and $R_i \in \calC$ has residue field $k_i$
for $i=0,1$, then for any compatible morphism $f: R_0 \to R_1$,
the group of $f$-equivariant automorphisms of $R_1$ is isomorphic to
$\Gal(k_1/k_0)$.
\end{cor}
\begin{proof}
Apply Lemma~\ref{lem:sep2} with $k_0 \to k_1$ the given embedding and
$k_1 \to k_1$ an element of $\Gal(k_1/k_0)$; the resulting $h$ is
the corresponding automorphism.
\end{proof}
\begin{cor} \label{cor:sep3}
If $k_1/k_0$ is finite separable, $R_i \in \calC$ has residue field $k_i$
for $i=0,1$, and $f: R_0 \to R_1$ is a compatible morphism, then
any compatible endomorphism of $R_0$ admits a unique $f$-equivariant
extension to $R_1$.
\end{cor}
\begin{proof}
If $\phi: R_0 \to R_0$ is the given endomorphism, apply
Lemma~\ref{lem:sep2} with $g = f \circ \phi$.
\end{proof}

For $m$ a nonnegative integer, let $\calO_m$ be a copy of $\calO$.
Then the assignment $k^{1/p^m} \leadsto \calO_m$ is functorial
via the morphism $\sigma_0^i$ compatible with $k^{1/p^m} \to
k^{1/p^{m+i}}$; thus we can define $\calO_\infty$ as the completed
direct
limit of the $\calO_m$.
For any finite separable extension $k_M$ of $k^{1/p^m}$, choose
$\calO_M$ in $\calC$ according to Lemma~\ref{lem:sep1},
to obtain a compatible morphism $\calO_m \to \calO_M$; note
that $\calO_M$ is unique up to \emph{canonical} isomorphism by
Lemma~\ref{lem:sep2}. Moreover, this assignment is functorial in $k_M$
(again by Lemma~\ref{lem:sep2}),
so again we may pass to infinite extensions by taking the completed
direct limit.

Now suppose $K$ is a valued field finite over $k((t))$
and $L,M,k_M, N,n$ are as in the definition of valued fields; recall that since
$n<\infty$, $L,M,k_M, N$ are uniquely determined by $K$.
Define $\calO_M$ associated to $k_M$ as above, define 
$\Gamma^M$ as the ring of
power series $\sum_{i \in \ZZ} a_i u^i$, with $a_i \in \calO_M$,
such that $|a_i| \to 0$ as $i \to -\infty$, and identify
$\Gamma^M / \pi \Gamma^M$ with $M = k_M((t))$
via the map $\sum_i a_i u^i \mapsto \sum_i \overline{a_i} t^i$.
Define $\Gamma^N$ as a copy of $\Gamma^M$, but with $\Gamma^M$ embedded
via $\sigma_0^n$ (which makes sense since $n<\infty$), and identify
the residue field of $\Gamma^N$ with $N$ compatibly.
Define $\Gamma^K$ as a copy of $\Gamma^N$ with its
residue field identified with $K$
via some continuous $k_M^{1/p^n}$-algebra
isomorphism $K \cong N$ (which exists because
both fields are power series fields over 
$k_M^{1/p^n}$ by the Cohen structure theorem).
Once this choice is made, there is a unique
(necessarily levelwise continuous) morphism $\Gamma^N \to \Gamma^K$ compatible
with the embedding $N \hookrightarrow K$.
The assignments of $\Gamma^M, \Gamma^N, \Gamma^K$ are functorial,
again by Lemma~\ref{lem:sep2},
so again we may extend the definition to infinite $K$ by completion.

Note that if $K/k((t))$ is finite, then $\Gamma^K$ is equipped
with a levelwise topology, and the embeddings provided by functoriality
are levelwise continuous. Moreover, $\sigma_0$ extends uniquely to
each $\Gamma^K$, and the functorial morphisms are $\sigma_0$-equivariant.

If $k$ and $K$ are perfect and
$\calO = C(k) = W(k)$,
then $\Gamma^K$ is canonically isomorphic to the
Witt ring $W(K)$. Under that isomorphism,
$\sigma_0$ corresponds to the Witt vector Frobenius, which sends
each Teichm\"uller lift to its $p$-th power.
For general $\calO$, we have $\Gamma^K \cong
W(K) \otimes_{W(k)} \calO$.

We will often fix a field $K$ (typically $k((t))$ itself) and
write $\Gamma$ instead of $\Gamma^K$. In this case, we will frequently
refer to $\Gamma^L$ for various canonical extensions $L$ of $K$, such as
the separable closure $K^{\sep}$,
the perfect closure $K^{\perf}$, and
the algebraic closure $K^{\alg}$.
In all of these cases, we will drop the $K$
from the notation where it is understood, writing $\Gamma^{\perf}$
for $\Gamma^{K^{\perf}}$ and so forth.

\subsection{Overconvergent rings} \label{sec:overcon}

Let $K$ be a valued field.
Let $v_K$ denote the valuation on $K$ extending the valuation on
$k((t))$, normalized so that $v_K(t) = 1$.
Again, let $q = p^f$, and put $\sigma = \sigma_0^f$ on
$\Gamma^K$.
We define a subring $\Gamma^K_{\con}$
of $\Gamma^K$ of ``overconvergent'' elements; the construction will not
look quite like the construction of $\Gamma^{k((t))}_{\con}$ from
Section~\ref{subsec:classical}, so we must
check that the two are consistent.

First assume $K$ is perfect. For
$x \in \Gamma^K[\fp]$, write $x = \sum_{i=m}^\infty \pi^i [\overline{x_i}]$,
 where
$m v_p(\pi) = v_p(x)$, each $\overline{x_i}$ belongs to $K$ and the
brackets denote Teichm\"uller lifts. For $n$ in the value group
of $\calO$, we define the ``partial valuations''
\[
v_n(x) = \min_{v_p(\pi^i) \leq n} \{ v_K(\overline{x_i}) \}.
\]
These partial valuations obey two rules
analogous to those for their na\"\i ve counterparts, plus a third
that has no analogue:
\begin{align*}
v_n(x+y) &\geq \min\{v_n(x), v_n(y) \} \\
v_n(xy) &\geq \min_m \{v_m(x) + v_{n-m}(y)\} \\
v_n(x^\sigma) &= q v_n(x).
\end{align*}
Again, equality holds in the first two lines if the minimum is achieved
exactly once.

For each $r>0$, let $\Gamma^K_r$ denote the subring
of $x \in \Gamma^K$ such that
$\lim_{n \to \infty} (r v_n(x) + n) = \infty$. On $\Gamma^K_r[\fp]
\setminus \{0\}$, we define the function
\[
w_r(x) = \min_n \{ r v_n(x) + n \};
\]
then $w_r$ is a nonarchimedean valuation by the same argument 
as for $w_r^{\naive}$ given in Section~\ref{subsec:classical}.
Define $\Gamma^K_{\con} = \cup_{r>0} \Gamma^K_r$. 

The rings $\Gamma^K_r$ will be quite useful,
but one must handle them with some caution, for the following reasons:
\begin{enumerate}
\item[(a)]
The map $\sigma: \Gamma^K \to \Gamma^K$ sends $\Gamma^K_{\con}$ into
itself, but does not send $\Gamma^K_r$ into itself; rather, it sends
$\Gamma^K_r$ into $\Gamma^K_{r/q}$. 
\item[(b)]
The ring $\GKcon$ is a discrete valuation ring, but the rings
$\Gamma^K_r$ are not.
\item[(c)]
The ring $\Gamma^K_r$ is complete for $w_r$, but not for
the $p$-adic valuation.
\end{enumerate}

For $K$ arbitrary, we want to define $\Gamma^K_{\con}$ as $\Galgcon \cap \GK$. 
This intersection is indeed a discrete valuation ring
(so again its fraction field is obtained by adjoining $\fp$), but it is
not clear that its residue field is all of $K$. Indeed, it is \emph{a
priori} possible that the intersection is no larger than $\calO$ itself!
In fact, this pathology does not occur, as we will see below.

To make that definition, we must also check that
$\Galgcon \cap \Gamma^{k((t))}$ coincides with
the ring $\Gamma^{k((t))}_{\con}$ defined earlier. This is obvious in a
special case: if $\sigma_0(u) = u^p$, then $u$ is a Teichm\"uller lift in
$\Galgcon$, and in this case one can check that the
partial valuations and na\"\i ve partial valuations coincide. In general
they do not coincide, but in a sense they are not too far apart.
The relationship might be likened to that between the na\"\i ve and canonical
heights on an abelian variety over a number field.

Put $z = u^\sigma/u^q - 1$. By the original definition of $\sigma$
on $\Gamma^{k((t))}$, $v_p(z) > 0$ and $z \in \Gamma^{k((t))}_{\con}$.
That means we can find $r>0$ such that
$r v^{\naive}_n(z) + n > 0$ for all $n$; for all $s \leq r$,
we then have $w^{\naive}_s(u^\sigma/u^q) = 0$.

\begin{lemma} \label{lem:naivesigma}
Choose $r>0$ such that
$r v^{\naive}_n(z) + n > 0$ for all $n$.
For $x = \sum_i x_i u^i$ in $\Gamma^{k((t))}_{r, \naive}$, if $0<s\leq qr$ and
$w^{\naive}_s(x) \geq c$, then
$w^{\naive}_{s/q}(x^\sigma) \geq c$.
\end{lemma}
\begin{proof}
We have
\begin{align*}
w^{\naive}_{s/q}(x_i^\sigma (u^i)^\sigma)
&= w^{\naive}_{s/q}(x_i u^{qi} (u^\sigma/u^q)^i) \\
&= w^{\naive}_{s/q}(x_i u^{qi}) + w^{\naive}_{s/q}((u^\sigma/u^q)^i) \\
&= w^{\naive}_s(x_i u^i)
\end{align*}
since $w^{\naive}_{s/q}(u^\sigma/u^q) = 0$ whenever $s/q \leq r$.

Given that $w^{\naive}_s(x) \geq c$, it follows that $w^{\naive}_{s/q}(x_i u^i)
\geq c$ for each $i$, and by the above argument, that $w^{\naive}_{s/q}
(x_i^\sigma (u^i)^\sigma) \geq c$. We conclude that
$w^{\naive}_{s/q}(x^\sigma) \geq c$, as desired.
\end{proof}

\begin{lemma} \label{lem:naiveineq}
Choose $r>0$ such that
$r v^{\naive}_n(z) + n > 0$ for all $n$.
For $x \in \Gamma^{k((t))}_{r,\naive}$, if
$sv^{\naive}_j(x) + j \geq c$ for all $j \leq n$, then
$sv_j(x) + j \geq c$ for all $j \leq n$.
\end{lemma}
\begin{proof}
Note that $v_0 = v_0^{\naive}$, so the desired result holds for $n=0$;
we prove the general result by induction on $n$. Suppose, as the induction
hypothesis, that if
$sv^{\naive}_j(x) + j \geq c$ for all $j < n$, then
$sv_j(x) + j \geq c$ for all $j < n$. Before deducing the desired
result, we first study the special case $x = u$ in detail (but using
the induction hypothesis in full generality).

Choose $i$ large enough that
\[
v_p([t] - (u^{\sigma^{-i}})^{q^i}) > n.
\]
Then
\begin{align*}
v_n(u) &\geq \min\{v_n([t]), v_{n}(u - [t])\} \\
&= \min\{1, v_n(u - (u^{\sigma^{-i}})^{q^i})\}.
\end{align*}
Applying $\sigma^i$ yields
\[
q^i v_n(u) \geq \min\{q^i, v_n(u^{\sigma^i} - u^{q^i})\}.
\]
Since $u \in \Gamma^{k((t))}_{r,\naive}$ and $w^{\naive}_r(u) = r$
trivially, we may apply Lemma~\ref{lem:naivesigma}
to $u, u^\sigma, \dots, u^{\sigma^{i-1}}$ in succession to
obtain
\[
w^{\naive}_{r/q^i}(u^{\sigma^i}) \geq r.
\]
Since $w^{\naive}_{r/q^i}(u^{q^i}) = r$, we conclude that
$w^{\naive}_{r/q^i}(u^{\sigma^i} - u^{q^i}) \geq r$.

Let $y = (u^{\sigma^i}-u^{q^i})/\pi$. Then for $j \leq n-v_p(\pi)$,
\begin{align*}
(r/q^i) v^{\naive}_{j}(y) + j
&= (r/q^i) v^{\naive}_{j+v_p(\pi)}(y \pi) + j+v_p(\pi) - v_p(\pi)\\
&\geq w^{\naive}_{r/q^i}(y \pi) - v_p(\pi) \\
&\geq r - v_p(\pi).
\end{align*}
By the induction hypothesis, we conclude that $(r/q^i) v_{n-v_p(\pi)}(y) + 
n - v_p(\pi) \geq r - v_p(\pi)$, and so $(r/q^i) v_n(y \pi) + n \geq r$. 
From above, we have
\begin{align*}
q^i v_n(u) &\geq \min\{q^i, v_n(u^{\sigma^i} - u^{q^i})\} \\
&\geq \min\{q^i, q^i - q^i n/r \} \\
&= q^i - q^i n/r.
\end{align*}
Thus $r v_n(u) + n \geq r$.
Since $v_n(u) \leq 1$, we also have
$s v_n(u) + n \geq s$ for $s \leq r$;
 that is, the desired conclusion holds for the
special case $x=u$. By the multiplication rule for partial valuations,
we also have $s v_n(u^i) + n \geq s i$ for all $i$.

With the case $x=u$ in hand, we now prove the desired conclusion for 
general $x$. We are given $s v^{\naive}_j(x) + j \geq c$ for $j \leq n$;
by the induction hypothesis, all that we must
prove is that $sv_n(x) + n \geq c$.

The assumption $sv^{\naive}_j(x) + j \geq c$ implies that
$s v^{\naive}_j(x_i u^i) + j \geq c$ for all
$j \leq n$, which is to say, if $v_p(x_i) \leq n$ then $s i + v_p(x_i) \geq c$.
For $j = v_p(x_i)$, we have
\begin{align*}
s v_n(x_i u^i) + n &= s v_{n-j}(u^i) + n-j + j\\
& \geq si + j \\
&\geq c.
\end{align*}
We conclude that $s v_n(x) + n \geq c$,
completing the induction.
\end{proof}
We next refine the previous result as follows.
\begin{lemma} \label{lem:compnaive}
Choose $r>0$ such that
$r v^{\naive}_n(z) + n > 0$ for all $n$.
If $x \in \Gamma^{k((t))}_{r, \naive}$, then for any $s \leq r$,
$\min_{j \leq n} \{ sv^{\naive}_j(x) + j \} =
\min_{j \leq n} \{ s v_j(x) + j \}$ for all $n$. In particular,
$w^{\naive}_s(x) = w_s(x)$.
\end{lemma}
That is, the na\"\i ve valuations $w^{\naive}_s$ are not so
simple-minded after all; as long as $s$ is not too large, they agree with
the more canonically defined $w_s$.
\begin{proof}
Lemma~\ref{lem:naiveineq} asserts that
$\min_{j \leq n} \{ sv_j(x) + j \} \geq
\min_{j \leq n} \{ s v^{\naive}_j(x) + j \}$, so it remains to prove
the reverse inequality, which we do by induction on $n$. If
$\min_{j \leq n} \{ s v^{\naive}_j(x) + j \}$ is achieved by some
$j < n$, then by the induction hypothesis,
\begin{align*}
\min_{j \leq n} \{ s v^{\naive}_j(x) + j\}
&= \min_{j \leq n - v_p(\pi)} \{ s v^{\naive}_j(x) + j\} \\
&\geq \min_{j \leq n - v_p(\pi)} \{ s v_j(x) + j \} \\
&\geq \min_{j \leq n} \{ s v_j(x) + j\}.
\end{align*}
Suppose then that $\min_{j \leq n} \{ s v^{\naive}_j(x) + j \}$ is 
achieved only for $j=n$.
Put $x = \sum x_i u^i$; 
by definition, $v^{\naive}_n(x)$ is the smallest integer $i$ with
$v_p(x_i) \leq n$.
In fact, we must have $v_p(x_i) = n$, or else we have
$s v^{\naive}_j(x) + j < s v^{\naive}_n(x) + n$ for $j = v_p(x_i)$.
Therefore $v_n(x_i u^i) = v_n^{\naive}(x_i u^i) = i$.

For $j < n$, $s v^{\naive}_j(x - x_i u^i) + j
= s v^{\naive}_j(x) + j > si + n$.
On the other hand, $v^{\naive}_n(x) = v^{\naive}_n(x_i u^i) = i$
and $v^{\naive}_n(x - x_i u^i) > i$. Thus for all $j \leq n$,
\[
s v^{\naive}_j(x - x_i u^i) + j > si + n;
\]
by Lemma~\ref{lem:naiveineq}, $s v_n(x-x_i u^i) + n > si + n$
and so $v_n(x - x_i u^i) > i = v_n(x_i u^i)$. Therefore
$v_n(x) = v_n(x_i u^i) = i$, so
\[
\min_{j \leq n} \{ s v_j(x) + j\} 
\leq s v_n(x) + n = si + n = \min_{j \leq n} \{ s v_j^{\naive}(x) + j\},
\]
yielding the desired inequality.
\end{proof}
\begin{cor} \label{cor:naiveeq}
We have $\Galgcon \cap \Gamma^{k((t))} = \Gamma^{k((t))}_{\con}$.
\end{cor}

We now define $\GKcon = \Galgcon \cap \GK$, and Corollary~\ref{cor:naiveeq}
assures that this definition is consistent with our prior definition
for $K = k((t))$.
To show that $\Galgcon \cap \GK$ is ``large'' for any
$K$, we need one more lemma, which will end up generalizing a
 standard fact about
$\Gamma^{k((t))}_{\con}$.
\begin{lemma} \label{lem:intclos}
For any valued field $K$, $\GKcon$ is Henselian.
\end{lemma}
\begin{proof}
By a lemma of Nagata \cite[43.2]{bib:nag}, it suffices to show that if
$P(x)=x^d + a_1 x^{d-1} + \cdots + a_d$
is a polynomial over $\GKcon$ such that $a_1 \not\equiv 0 \pmod{\pi}$
and $a_i \equiv 0 \pmod{\pi}$ for $i>1$, then $P(x)$ has a root
$y$ in $\GKcon$ such that $y \equiv -a_1 \pmod{\pi}$.
By replacing $P(x)$ by $P(-x/a_1)$, we may reduce to the
case $a_1 = -1$; by Hensel's lemma, $P$ has a root $y$ in $\GK$
congruent to $1$ modulo $\pi$, and $P'(y) \equiv dy^{d-1} - (d-1)y^{d-2}
\equiv 1 \pmod{\pi}$.

Choose a constant $c>0$ such that $v_n(a_i) \geq -cn$
for all $n$, and define the sequence 
$\{y_j\}_{j=0}^\infty$ by the Newton iteration, putting $y_0 = 1$
and $y_{j+1} = y_j - P(y_j)/P'(y_j)$. Then $\{y_j\}$ converges $\pi$-adically
to $y$; we now show by induction on $j$ that
$v_n(y_j) \geq -cn$ for all $n$
and all $j$. Namely, this is obvious for $y_0$, and given $v_n(y_j) \geq -cn$
for all $n$, it follows that $v_n(P(y_j)) \geq -cn$, $v_n(P'(y_j)) \geq -cn$,
and $v_n(1/P'(y_j)) \geq -cn$ (the last because $v_0(P'(y_j)) = 0$). These
together imply $v_n(y_{j+1}) \geq -cn$ for all $n$, completing the induction.
We conclude that $y \in \GKcon$ and $\GKcon$ is Henselian, as desired.
\end{proof}

We can now prove the following. 
\begin{prop}
For any valued field $K$, $\GKcon$ has
residue field $K$.
\end{prop}
\begin{proof}
We have already shown this for $K= k((t))$ by Corollary~\ref{cor:naiveeq}.
If $K/k((t))$ is finite, then $K$ uniquely determines $L,M,N$ as in the
definition of valued fields. Now $M = k_M((t))$ for some finite extension
$k_M$ of $k$, so Corollary~\ref{cor:naiveeq} also implies that 
$\Gamma^M_{\con}$ 
has residue field $M$. Also, $N = M^{1/p^n}$ for some integer $n$,
so for any $\overline{x} \in M$, we can find $y \in \Gamma^M_{\con}$
which lifts $\overline{x}^{p^n}$, and then $y^{\sigma^{-n}} \in 
\Gamma^N_{\con}$ lifts $\overline{x}$.

Choose a monic polynomial $P$ over $\Gamma^{N}_{\con}$ lifting
a monic separable polynomial $\overline{P}$ for which $K \cong 
N[x]/(\overline{P}(x))$ (again, possible by the primitive element
theorem). The reduction is a separable polynomial, so by
Hensel's lemma $P$ has a root $y$ in $\GK$, and $\GK \cong 
\Gamma^{N}[y]/(P(y))$. But since $\Galgcon$ is Henselian and
$P$ has coefficients in $\Galgcon$, $y \in \Galgcon$. Thus the residue
field of $\GKcon$ contains $N$ and $y$, and hence is all of $K$.

This concludes the proof for $K$ finite over $k((t))$. A general valued
field $K$ is the union of its finite valued subfields $K_1$,
and $\GKcon$ contains (but does not equal) the
direct limit of the $\Gamma^{K_1}_{\con}$. Thus its residue field contains
the union of the $K_1$, and hence is equal to $K$.
\end{proof}
If $L/K$ is a finite extension of valued fields, then $\GLcon/\GKcon$
is an unramified extension of discrete valuation rings, and the corresponding
residue field extension is $L/K$, so is finite. Thus $\GLcon$ is 
integral over $\GKcon$.

\subsection{Analytic rings: generalizing the Robba ring}

In this section, we generalize the construction of the
Robba ring. Besides the classical case where $K$ is a finite extension
of $k((t))$, we will be especially interested in the case $K = k((t))^{\alg}$,
which will give a sort of ``maximal unramified extension'' of the standard Robba ring.

\begin{prop} \label{prop:uniform}
Suppose the valued field $K$ is either
\begin{itemize}
\item[(a)] finite over $k((t))$ or
\item[(b)] perfect.
\end{itemize}
Then there exists $r>0$ such that $\GK_r = \Galg_r \cap \GK$ has
units congruent to every nonzero element of $K$.
\end{prop}
\begin{proof}
For (a), let $u$ be a lift to $\GKcon$ of a uniformizer 
$\overline{u}$ of $K$,
and choose $r>0$ so that $u$ is a unit in $\GK_r$. Let $\calO'$
be the integral closure of $\calO$ in $\GK$; its residue field
is the integral closure $k'$ of $k$ in $K$.

For any $c_i \in \calO'$, the series $1 + \sum_{i=1}^\infty c_i u^i$
converges with respect to $w_r$ (hence levelwise)
to a unit of $\GK_r$, because we can formally invert the
series and the result also converges with respect to $w_r$.
Any nonzero element of $K$ can be written as a nonzero element of $k'$ times
a power of $\overline{u}$ times a series in $\overline{u}$ with
leading term 1, thus can be lifted
as an invertible element of $\calO'$ times a power of $u$ times a series
of the form $1 + \sum_{i=1}^\infty c_i u^i$. The result is invertible
in $\GK_r$, as desired.

For (b), we can choose any $r>0$, since every Teichm\"uller lift 
belongs to $\GK_r$.
\end{proof}
Note that the conclusion of the proposition need not hold for other
valued fields. For example, it fails for $K = k((t))^{\sep}$ if
$\sigma_0(u) = u^p$ for some $u \in \Gamma^{k((t))}_{\con}$
lifting $t$:
define a sequence $\{y_i\}_{i=1}^\infty$
of elements of $K$ by setting $y_i$ to be a root of
$y_{i}^p - y_{i} = u^{-i}$. Then it can be shown that
$y_{i}$ has a lift in $\GK_r$ only if $r < \frac{1}{i} (p/(p-1))^2$,
so there is no way to choose $r$ uniformly.

For the rest of this section, we assume that the hypotheses of
Proposition~\ref{prop:uniform} ace satisfied.
Recall that for $0 < s \leq r$, we have defined
the valuation $w_s$ on $\Gamma^K_r[\fp]$ by
\[
w_s(x) = \min_{n} \{ n + s v_n(x) \},
\]
the minimum taken as $n$ runs over the value group of $\calO$.
We define a corresponding norm $|\cdot|_s$ by $|x|_s = p^{-w_s(x)}$.

While $\Gamma^K_r$ is complete under $|\cdot|_r$, $\Gamma^K_r[\fp]$ is
not, so we can attempt to complete it. In fact,
we can define a Fr\'echet topology on $\Gamma^K_r[\fp]$ using
the $w_s$ for $0 < s \leq r$, and define $\Gamma^K_{\an, r}$
as the Fr\'echet completion of $\Gamma^K_r[\fp]$.
That is, $\Gamma^K_{\an,r}$ consists of equivalence classes of sequences
of elements of $\Gamma^K_r[\fp]$ which are simultaneously
Cauchy for all of the norms $|\cdot|_s$.

Set $\Gamma^K_{\an,\con} = \cup_{r>0} \Gamma^K_{\an,r}$.
Echoing a warning from the previous section, we note that
$\Gamma^K_{\an,\con}$ admits an action of $\sigma$, but
each $\Gamma^K_{\an,r}$ is mapped not into itself, but into
$\Gamma^K_{\an,r/q}$. More precisely, we have
$w_{r/q}(x^\sigma) = w_r(x)$ for all $x \in \Gamma^K_{\an,r}$.

In case $K = k((t))$, we defined another ring called $\GKancon$
in Section~\ref{subsec:classical}. Fortunately, these rings coincide:
for $r$ sufficiently small, by Corollary~\ref{cor:naiveeq}
we have $\GK_{r} = \GK_{r,\naive}$ and so $\GK_{\an,r} =
\GK_{\an,r,\naive}$ by Proposition~\ref{prop:freccomp}.

Since $\Gamma^K_{\an,\con}$ is defined from $\Gamma^K_{\con}$ by a
canonical completion process, it inherits as much functoriality as is
possible
given the restricted applicability of Proposition~\ref{prop:uniform}.
For example, if $L/K$ is a finite extension, then $\GLancon$ is an
integral extension of $\GKancon$; in fact, one has a canonical
identification of $\GLancon$ with $\GLcon \otimes_{\GKcon} \GKancon$,
which in case $L/K$ is Galois gives an action of $\Gal(L/K)$ on $\GLancon$ with fixed ring
$\GKancon$. 
Likewise, if $K$ is perfect, then the union
$\cup_L \Gamma^L_{\an,r}$ running over all finite subextensions $L$ of $K$
is dense in $\GK_{\an,r}$ for each $r>0$, so $\cup_L \GLancon$ is
dense in $\GKancon$.

We can extend the functions $v_n$ to $\Gamma^K_{\an,r}$ by continuity: if
$x_i \to x$ in the Fr\'echet topology, then $v_n(x_i)$ either stabilizes
at some finite value or tends to $+\infty$ as $i \to \infty$, and we may
put $v_n(x) = \lim_{i \to \infty} v_n(x_i)$. Likewise, we can extend
the functions $w_s$ to $\GK_{\an,r}$ by continuity, and again one has the
formula
\[
w_s(x) = \min_n \{n + s v_n(x)\},
\]
as $n$ runs over the value group of $\calO$. One also has
\[
\lim_{n\to \pm \infty} (n + s v_n(x)) = \infty
\]
for any $0<s<r$. For $n \to -\infty$, this follows from the corresponding
limiting statement for $s=r$. For $n \to \infty$, note that if the
limit did not tend to infinity, $x$ could not be written as a limit
under $|\cdot|_s$ of elements of $\GK_{r}[\fp]$.

It is not so easy to prove anything about the ring $\GK_{\an,\con}$ just
from the above definition, since it is inconvenient to even write down
elements of this ring. To this end, we isolate a special class of elements,
which we call semiunits, and use them as building blocks to represent
more general ring elements.

We define a \emph{semiunit} of $\GK_r$ (resp.\ of $\GK_{\an,r}$)
as an element $u$ of $\GK_r$ (resp.\ of $\GK_{\an,r}$)
which is either zero, or which satisfies the following conditions:
\begin{enumerate}
\item[(a)] $v_n(u) = \infty$ for $n < 0$;
\item[(b)] $v_0(u) < \infty$;
\item[(c)] $r v_n(u) + n > r v_0(u)$ for $n > 0$.
\end{enumerate}
In particular, if $u \in \GK_r$, then $u$ is a semiunit if either
$u=0$ or $u$ is a unit in $\GK_r$, hence the terminology.
In particular, under the condition of Proposition~\ref{prop:uniform},
every element of $K$ lifts to a semiunit in $\GK_r$.
Note that if $u$ is a
semiunit in $\GK_{\an,r}$, it is also a semiunit in $\GK_{\an,s}$
for any $0<s<r$. Also beware that if $K/k((t))$ is infinite,
a semiunit in $\GK_{\an,r}$
need not belong to $\GK_r$ even though $v_p(u) \geq 0$. (If $R$ is
the subring of $x \in \GK_{\an,r}$ with $v_p(x) \geq 0$,
then $R/\pi R$ is isomorphic to the completion of $K$ with respect to
$v_K$.)

If $K$ is perfect, we define a \emph{strong semiunit}
of $\GK_{r}$ (resp.\ of $\GK_{\an,r}$)
 as an element $u$ of $\GK_r$ (resp.\ of $\GK_{\an,r}$)
which is either zero,
or satisfies the following conditions:
\begin{enumerate}
\item[(a)] $v_n(u) = \infty$ for $n < 0$;
\item[(b)] $v_0(u) < \infty$;
\item[(c)] $v_n(u) = v_0(u)$ for $n > 0$.
\end{enumerate}
Every Teichm\"uller lift is a strong semiunit, so every
element of $K$ lifts to a strong semiunit in $\GK_r$.

Let $\{u_i\}_{i=-\infty}^\infty$ be a doubly infinite sequence of
semiunits in $\GK_r$ (resp.\ in $\GK_{\an,r}$).
Then we say $\{u_i\}$ is a \emph{semiunit decomposition}
of $x$ in $\GK_r$ (resp.\ in $\GK_{\an,r}$)
if $w_r(u_i \pi^i) \leq w_r(u_j \pi^j)$ whenever
$i>j$ and $u_i, u_j\neq 0$, and if
$\sum_{i=-M}^N u_i \pi^i$ converges to $x$ in the Fr\'echet
topology as $M, N \to \infty$. We express this more succinctly by
saying that $\sum u_i \pi^i$ is a semiunit
decomposition of $x$. Analogously,
if $K$ is perfect and the $u_i$ are strong semiunits, we say
$\sum u_i \pi^i$ is a \emph{strong semiunit
decomposition} of $x$ if $v_0(u_i) < v_0(u_j)$ whenever $i > j$
and $u_i, u_j \neq 0$, and if $\sum_{i=-M}^N u_i \pi^i$ converges
to $x$ in the Fr\'echet topology as $M, N \to \infty$.

If $\sum u_i \pi^i$ is a semiunit decomposition of $x \in \GK_{\an,r}$,
then for each $i$ such that $u_i \neq 0$, we may set $n = i v_p(\pi)$
and obtain $r v_n(x) + n = r v_n(u_i \pi^i) + n$, that is,
$v_n(x) = v_n(u_i \pi^i)$.
Since $r v_n(x) + n \to \infty$ as $n \to \infty$ for any $x \in \GK_{\an,r}$,
we must then have $u_i = 0$ for
$i$ sufficiently large. 
There is no analogous phenomenon for strong
semiunit decompositions, however: for each $i$ such that $u_i \neq 0$,
we set $n = i v_p(\pi)$ and obtain $v_n(x) = v_n(u_i \pi^i)$, but 
$v_n(x)$ may continue to decrease forever as $n \to \infty$, so
the $u_i$ need not eventually vanish.

\begin{lemma} \label{lem:semidec1}
Each element $x$ of $\GK_r$ admits a semiunit decomposition.
If $K$ is perfect,
each element $x$ of $\GK_r$ admits a strong semiunit decomposition.
\end{lemma}
\begin{proof}
Without loss of generality (by dividing by a suitable power of $\pi$),
we may reduce to the case where $x \not\equiv 0
\pmod{\pi}$.
We define a sequence of semiunits $\{y_i\}_{i=0}^\infty$ such that
$x \equiv \sum_{i=0}^j y_i \pi^i \pmod{\pi^{j+1}}$, as follows.
Let $y_0$ be a semiunit congruent to $x$ modulo $\pi$.
Given $y_0, \dots,  y_j$, let $y_{j+1}$ be a semiunit congruent to
$(x - \sum_{i=0}^j y_i \pi^i)/\pi^{j+1}$ modulo $\pi$.

The sum $\sum_{i=0}^\infty y_i \pi^i$ now converges to $x$, but 
we do not have the necessary comparison between $w_r(y_i \pi^i)$
and $w_r(y_j \pi^j)$, so we must revise the decomposition.
We say $i$ is a \emph{corner} if $w_r(y_i \pi^i) = 
\min_{j\leq i} \{w_r(y_j \pi^j)\}$. We now set $u_i = 0$ if $i$ is
not a corner; if $i$ is a corner, let $l$ be the next largest corner
(or $\infty$ if there is none), and put
$u_i = \sum_{j=i}^{l-1} y_j \pi^{j-i}$. By the definition of 
a corner, $w_r(y_j \pi^{j-i}) > w_r(y_i)$ for $i < j < l$, 
so $u_i$ is a semiunit.
Moreover, if $i$ and $j$ are corners and $i>j$, then
$w_r(u_i \pi^i) = w_r(y_i \pi^i) \leq w_r(y_j \pi^j) = w_r(u_j \pi^j)$;
and the sum $\sum_{i=0}^\infty u_i \pi^i$ is merely the sum
$\sum_{i=0}^\infty y_i \pi^i$ with the terms regrouped, so
it still converges to $x$. Thus $\sum_{i=0}^\infty u_i \pi^i$
is a semiunit decomposition of $x$.

If $K$ is perfect, we perform the revision slightly differently.
We say $i$ is a corner if $v_0(y_i) < v_0(y_j)$ for all $j<i$.
Again, we set $u_i=0$ if $i$ is not a corner, and if $i$ is a corner
and $l$ is the next largest corner, we set $u_i = \sum_{j=i}^{l-i}
y_j \pi^{j-i}$. Clearly $u_i$ is a strong semiunit for each $i$,
and the sum $\sum_{i=0}^\infty u_i \pi^i$ converges
to $x$. If $i>j$ are corners, then $v_0(u_i) = v_0(y_i) < v_0(y_j)
= v_0(u_j)$. Thus $\sum_{i=0}^\infty u_i \pi^i$ is a strong
semiunit decomposition of $x$.
\end{proof}

\begin{prop} \label{prop:semidec2}
Every element of $\GK_{\an,r}$ admits a semiunit decomposition.
\end{prop}
\begin{proof}
For $x \in \GK_{\an,r}$, let $\sum_{l=0}^\infty x_l$ be a series
of elements of $\GK_r[\fp]$ that converges under $|\cdot|_r$ to $x$,
such that $w_r(x_l) < w_r(x_{l+1})$. (For example,
choose $x_0$ such that $w_r(x - x_0) > w_r(x)$, then choose
$x_1$ such that $w_r(x - x_0 - x_1) > w_r(x -x_0)$, and so forth.)

For $l=0,1,\dots$ and $i \in \ZZ$, we define elements
$y_{il}$ of $\GK_r[\fp]$ recursively in $l$, such that for any
$l$, only finitely many of the $y_{il}$ are nonzero,
as follows.
Apply Lemma~\ref{lem:semidec1} (after multiplying by a suitable
power of $\pi$)
to produce a semiunit decomposition of $x_0 + \cdots + x_l
- \sum_{j<l} \sum_i y_{ij} \pi^i$.
For each of the finitely
many terms  
$u_i \pi^i$ of this decomposition with $u_i \neq 0$ and
$w_r(u_i \pi^i) < w_r(x_{l+1})$, put $y_{il} = u_i$;
for all other $i$, put $y_{il} = 0$. Then
\[
w_r\left(x_0 + \cdots + x_l - \sum_{j\leq l} \sum_i y_{ij}\pi^i\right) 
\geq w_r(x_{l+1}).
\]
In particular, the doubly infinite sum $\sum_l \sum_i y_{il} \pi^i$
converges under $|\cdot|_r$ to $x$.
If we set
$z_i = \sum_l y_{il}$, the series $\sum_i z_i \pi^i$
converges under $|\cdot|_r$ to $x$.

Note that $w_r(x_l) \leq w_r(y_{il} \pi^i) < w_r(x_{l+1})$ whenever
$y_{il} \neq 0$. Thus for any fixed $i$, the values of $w_r(y_{il}
\pi^i)$,
taken over all $l$ such that $y_{il} \neq 0$, form a strictly
increasing sequence. If $j$ is the first such index, we
then have $w_r(y_{ij} \pi^i) < w_r(\sum_{l>j} y_{il} \pi^i)$,
and so $z_i$ is a semiunit.

Define $u_i$ to be zero if $w_r(z_i \pi^i) > w_r(z_j \pi^j)$ for some
$j < i$; otherwise, let $l$ be the smallest integer greater than
$i$ such that
$w_r(z_l \pi^l) \leq w_r(z_i \pi^i)$ (or $\infty$ if none exists),
and put $u_i = \sum_{j=i}^{l-1}
z_j \pi^{j-i}$. Then the series $\sum_i u_i \pi^i$ also converges
under $|\cdot|_r$ to $x$, and if $u_i \neq 0$, then
$w_r(u_i \pi^i) = w_r(z_i \pi^i)$. It follows that
$w_r(u_i \pi^i) \leq w_r(u_j \pi^j)$ whenever $i > j$
and $u_i, u_j \neq 0$. This in turn implies that if $u_i \neq 0$
and $n = v_p(\pi^i)$, then $v_n(u_i \pi^i) = v_n(x)$.

We finally check that $\sum_i u_i \pi^i$ converges under
$|\cdot|_s$ for $0<s<r$.
The fact that $s v_{n}(x) + n \to \infty$ as $n \to \pm \infty$
implies that $s v_{v_p(\pi^i)}(u_i \pi^{i}) + v_p(\pi^i)
\to \infty$ as $i \to \pm \infty$. Since $u_i$ is a
semiunit, $w_s(u_i \pi^i) = s v_{v_p(\pi^i)}(u_i \pi^i) + v_p(\pi^i)$,
so $w_s(u_i \pi^i) \to \infty$ as $i \to \pm \infty$.
Thus the sum $\sum_i u_i \pi^i$ converges under $|\cdot|_s$
for $0 < s< r$, and the limit must equal $x$ because the sum
converges to $x$ under $|\cdot|_r$.
Therefore $\sum_i u_i \pi^i$ is a semiunit decomposition,
as desired.
\end{proof}

\begin{prop} \label{prop:semidec3}
If $K$ is perfect,
every element of $\GK_{\an,r}$ admits a strong semiunit decomposition.
\end{prop}
\begin{proof}
As in the previous proof,
for $x \in \GK_{\an,r}$, let $\sum_{l=0}^\infty x_l$ be a series
of elements of $\GK_r[\fp]$ that converges under $|\cdot|_r$ to $x$,
such that $w_r(x_l) < w_r(x_{l+1})$.

For $l=0,1,\dots$ and $i \in \ZZ$, we define elements
$y_{il}$ of $\GK_r[\fp]$ recursively in $l$, such that for any
$l$, only finitely many of the $y_{il}$ are nonzero,
as follows.
Apply Lemma~\ref{lem:semidec1} to produce a
strong semiunit decomposition of $x_0 + \cdots + x_l
- \sum_{j<l} \sum_i y_{ij} \pi^i$. For each of the finitely
many terms  
$u_i \pi^i$ of this decomposition with $u_i \neq 0$ and
$w_r(u_i \pi^i) < w_r(x_{l+1})$, put $y_{il} = u_i$;
for all other $i$, put $y_{il} = 0$. Then
\[
w_r\left(x_0 + \cdots + x_l - \sum_{j\leq l} \sum_i y_{ij}\pi^i\right) 
\geq w_r(x_{l+1}).
\]
In particular, the doubly infinite sum $\sum_l \sum_i y_{il} \pi^i$
converges under $|\cdot|_r$ to $x$.
If we set
$z_i = \sum_l y_{il}$, and the series $\sum_i z_i \pi^i$
converges under $|\cdot|_r$ to $x$.

Note that $w_r(x_l) \leq w_r(y_{il} \pi^i) < w_r(x_{l+1})$ whenever
$y_{il} \neq 0$. Thus for any fixed $i$, the values of $v_0(y_{il})$,
taken over all $l$ such that $y_{il} \neq 0$, form a strictly
increasing sequence. If $j$ is the first such index, we
then have $v_0(y_{ij}) < v_0(\sum_{l>j} y_{il})$, and so $z_i$
is a strong semiunit.

Define $u_i$ to be zero if $v_0(z_i) \geq v_0(z_j)$ for some
$j < i$; otherwise, let $l$ be the smallest integer such that
$v_0(z_l) < v_0(z_i)$ (or $\infty$ if none exists), 
and put $u_i = \sum_{j=i}^{l-1}
z_j \pi^{j-i}$. Then the series $\sum_i u_i \pi^i$ also converges
under $|\cdot|_r$ to $x$, and if $u_i \neq 0$, then
$v_0(u_i) = v_0(z_i)$. It follows that
$v_0(u_i) < v_0(u_j)$ whenever $i > j$
and $u_i, u_j \neq 0$. This in turn implies that if $u_i \neq 0$
and $n = v_p(\pi^i)$, then $v_n(u_i \pi^i) = v_n(x)$.

We finally check that $\sum_i u_i \pi^i$ converges under
$|\cdot|_s$ for $0<s<r$, by the same argument as in the
previous proof. Namely,
the fact that $s v_{n}(x) + n \to \infty$ as $n \to \pm \infty$
implies that $s v_{v_p(\pi^i)}(u_i \pi^{i}) + v_p(\pi^i)
\to \infty$ as $i \to \pm \infty$. Since $u_i$ is a strong
semiunit, $w_s(u_i \pi^i) = s v_{v_p(\pi^i)}(u_i \pi^i) + v_p(\pi^i)$,
so $w_s(u_i \pi^i) \to \infty$ as $i \to \pm \infty$.
Thus the sum $\sum_i u_i \pi^i$ converges under $|\cdot|_s$
for $0 < s< r$, and the limit must equal $x$ because the sum
converges to $x$ under $|\cdot|_r$.
Therefore $\sum_i u_i \pi^i$ is a strong semiunit decomposition,
as desired.
\end{proof}

Although (strong) semiunit decompositions are not unique, in a certain
sense the ``leading terms'' are unique. To make sense of this 
remark, we first need a ``leading coefficient map'' for $K$.
\begin{lemma} \label{lem:leadcoef}
For $K$ a valued field, there exists a homomorphism
$\lambda: K^* \to (k^{\alg})^*$ such that $\lambda(c) = c$ for all
$c \in k^{\alg} \cap K$ and $\lambda(x) = 1$ if $v_K(x-1) >0$.
\end{lemma}
For instance, if $K = k((t))$, we could take $\lambda(x)$ to be the
leading coefficient of $x$.

\begin{proof}
There is no loss of generality in enlarging $K$, so we may assume
$K = k((t))^{\alg}$. Define $t_0 = t$, and for $i>0$, let $t_i$ be
an $i$-th root of $t_{i-1}$. With this choice, for any $d \in \QQ$
we can define $t^d$ as $t_i^{i!d}$ for any $i \geq d$; the expression
does not depend on $i$.

Now for each $x \in K^*$, there exists a unique $c \in k^{\alg}$
such that
\[
v_K\left( \frac{x}{ct^{v_K(x)}} - 1 \right) > 0;
\]
set $\lambda(x) = c$.
\end{proof}

Choose a map $\lambda$ as in Lemma~\ref{lem:leadcoef};
we define the leading terms map
$L_r: \GK_{\an,r} \to \cup_{n=1}^\infty k^{\alg}[t^{1/n},t^{-1/n}]$ 
as follows.
For $x \in \GK_{\an,r}$ nonzero, find a finite sum
$y = \sum_j u_j \pi^j$ such that each $u_j$ is a semiunit,
$w_r(u_j \pi^j) = w_r(x)$ for all $j$ such that $u_j \neq 0$,
and $w_r(x - y) > w_r(x)$. Then put $L_r(x) =
\sum_j \lambda(\overline{u_j}) t^{v_0(u_j)}$;
this definition does not depend on the choice of $y$.
Moreover, the leading terms map is multiplicative, that is,
$L_r(xy) = L_r(x) L_r(y)$.

We define the \emph{upper degree} and \emph{lower degree} of a nonzero element
of $\cup_{n=1}^\infty k^{\alg}[t^{1/n},t^{-1/n}]$ as the largest and smallest
powers of $t$, respectively, occurring in the element; we define the
\emph{length} of an element as the upper degree minus the lower degree. We
extend all of these definitions to $\GK_{\an,r}$ through the map $L_r$.

Warning: if $K$ is not finite over $k((t))$, then the subring
of $x \in \GKancon$ with $v_n(x) = \infty$ for
$n < 0$ is a complete discrete valuation ring containing $\GKcon$, but it
is actually much bigger than $\GKcon$. In fact, its residue field
is the completion of $K$ with respect to the valuation $v_K$.

As noted earlier, a theorem of Lazard asserts that 
$\GKancon$ is a B\'ezout ring (every finitely generated ideal
is principal) for $K = k((t))$; the same is true for $K$ a finite extension of
$k((t))$, since $K \cong k'((t'))$ for some uniformizer $t'$ and some
field $k'$. We will generalize the B\'ezout property to $\GKancon$ for
$K/k((t))$ infinite in 
Section~\ref{subsec:anbezout2}; for now, we deduce from Lemma~\ref{lem:gal}
the following descent lemma
for $\sigma$-modules. (The condition on $G$-stable
ideals is satisfied because $G = \Gal(L/K)$ here is finite.)
\begin{cor} \label{cor:galois2}
Let $L/K$ be a finite Galois extension of valued fields finite over $k((t))$.
Let $M$ be a $\sigma$-module over $\GKancon$ and $N$
a saturated $\sigma$-submodule of $M \otimes_{\GKancon} \GLancon$ 
stable under $\Gal(L/K)$.
Then $N$ is equal to $P \otimes_{\GKancon} \GLancon$ 
for some saturated $\sigma$-submodule $P$ of $M$.
\end{cor}

\subsection{Some $\sigma$-equations}

We record here the behavior of some simple equations involving $\sigma$. 
For starters, we have the following variant of Hensel's lemma.
\begin{prop} \label{prop:hensel}
Let $R$ be a complete discrete valuation ring, unramified over 
$\calO$, with separably closed residue field, and let
$\sigma$ be a $q$-power Frobenius lift. For
$c_0, \dots, c_n \in R$ with $c_0$ not divisible by $\pi$
and $x \in R$, define $f(x) = c_0 x + c_1 x^\sigma + \cdots
+ c_n x^{\sigma^n}$. Then for any $x, y \in R$ for which $f(x)
\equiv y \pmod{\pi}$, there exists $z \in R$ congruent to
$x$ modulo $\pi$ for which $f(z) = y$. Moreover, if $R$ has
algebraically closed residue field, then the same holds if any
of $c_0, \dots, c_n$ is not divisible by $\pi$.
\end{prop}
\begin{proof}
Define a sequence $\{z_l\}_{l=1}^\infty$ of elements of $R$ such
that $z_1 = x$, $z_{l+1} \equiv z_l \pmod{\pi^l}$ and
$f(z_l) \equiv y \pmod{\pi^{l}}$; then the limit $z$ of the $z_l$ will
have the desired property. Given $z_l$, put $a_l
= (y - f(z_l))/\pi^l$, and choose $b_l \in R$ such that
\[
c_0 b_l + c_1 b_l^q (\pi^\sigma/\pi)^l + \cdots + c_n b_l^{q^n} (\pi^{\sigma^n}/\pi)^l \equiv a_l
\pmod{\pi};
\]
this is possible because either $R$ has algebraically closed residue field,
or $c_0 \neq 0$ and  the polynomial at left must be separable.
Put $z_{l+1} = z_l + \pi^l b_l$; then $f(z_{l+1}) \equiv f(z_l) + f(\pi^l b_l)
\equiv y \pmod{\pi^{l+1}}$,
as desired.
\end{proof}

We next consider similar equations over some other rings. The
following result will be vastly generalized by
Proposition~\ref{prop:descfilt} later.
\begin{prop} \label{prop:rank1}
Suppose $x \in \Galgcon$ (resp.\ $x \in \Galgancon$ with $v_n(x) = \infty$
for $n < 0$) is not congruent to $0$ modulo $\pi$. Then
there exists a nonzero $y \in \Galgcon$ (resp.\ $y \in \Galgancon$
with $v_n(y) = \infty$ for $n < 0$) such that $y^\sigma = xy$.
\end{prop}
\begin{proof}
Put $R = \Galgcon$ (resp.\ let $R$ be the subring of
$x \in \Galgancon$ with $v_n(x) = \infty$ for $n<0$) and let $S$
be the completion of $R$.
By Proposition~\ref{prop:hensel}, we can find nonzero $y \in S$
such that $y^\sigma = xy$; we need to show that $y \in R$. Choose $r>0$ and
$c \in \RR$ such that $r v_n(x) + n \geq c$ for all $n$.
We then show that $r (q-1) v_n(y) + n \geq c$ by
induction on $n$. We have
\[
q v_n(y) = v_n(y^\sigma) \geq \min_{m \leq n} \{ v_m(x) + v_{n-m}(y) \}.
\]
If the minimum is achieved for $m = 0$ (which includes the base case
$n=0$), then $(q-1) v_n(y) \geq v_0(x)$,
so $r (q-1) v_n(y) + n \geq r(q-1) v_0(x) + n \geq c$.
If the minimum is achieved for some $m>0$, then by the
induction hypothesis
\begin{align*}
r (q-1) v_n(y) &\geq \frac{r(q-1)}{q} v_m(x) + \frac{r(q-1)}{q} 
v_{n-m}(y) \\
&\geq \frac{(q-1)(c - m)}{q} + \frac{(c-n+m)}{q} \\
&\geq \frac{(q-1)(c-n)}{q} + \frac{(c-n)}{q} \geq c-n,
\end{align*}
so $r(q-1) v_n(y) + n \geq c$. Thus the induction goes through,
and demonstrates that $y \in R$, as desired.
\end{proof}


Finally, we consider a class of equations involving the analytic rings.
We suppress $K$ from all
superscripts for convenience, writing $\Gcon$ for $\GKcon$ and so forth.
\begin{prop} \label{prop:ansigeq}
Let $K$ be a valued field (satisfying the condition of 
Proposition~\ref{prop:uniform} in case $\Gancon$ is referenced).
\begin{enumerate}
\item[(a)]
Assume $K$ is separably closed (resp.\ algebraically closed).
For $\lambda \in \calO$ a unit
and $x \in \Gcon$ (resp.\ $x \in \Gancon$), 
there exists $y \in \Gcon$ (resp.\ $y \in \Gancon$)
such that $y^\sigma - \lambda y = x$. Moreover, if
$x \in \Gcon[\fp]$, then any such $y$ belongs to $\Gcon[\fp]$.
\item[(b)]
Assume $K$ is perfect.
For $\lambda \in \calO$ not a unit
and $x \in \Gcon$ (resp.\ $x \in \Gancon$),
there exists $y \in \Gcon$ (resp.\ $y \in \Gancon$)
such that $y^\sigma - \lambda y = x$.
Moreover, we can take $y$ nonzero in $\Gancon$ even if
$x=0$.
\item[(c)]
For $\lambda \in \calO$ not a unit
and $x \in \Gancon$, there is at most one $y \in \Gancon$
such that $\lambda y^\sigma - y = x$, and if $x \in \Gcon$,
then $y \in \Gcon$ as well.
\item[(d)]
For $\lambda \in \calO$ not a unit
and $x \in \Gancon$ such that $v_n(x) \geq 0$ for all $n$, there exists
$y \in \Gancon$ such that $\lambda y^\sigma - y = x$.
\end{enumerate}
\end{prop}
\begin{proof}
\begin{enumerate}
\item[(a)]
If $x \in \Gcon$, then Proposition~\ref{prop:hensel} implies that
there exists $y \in \Gamma$ such that $y^\sigma - \lambda y = x$.
To see that in fact $y \in \Gcon$, note that if $v_n(y) \leq 0$, the fact that
\[
q v_n(y) = v_n(y^\sigma) = v_n(\lambda y + x) \geq \min\{v_n(x), v_n(y)\}
\]
implies that $q v_n(y) \geq v_n(x)$; while if $v_n(y) > 0$, the fact that
\[
v_n(y) = v_n(\lambda y) = v_n(y^\sigma - x) \geq \min\{q v_n(y), v_n(x) \}
\]
implies that $v_n(y) \geq v_n(x)$, which also implies $q v_n(y) \geq v_n(x)$.
Hence $y \in \Gcon$ and $w_{qr}(y) \geq w_r(x)$.

For $x \in \Gancon$ (with $K$ algebraically closed),
choose $r>0$ such that $x \in \Gamma_{\an,r}$,
and let $x = \sum_{i=-\infty}^\infty u_i \pi^i$ be a strong semiunit
decomposition. As above, there exists $y_i \in \Gamma_{\an,qr}$ with
$y_i^\sigma - \lambda(\pi/\pi^\sigma)^i y_i = u_i (\pi/\pi^\sigma)^i$
such that $v_n(y_i) = \infty$ for $n<0$
and $w_{qr}(y_i) \geq w_r(u_i)$.
This implies that $\sum_{i=-\infty}^\infty y_i \pi^i$ converges with respect
to $|\cdot|_s$ for $0<s\leq r$; let $y$ be its Fr\'echet
limit.
Then
\begin{align*}
y^\sigma - \lambda y &=
\sum_i y_i^\sigma (\pi^i)^\sigma -  \lambda y_i \pi^i \\
&= \sum_i \lambda y_i \pi^i + u_i \pi^i - \lambda y_i \pi^i \\
&= \sum_i u_i \pi^i = x,
\end{align*}
so $y$ is the desired solution.

To verify the last assertion, suppose $x \in \Gcon[\fp]$ and $y \in \Gancon$
satisfy $y^\sigma - \lambda y = x$. By what we have shown above, there also
exists $z \in \Gcon[\fp]$ such that $z^\sigma - \lambda z = x$, so
$(y-z)^\sigma = \lambda(y-z)$.
This equation yields $q v_n(y-z) = v_n(y-z)$ for all $n$, so
$v_n(y-z) = 0$ or $\infty$ for all $n$. 
We cannot have $v_n(y-z) = 0$
for all $n$, so there is a smallest such $n$; we may assume $n=0$
without loss of generality.
Let $\calO'$ be the
completed integral closure of $\calO$ in $\Gcon$. Then every solution $w$
of $w^\sigma = \lambda w$ in $\Gancon$ with $v_n(w)=\infty$ for $n<0$
is congruent to some element of $\calO'$ modulo $\pi$. In particular,
we can find $c_0, c_1, \dots \in \calO'$ such that
$\sum_{j=0}^l c_j \pi^j \equiv y-z \pmod{\pi^{l+1}}$, since
once $c_0, \dots, c_l$ have been computed, we can
take $w = (y-z)\pi^{-l-1} - \sum_{j=0}^l c_j \pi^{j-l-1}$, and there
must be some $c_{l+1} \in \calO'$ congruent to $w$ modulo $\pi$.
Thus
$y-z \in \calO' \subseteq \Gcon[\fp]$, so $y \in \Gcon[\fp]$.

\item[(b)]
If $x \in \Gcon$, then the series
\[
\sum_{i=0}^\infty \lambda^{\sigma^{-1}}\cdots \lambda^{\sigma^{-i}}
x^{\sigma^{-i-1}}
\]
converges $\pi$-adically to an element $y \in \Gamma$
satisfying
\begin{align*}
y^\sigma - \lambda y &= 
\sum_{i=0}^\infty 
\lambda \lambda^{\sigma^{-1}} \cdots \lambda^{\sigma^{-i+1}}
x^{\sigma^{-i}}
- \sum_{i=0}^\infty
\lambda \lambda^{\sigma^{-1}} \cdots \lambda^{\sigma^{-i}}
x^{\sigma^{-i-1}} \\
&= 
\sum_{i=0}^\infty 
\lambda \lambda^{\sigma^{-1}} \cdots \lambda^{\sigma^{-i+1}}
x^{\sigma^{-i}}
- \sum_{i=1}^\infty
\lambda \lambda^{\sigma^{-1}} \cdots \lambda^{\sigma^{-i+1}}
x^{\sigma^{-i}} \\
&= x.
\end{align*}
To see that in fact $y \in \Gcon$, choose $r>0$ and $c \leq 0$ such that
$w_r(x) \geq c$, that is, $r v_n(x) + n \geq c$ for all $n \geq 0$.
If $v_n(x) \leq 0$, then
$r v_n(x^{\sigma^{-i}}) + n = (r/q^i) v_n(x) + n \geq r v_n(x) + n
\geq c$; if $v_n(x) \geq 0$, then $r v_n(x^{\sigma^{-i}}) + n 
\geq 0 \geq c$. In any case,
we have $w_r(x^{\sigma^{-i}}) \geq c$ for all $i$.
Since $w_r(\lambda^{\sigma^{-i}}) = w_r(\lambda) > 0$ for all $i$,
we conclude that the series defining $y$ converges under
$|\cdot|_r$, and so its limit $y$ in $\Gamma$ must actually lie
in $\Gcon$.

Suppose now that $x \in \Gancon$;
by Proposition~\ref{prop:semidec3}, there exists a strong semiunit
decomposition $x = \sum_n \pi^n u_n$ of $x$.
Let $N$ be the largest value of $n$ for which $v_0(u_n) \geq 0$,
and put
\[
x_- = \sum_{n=-\infty}^N \pi^n u_n,
\qquad
x_+ = \sum_{n=N+1}^\infty \pi^n u_n.
\]
As above, we can construct $y_{+} \in \Gcon[\fp]$ so that
$y_+^\sigma - \lambda y_+ = x_+$.
As for $x_-$, let $m$ be the greatest integer less than or equal
to $N$ for which $u_m \neq 0$. For any fixed $r$,
$w_r(x_-^{\sigma^i}) = w_r((u_m \pi^m)^{\sigma^i})$ for $i$
sufficiently large.
The series
\[
- \sum_{i=0}^\infty (\lambda \lambda^\sigma \cdots
\lambda^{\sigma^i})^{-1} x_-^{\sigma^i}
\]
then converges under $|\cdot|_r$, since
\begin{align*}
w_r((\lambda \lambda^\sigma \cdots
\lambda^{\sigma^i})^{-1} x_-^{\sigma^i}) &=
-(i+1) w_r(\lambda) + w_r(x_-^{\sigma^i}) \\
&=
-(i+1) w_r(\lambda) + r q^i v_0(u_m) + m v_p(\pi)
\end{align*}
tends to
infinity with $i$. Since this holds for every $r$, the series
converges in $\Gancon$ to a limit $y_-$, which satisfies
\begin{align*}
y_-^\sigma - \lambda y_- &=
- \sum_{i=0}^\infty (\lambda^\sigma \cdots
\lambda^{\sigma^{i+1}})^{-1} x_-^{\sigma^{i+1}}
+
\sum_{i=0}^\infty (\lambda^\sigma \cdots
\lambda^{\sigma^i})^{-1} x_-^{\sigma^i}
\\
&=
- \sum_{i=1}^\infty (\lambda^\sigma \cdots
\lambda^{\sigma^{i}})^{-1} x_-^{\sigma^{i}}
+
\sum_{i=0}^\infty (\lambda^\sigma \cdots
\lambda^{\sigma^i})^{-1} x_-^{\sigma^i}
\\
&= x_-.
\end{align*}
We conclude that $y = y_+ + y_-$ satisfies $y^\sigma - \lambda y = x$.

To prove the final assertion, let $u$ be any strong semiunit with
$v_0(u) > 0$, and set
\[
y = \sum_{i=0}^\infty \lambda^{\sigma^{-1}} \cdots \lambda^{\sigma^{-i}}
u^{\sigma^{-i-1}}
+ \sum_{i=0}^\infty (\lambda \lambda^{\sigma}\cdots 
\lambda^{\sigma^{i}})^{-1}
u^{\sigma^i};
\]
then the above arguments show that both series converge and $y^\sigma - \lambda
y = u - u = 0$.

\item[(c)]
We prove the second assertion first. Namely, assume $x \in \Gcon$
and $y \in \Gancon$ satisfy $\lambda y^\sigma - y = x$; we show that
$y \in \Gcon$. First suppose $0 < v_n(y) < \infty$ for some $n < 0$.
Then
\[
v_n(y) = v_n(y+x) = v_n(\lambda y^\sigma)
\geq v_n(y^\sigma) = q v_n(y),
\]
contradiction. 
Thus $v_n(y)$ is either nonpositive or $\infty$
for all $n<0$. We cannot have $v_n(y) < 0$ for all $n$, since for
some $r>0$ we have $r v_n(y) + n \to \infty$ as $n \to -\infty$.
Thus $v_n(y) = \infty$ for some $y$. (Beware: this is not enough
\emph{a priori} to imply that $y \in \Gcon[\fp]$ if $K$ is infinite
over $k((t))$.) Choose $n$ minimal such that $v_n(y) < \infty$. If
$n < 0$, then $v_n(y) = v_n(y+x) = v_n(\lambda y^\sigma) = \infty$,
contradiction. Thus $n \geq 0$. We can now show that
$y$ is congruent modulo $\pi^i$ to an element of $\Gcon$, by
induction on $i$. The base case $i=0$ is vacuous; given
$y \equiv y_i \pmod{\pi^i}$ for $y_i \in \Gcon$, we have
\[
y = -x + \lambda y^\sigma \equiv -x + \lambda y_i^\sigma \pmod{\pi^{i+1}}.
\]
Thus the induction follows. Since $y$ is the $\pi$-adic limit of 
elements of $\Gcon$, we conclude $y \in \Gcon$.

For the first assertion, suppose $x \in \Gancon$ and $y_1, y_2 \in
\Gancon$ satisfy $\lambda y_i^\sigma - y_i = x$ for $i=1, 2$. Then
$\lambda (y_1-y_2)^\sigma - (y_1-y_2) = 0$; by the previous paragraph,
this implies $y_1-y_2 \in \Gcon$. But then $v_p(y_1-y_2)
= v_p(\lambda) + v_p((y_1-y_2)^\sigma)$, a contradiction unless
$y_1-y_2=0$.

\item[(d)]
Since $v_n(x) \geq 0$ for all $n$, we have $v_n(x^{\sigma^i}) = q^i v_n(x) \geq v_n(x)$
for all nonnegative integers $i$. Thus
$w_s(x^{\sigma^i}) \geq w_s(x)$ for all $s$, so the series
\[
y = - \sum_{i=0}^\infty \lambda \lambda^{\sigma}\cdots \lambda^{\sigma^{i-1}}
x^{\sigma^{i}}
\]
converges with
respect to each of the norms $|\cdot|_s$,
and
\begin{align*}
\lambda y^\sigma - y &=
- \sum_{i=0}^\infty \lambda \lambda^{\sigma}\cdots \lambda^{\sigma^{i}}
x^{\sigma^{i+1}}
+ \sum_{i=0}^\infty \lambda \lambda^{\sigma}\cdots \lambda^{\sigma^{i-1}}
x^{\sigma^{i}} \\
&=
- \sum_{i=1}^\infty \lambda \lambda^{\sigma}\cdots \lambda^{\sigma^{i-1}}
x^{\sigma^{i}}
+ \sum_{i=0}^\infty \lambda \lambda^{\sigma}\cdots \lambda^{\sigma^{i-1}}
x^{\sigma^{i}} \\
&= x,
\end{align*}
so $y$ is the desired solution.
\end{enumerate}
\end{proof}

\subsection{Factorizations over analytic rings}
\label{subsec:anbezout}

We assume that the valued field
$K$ satisfies the conditions of Proposition~\ref{prop:uniform},
so that the ring $\Gancon = \GKancon$ is defined.
As noted earlier, $\Gancon$ is not Noetherian
even for $K = k((t))$, but in this case Lazard
\cite{bib:laz} proved that $\Gancon$ is a
B\'ezout ring, that is, a ring in which every finitely generated ideal
is principal. 
In this section and the next, we generalize Lazard's result 
as follows.
\begin{theorem} \label{thm:bezout}
Suppose the conclusion of
Proposition~\ref{prop:uniform} is satisfied for the valued field $K$
and the positive number $r$.
Then every finitely generated ideal in $\Gamma_{\an,r} = 
\GK_{\an,r}$ is principal.
In particular, every finitely generated ideal in $\Gancon$ is
principal.
\end{theorem}
Our approach resembles that of Lazard,
with ``pure elements'' standing in for the divisors in his theory.
The approach requires a number of auxiliary results on
factorizations of elements of $\Gancon$; for the most
part (specifically, excepting Section~\ref{subsec:approx}),
only Theorem~\ref{thm:bezout} will be used in the sequel, not
the auxiliary results.

\begin{figure}[ht] 
\begin{center}
\setlength{\unitlength}{0.00087489in}
\begingroup\makeatletter\ifx\SetFigFont\undefined
\def\x#1#2#3#4#5#6#7\relax{\def\x{#1#2#3#4#5#6}}%
\expandafter\x\fmtname xxxxxx\relax \def\y{splain}%
\ifx\x\y   
\gdef\SetFigFont#1#2#3{%
  \ifnum #1<17\tiny\else \ifnum #1<20\small\else
  \ifnum #1<24\normalsize\else \ifnum #1<29\large\else
  \ifnum #1<34\Large\else \ifnum #1<41\LARGE\else
     \huge\fi\fi\fi\fi\fi\fi
  \csname #3\endcsname}%
\else
\gdef\SetFigFont#1#2#3{\begingroup
  \count@#1\relax \ifnum 25<\count@\count@25\fi
  \def\x{\endgroup\@setsize\SetFigFont{#2pt}}%
  \expandafter\x
    \csname \romannumeral\the\count@ pt\expandafter\endcsname
    \csname @\romannumeral\the\count@ pt\endcsname
  \csname #3\endcsname}%
\fi
\fi\endgroup
{\renewcommand{\dashlinestretch}{30}
\begin{picture}(3342,2321)(0,-10)
\put(237,1871){\blacken\ellipse{128}{128}}
\put(237,1871){\ellipse{128}{128}}
\put(687,971){\blacken\ellipse{128}{128}}
\put(687,971){\ellipse{128}{128}}
\put(1137,476){\blacken\ellipse{128}{128}}
\put(1137,476){\ellipse{128}{128}}
\put(552,1646){\blacken\ellipse{128}{128}}
\put(552,1646){\ellipse{128}{128}}
\put(1362,1286){\blacken\ellipse{128}{128}}
\put(1362,1286){\ellipse{128}{128}}
\put(2082,701){\blacken\ellipse{128}{128}}
\put(2082,701){\ellipse{128}{128}}
\put(2082,296){\blacken\ellipse{128}{128}}
\put(2082,296){\ellipse{128}{128}}
\put(2937,71){\blacken\ellipse{128}{128}}
\put(2937,71){\ellipse{128}{128}}
\path(237,1871)(687,971)(1137,476)(2937,71)
\path(12,971)(3162,971)
\blacken\path(3042.000,941.000)(3162.000,971.000)(3042.000,1001.000)(3042.000,941.000)
\path(1587,71)(1587,2006)
\blacken\path(1617.000,1886.000)(1587.000,2006.000)(1557.000,1886.000)(1617.000,1886.000)
\put(1542,2141){\makebox(0,0)[lb]{\smash{{{\SetFigFont{12}{14.4}{rm}$n$}}}}}
\put(3342,926){\makebox(0,0)[lb]{\smash{{{\SetFigFont{12}{14.4}{rm}$v_n$}}}}}
\end{picture}
}
\end{center}
\caption{An example of a Newton polygon}
\label{fig:newton}
\end{figure}
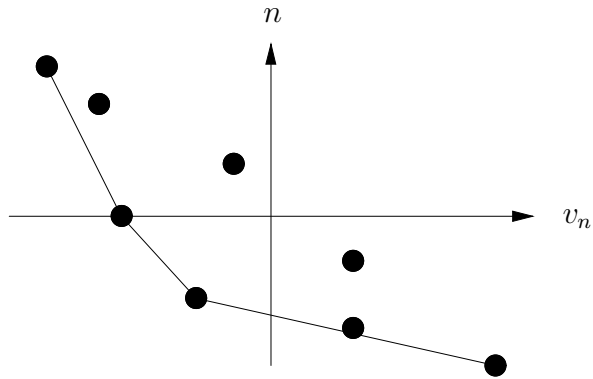

For $x \in \Gamma_{\an,r}$ nonzero, define the \emph{Newton
polygon} of $x$ as the lower convex hull of the set of points $(v_n(x),
n)$, minus
any segments of slopes less than $-r$ on the left end and/or
any segments of nonnegative slope on the right end of the polygon;
see Figure~\ref{fig:newton} for an example.
Define the \emph{slopes} of $x$ as the negatives of the
slopes of the Newton polygon of $x$. (The negation is to ensure that the
slopes of $x$ are positive.) Also define the \emph{multiplicity} of a slope
$s \in (0, r]$
of $x$ as the positive difference in $y$-coordinates between the endpoints of the
segment of the Newton polygon of slope $-s$, or 0 if there is no
such segment.
If $x$ has only one slope $s$, we say $x$ is \emph{pure} (of slope $s$).
(Beware: this notion of slope differs from the slope of an eigenvector
of a $\sigma$-module introduced in Section~\ref{subsec:sigmod}, and
the Newton polygon here does not correspond to either the generic or
special Newton polygons we define later.)

\begin{lemma}
The multiplicity of $s$ as a slope of $x$ is equal to 
$s$ times
the length (upper degree minus lower degree)
of $L_s(x)$, where $L_s$ is the leading terms map in
$\Gamma_{\an,s}$.
\end{lemma}
\begin{proof}
Let $\sum_i u_i \pi^i$ be a semiunit decomposition of $x$.
Let $S$ be the set of $l$ which achieve $\min_l \{w_s(u_l \pi^l)\}$, and let $i$ and $j$
be the smallest and largest elements of $S$;
then $L_s(x) = \sum_{l \in S} \lambda(\overline{u_l}) t^{v_0(u_l)}$
and the length of $L_s(x)$ is equal to $v_0(u_i) - v_0(u_j)$.

We now show that the endpoints of the segment of the Newton polygon
of $x$ of slope $-s$ are $(v_0(u_i), v_p(\pi^i))$ and 
$(v_0(u_j), v_p(\pi^j))$. First of all, for $n = v_p(\pi^i)$,
we have $s v_n(x) + n = s v_0(u_i) + i v_p(\pi) = w_s(u_i \pi^i)$;
likewise for $n = v_p(\pi^j)$. Next, we note that
$w_s(x) \geq \min_l \{w_s(u_l \pi^l)\} = w_s(u_i \pi^i)$.
Thus for any $n$, $s v_n(x) + n \geq w_s(u_i \pi^i)$; this
means that
the line through $(v_0(u_i), v_p(\pi^i))$ and $(v_0(u_j), v_p(\pi^j))$
is a lower supporting line for the set of points $(v_n(x), n)$.
Finally, note that for $n < v_p(\pi^i)$,
\begin{align*}
s v_n(x) + n &\geq \min_{l < i} \{ s v_n(\pi^l u_l) + n \} \\
&\geq \min_{l < i} \{w_s(\pi^l u_l)\} \\
&> w_s(x);
\end{align*}
while for $n > v_p(\pi^j)$,
\begin{align*}
s v_n(x) + n &\geq \min \{ \min_{l \in [i,j]} \{
s v_n(\pi^l u_l)+n \}, \min_{l \notin [i,j]} \{ s v_n(\pi^l u_l) +n\}
\} \\
&\geq \min \{ \min_{l \in [i,j]} \{ s v_n(\pi^l u_l) +n\},
\min_{l \notin [i,j]} \{w_s(\pi^l u_l)\} \}.
\end{align*}
For $l \in [i,j]$, $n > v_p(\pi^l)$ and $u_l$ is a semiunit,
so $s v_n(\pi^l u_l)+n > w_s(\pi^l u_l) \geq w_s(x)$;
for $l \notin [i,j]$,
$w_s(\pi^l u_l) > w_s(x)$ by the choice of $i$ and $j$. Putting
the inequalities together, we again conclude $s v_n(x) + n
> w_s(x)$.

Therefore the endpoints of the segment of the Newton polygon of $x$ of
slope
$-s$ are
$(v_0(u_i), v_p(\pi^i))$ and $(v_0(u_j), v_p(\pi^j))$.
Thus the multiplicity of $s$ as a slope of $x$ is
$v_p(\pi^j) - v_p(\pi^i) = s(v_0(u_i) - v_0(u_j))$,
which is indeed $s$ times the length of $L_s(x)$, as claimed.
\end{proof}
\begin{cor}
Let $x$ and $y$ be nonzero elements of $\Gamma_{\an,r}$. Then the
multiplicity of a slope $s$ of $xy$ is the sum of its multiplicities as
a slope of $x$ and of $y$.
\end{cor}
\begin{proof}
This follows immediately from the previous lemma and
the multiplicativity 
of the leading terms map $L_s$.
\end{proof}
\begin{cor} \label{cor:units}
The units of $\Gancon$ are precisely those $x$ with
$v_n(x) = \infty$ for some $n$.
\end{cor}
\begin{proof}
A unit of $\Gancon$ must also be a unit in $\Gamma_{\an,r}$ for some
$r$, and a unit of
$\Gamma_{\an,r}$ must have all slopes of multiplicity zero. (Remember,
in $\Gamma_{\an,r}$, any slopes greater than $r$ are disregarded.)
If $v_n(x) < \infty$ for all $n$, then $x$ has infinitely many
different slopes, so it still has slopes of nonzero multiplicity
in $\Gamma_{\an,r}$ for
any $r$, and so can never become a unit.
\end{proof}
We again caution that the condition $v_n(x) = \infty$ does not
imply that $x \in \Gcon[\fp]$, if $K$ is not finite over $k((t))$.

It will be convenient to put elements of $x$ into a standard
(multiplicative) form, so we make a statement to this effect
as a lemma.
\begin{lemma} \label{lem:sorite}
For any $x \in \Gamma_{\an,r}$ nonzero, there exists
a unit $u \in \Gamma_{\an,r}$ such that $ux$ admits
a semiunit decomposition $\sum_i u_i \pi^i$ with
$u_0 = 1$ and $u_i = 0$ for $i>0$. Moreover, for such $u$,
we have
\begin{enumerate}
\item[(a)] $v_0(ux) = 0$; 
\item[(b)] $w_r(ux) = 0$;
\item[(c)] $r v_n(ux) + n > 0$ for $n > 0$;
\item[(d)] the Newton polygon of $ux$ begins at $(0,0)$.
\end{enumerate}
\end{lemma}
\begin{proof}
By Proposition~\ref{prop:semidec2}, we can find
a semiunit decomposition $\sum_i u'_i \pi^i$
of $x$;
then $u'_i = 0$ for $i$ sufficiently large. Choose the
largest $j$ such that $u'_j \neq 0$, and put $u = 
\pi^{-j} (u'_j)^{-1}$. Then $ux$ admits the semiunit decomposition
$\sum_i u_i \pi^i$ with $u_i = u'_{i+j}/u'_j$, so $u_0 = 1$
and $u_i =0$ for $i>0$.

To verify (a), note that
$r v_0(ux) \geq \min_i \{r v_0(u_i \pi^i)\} \geq 0$,
and the minimum is only achieved for $i=0$: for $i<0$,
$r v_0(u_i \pi^i) > w_r(u_i \pi^i) \geq 0$ since $\sum_i
u_i \pi^i$
is a semiunit decomposition of $ux$. Thus $r v_0(ux) = 0$, whence (a).

To verify (b), note that
$w_r(ux) \geq \min_i \{w_r(u_i \pi^i)\} = 0$, whereas
$w_r(ux) \leq r v_0(ux) = 0$ from (a).

To verify (c), note that for $n > 0$ and $m = v_p(\pi^i)$,
$r v_n (u_i \pi^i) + n > r v_m(u_i \pi^i) + m \geq w_r(u_i \pi^i)
\geq 0$, so
$r v_n(ux) + n \geq \min_{i \leq 0} \{r v_n(u_i \pi^i) + n\} >0$.

To verify (d), first note that
the line through $(0,0)$ of slope $-r$ is a lower supporting
line of the set of points $(v_n(ux), n)$, since 
$r v_n(ux) + n \geq w_r(ux) \geq 0$ for $n \leq 0$.
Thus $(0,0)$ lies on the Newton polygon, and the slope of
the segment of the Newton polygon just to the right of $(0,0)$
is at least $-r$.
We also have $r v_n(ux) + n > 0$ for $n>0$, so the slope of the
segment of the Newton polygon just to the left of $(0,0)$,
if there is one, must be less than $-r$. Thus the first segment
of slope at least $-r$ does indeed begin at $(0,0)$, as desired.
\end{proof}

The next lemma may be viewed as a version of the Weierstrass preparation
theorem.
\begin{lemma} \label{lem:firstslope}
Let $x$ be a nonzero element of $\Gamma_{\an,r}$ whose largest
slope is $s_1$ with multiplicity $m > 0$. Then there exists $y
\in \Gamma_{\an,r}$, pure of slope
$s_1$ with multiplicity $m$, which divides $x$.  
\end{lemma}
\begin{proof}
If $x$ is pure of slope $s_1$, there is nothing to prove. So
assume that $x$ is not pure, and let $s_2$ be the second largest
slope of $x$.

By Lemma~\ref{lem:sorite}, there exists a unit $u \in \Gamma_{\an,r}$
such that $ux$ admits a semiunit decomposition 
$\sum_i u_i \pi^i$ with $u_0 =1$ and $u_i = 0$ for $i>0$,
the slopes of $x$ and $ux$ occur with the same multiplicities,
and the first segment of the Newton polygon of $ux$
has left endpoint $(0,0)$. Since that segment has
slope $-s_1$ and multiplicity $m$, its right endpoint is
$(m/s_1, -m)$. Put $M = -m/v_p(\pi)$; then $w_{s_1}(u_M \pi^M) = 
0$ and $w_{s_1}(u_i \pi^i) > 0$ for $i < M$.

We first construct a sort of ``Mittag-Leffler'' decomposition
of $ux$. Put $X = u x \pi^{-M} u_M^{-1}$, and
set $y_0 = z_0 = 1$. Given $y_l$ and $z_l$ for some
$l$, let $\sum_i w_i \pi^i$ be a semiunit decomposition
of $X - y_l z_l$. Put
\begin{align*}
y_{l+1} &= y_l + \sum_{v_0(w_i) < 0} w_i \pi^i \\
z_{l+1} &= z_l + \sum_{v_0(w_i) \geq 0} w_i \pi^i.
\end{align*}
Given $s$ with $s_2 < s < s_1$, put
$c_s = w_s(X - 1)$, so that $c_s > 0$.
We show that for each $l$, $w_s(y_l - 1)\geq c_s$, $w_s(z_l - 1) \geq c_s$,
and $w_s(X - y_l z_l) \geq (l+1)c_s$.
These inequalities are clear for $l=0$. If they hold for $l$, then
\begin{align*}
w_s(y_{l+1}-1) &\geq \min\{w_s(y_l - 1), w_s(y_{l+1}-y_l)\} \\
&\geq \min\{c_s, (l+1)c_s\} = c_s,
\end{align*}
and similarly $w_s(z_{l+1}-1) \geq c_s$. As for the third inequality,
note that
\begin{align*}
X - y_{l+1} z_{l+1} &=
X - y_l z_l + y_l (z_l - z_{l+1}) + z_{l+1} (y_l - y_{l+1}) \\
&= (y_l - 1)(z_l - z_{l+1}) + (z_{l+1} - 1)(y_l - y_{l+1})
\end{align*}
since $X - y_l z_l = (y_{l+1} - y_l) + (z_{l+1} - z_l)$.
Since $w_s(y_l - y_{l+1}) \geq (l+1)c_s$ and
$w_s(z_l - z_{l+1}) \geq (l+1)c_s$, we conclude that
\begin{align*}
w_s(X - y_{l+1} z_{l+1}) &\geq
 \min\{w_s((y_l - 1)(z_l - z_{l+1})), w_s((z_{l+1} - 1)(y_l - y_{l+1}))\} \\
&\geq \min\{c_s + (l+1)c_s, c_s + (l+1)c_s \} \\
&= (l+2)c_s,
\end{align*}
as desired. This completes the induction.

We do not yet know that either $\{y_l\}$ or $\{z_l\}$ converges
in $\Gamma_{\an,r}$; to get to that point, we need to play the
two sequences off of each other.
Suppose $s_3$ satisfies $s_2 < s_3 < s_1$.
Note that to get from $y_{l}$ to $y_{l+1}$, we add terms
of the form $w_i \pi^i$, with $w_i$ a semiunit,
for which $v_0(w_i) < 0$ but $s v_0(w_i) + v_p(\pi^i) \geq (l+1)c_s$
for $s_2 < s < s_1$.
This implies that
\[
s v_0(w_i) + v_p(\pi^i) \geq (l+1)c_{s_3}
\]
for all $s \leq s_3$. In particular, $w_s(y_{l+1} - y_l) \to \infty$
as $l \to \infty$, so $\{y_l\}$ converges to a limit in
$\Gamma_{\an,s}$ for any $s \leq s_3$. Moreover, for $s \leq s_3$, we have
\begin{align*}
w_s(y_{l+1} - 1) &\geq \min\{w_s(y_l - 1), w_s(y_{l+1} - y_l)\} \\
&\geq \min\{w_s(y_l - 1), (l+1)c_{s_3}\}
\end{align*}
so by induction on $l$, $w_s(y_l - 1) \geq c_{s_3}$.
Hence $y$ and each of the $y_l$ are units
in $\Gamma_{\an,s_3}$, for any $s_3 < s_1$.

On the flip side, to get from $z_l$ to $z_{l+1}$, we add terms
of the form $w_i \pi^i$, with $i$ a semiunit, for which
$v_0(w_i) > 0$ but $s v_0(w_i) + v_p(\pi^i) \geq (l+1)c$
for $s_2 < s < s_1$. This implies that $s v_0(w_i) + v_p(\pi^i)
\geq (l+1)c_{s_3}$ for all $s \geq s_3$. As in the previous paragraph,
we deduce $w_s(z_{l+1} - z_l) \to \infty$
and $w_s(z_l-1) > 0$ for $s_2 < s \leq r$.

Put $z = X y^{-1}$ in $\Gamma_{\an,s_3}$.
Since $w_{s_3}(y_l) = 0$ for all $l$, we have
\begin{align*}
w_{s_3}(z_l - z) &= w_{s_3}(y y_l z_l - y y_l z) \\
&= w_{s_3}(y (y_l z_l - X) + (y-y_l)X) \\
&\geq \min\{ w_{s_3}(y (y_l z_l - X)), w_{s_3}((y-y_l)X)\}
\end{align*}
and both terms in braces tend to infinity with $l$.
Thus $z_l \to z$ under $|\cdot|_{s_3}$.

For $s_3 \leq s \leq r$, since $s v_n(z_l - z) + n \geq (l+1)c_{s_3}$
and $s v_n(z_l) + n \to \infty$ as $n \to \pm \infty$,
for any given $l$ we have $s v_n(z) + n \geq (l+1)c_{s_3}$ 
for all but finitely many
$n$. Since this holds for any $l$, we have $s v_n(z) + n
\to \infty$ as $n \to \pm \infty$. As we already have
$z \in \Gamma_{\an,s_3}$, this is enough to imply $z \in \Gamma_{\an,r}$.
Meanwhile, put 
\[
a_l = X(1 + (1-z) + \cdots + (1-z)^l) = y(1 - (1 - z)^{l+1}),
\]
so
that $w_s(a_l - y) = (l+1)w_s(1-z)$ for $s_2 < s < s_1$. In
particular, for each $n$,
$v_n(a_m - y) \to \infty$ as $m \to \infty$,
and so the inequalities
\[
v_n(a_l - y) \geq \min\{v_n(a_l - a_{l+1}), \dots,
v_n(a_{m-1} - a_m), v_n(a_m - y)\}
\]
for each $m$ yield, in the limit as $m \to \infty$, the inequality
\[
v_n(a_l - y) \geq \min\{v_n(a_l - a_{l+1}), v_n(a_{l+1} - a_{l+2}), \dots \}.
\]
Now
$w_s(a_{l+1} - a_l) = w_s(X(1-z)^{l+1}) = 
w_s(X) + (l+1)w_s(1-z)$ for $s_2 < s \leq r$, so
$s v_n(a_{l+1} - a_l) + n \geq w_s(X) + (l+1)w_s(1-z)$. We conclude that
\[
s v_n(a_l - y) + n \geq w_s(X) + (l+1) w_s(1-z),
\]
so that $s v_n(y) + n \geq w_s(X) + (l+1) w_s(1-z)$ for all but finitely
many $n$. Therefore $s v_n(y) + n \to \infty$ as $n \to \pm \infty$
for $s_2 < s \leq r$. Again, since we already have $y \in \Gamma_{\an,s_3}$,
we deduce that $y \in \Gamma_{\an,r}$.

Since $y$ is a unit in $\Gamma_{\an,s}$ for any $s<s_1$,
it has no slopes less than $s_1$. Since $w_s(1-z) > 0$ for
$s_2 < s \leq r$, $z$ has no slopes greater than $s_2$. Since the
slopes of $y$ and $z$ together must comprise the slopes of $x$,
$y$ must have $s_1$ as a slope with multiplicity $m$ and no other
slopes, as desired.
\end{proof}

A \emph{slope factorization} of a nonzero element $x$ of $\Gamma_{\an,r}$ is
a Fr\'echet-convergent product $x = \prod_{j=1}^N x_j$ for 
$N$ a positive integer or $\infty$, where each $x_j$ is pure
and the slopes $s_j$ of $x_j$ satisfy $s_1 > s_2 > \cdots$.

\begin{lemma} \label{lem:factor}
Every nonzero element of $\Gamma_{\an,r}$ has a slope factorization.
\end{lemma}
\begin{proof}
Let $x$ be a nonzero element of $\Gamma_{\an,r}$ with slopes $s_1, s_2, \dots$.
By Lemma~\ref{lem:firstslope}, we can
find $y_1$ pure of slope $s_1$ dividing $x$ such that
 $x/y_1$ has largest slope $s_2$.
Likewise, we can find $y_2$ pure of slope $s_2$
such that $y_2$ divides $x/y_1$, $y_3$ pure of slope
$s_3$ such that $y_3$ divides
$x/(y_1y_2)$, and so on. 

If there are $N<\infty$ slopes, then
$x$ and $y_1\cdots y_N$ have the same slopes, so $x/(y_1 \cdots y_N)$
must be a unit $u$, and $x = (uy_1)y_2\cdots y_N$ is a slope factorization.
Suppose instead there are infinitely many slopes; then $s_i \to 0$
as $i \to \infty$.
By Lemma~\ref{lem:sorite}, for each $i$ we can find a unit $a_i$ such that
$a_i y_i$ admits a semiunit decomposition $\sum_j u_{ij} \pi^j$ with
$u_{i0} = 1$ and $u_{ij} = 0$ for $j>0$. For $j<0$,
$s v_n(u_{ij} \pi^j) + n$ is minimized for $n = v_p(\pi^j) < 0$ because
$u_{ij}$ is a semiunit; for $i$ sufficiently large, we have $s \geq s_i$,
so
\begin{align*}
s v_{v_p(\pi^j)}(u_{ij} \pi^j) + v_p(\pi^j) &= \frac{s}{s_i} (s_i v_{v_p(\pi^j)}(u_{ij} \pi^j) + v_p(\pi^j))
+ \left( \frac{s}{s_i} - 1 \right) (-v_p(\pi^j)) \\
&\geq \left( \frac{s}{s_i} - 1 \right) (-j) v_p(\pi),
\end{align*}
which tends to infinity as $i \to \infty$. Hence
$w_s(a_i y_i - 1) \to \infty$ as $i \to \infty$; if we put
$z_j = \prod_{i=1}^j a_i y_i$; then
$\{z_j\}$ converges to a limit $z$, and $\{x/z_j\}$ converges
to a limit $u$, such that $uz = x$.
The slopes of $z$ coincide with the slopes of $x$, so $u$ must
be a unit, and
$(u a_1 y_1) \prod_{i>1} (a_i y_i)$ is a slope factorization of $x$.
\end{proof}

\begin{lemma} \label{lem:rollup}
Let $x$ be an element of $\Gamma_{\an,r}$ which is pure of slope
$s$ and multiplicity $m$. Then for every $y \in \Gamma_{\an,r}$,
there exists $z \in \Gamma_{\an,r}$ such that:
\begin{enumerate}
\item[(a)] $y-z$ is divisible by $x$;
\item[(b)] $w_s(z) \geq w_s(y)$;
\item[(c)] $v_n(z) = \infty$ for $n < 0$.
\end{enumerate}
\end{lemma}
\begin{proof}
Put $M = m/v_p(\pi)$.
By Lemma~\ref{lem:sorite}, there exists a unit $u \in \Gamma_{\an,r}$
such that $xu$ admits a semiunit decomposition $\sum_{i=-M}^0
x_i \pi^i$ with $sv_0(x_{-M}) = m$.
Note that
\[
w_r(x_{-M} \pi^{-M}) = r v_{-m}(x_{-M} \pi^{-M}) - m
= m\left(\frac{r}{s} - 1\right).
\]
Let $\sum_i y_i \pi^i$ be a semiunit decomposition of $y$.

We define the sequence $\{c_l\}_{l=0}^\infty$ of elements of $\Gamma_{\an,r}$
such that $v_n(c_l) = \infty$ for $n < 0$, 
$w_r(c_l) \geq -l (v_p(\pi)+m(r/s-1))$, 
$w_s(c_l) \geq -l v_p(\pi)$, and
\[
c_l \equiv \pi^{-l} \pmod{x}.
\]
Put $c_0 = 1$ to start. Given $c_l$, let $\sum_{i} u_i \pi^i$ 
be a semiunit decomposition of $c_l$; since $v_n(c_l) = \infty$
for $n<0$,
we have $u_i = 0$ for $i < 0$. Now set
\[
c_{l+1} = \pi^{-1} (c_l - u x x_{-M}^{-1} \pi^{M} u_0 ).
\]
The congruence $c_{l+1} \equiv \pi^{-1} c_l \equiv \pi^{-l-1}
\pmod{x}$ is clear from the definition.
Since $ux x_{-M}^{-1} \pi^{M} \equiv 1 \pmod{\pi}$, the term
in parentheses has positive valuation, so $v_n(c_{l+1}) = \infty$
for $n<0$. 
Since $w_s(u x) = w_s(x_{-M} \pi^{-M}) = 0$ and $w_s(u_0) \geq w_s(c_l)$,
we have $w_s(c_{l+1}) \geq w_s(\pi^{-1}c_l) \geq -(l+1) v_p(\pi)$.
Finally, $w_r(u_0) \geq w_r(c_l)$, $w_r(ux) = 0$
and $w_r(x_{-M} \pi^{-M}) = m(r/s - 1)$, so
\begin{align*}
w_r(c_{l+1}) &\geq w_r(\pi^{-1}) +
\min\{w_r(c_l), w_r( u x x_{-M}^{-1} \pi^{M} u_0) \} \\
&\geq -v_p(\pi) + w_r(c_l) - m(r/s-1) \\
&\geq -(l+1)(m(r/s-1) + v_p(\pi)).
\end{align*}

We wish to show that $\sum_{i=-\infty}^{-1} y_i c_{-i}$ converges, so that
its limit is congruent to $\sum_{i=-\infty}^{-1} y_i \pi^i$ modulo $x$.
To this end, choose $t>0$ large enough that
\[
tr v_p(\pi) > m(r/s-1) + v_p(\pi).
\]
Then $(1/t) v_n(y) + n \to \infty$ as $n \to -\infty$,
so in particular there exists $c>0$ such that 
$(1/t) v_n(y) \geq -n - c$ for $n < 0$. 
For $n = v_p(\pi^i)$ where $i<0$ and $y_i \neq 0$, we have
$v_n(y) = v_0(y_i)$, so we have $v_0(y_i) \geq -ti v_p(\pi) - tc$.
Then
\begin{align*}
w_r(y_i c_{-i}) &= w_r(y_i) + w_r(c_{-i}) \\
&= r v_0(y_i) + w_r(c_{-i}) \\
&\geq -tri v_p(\pi) - trc + i (m(r/s-1) + v_p(\pi))
\end{align*}
which tends to infinity as $i \to -\infty$.
Thus $\sum_{i=-\infty}^{-1} y_i c_{-i}$ converges under $|\cdot|_r$;
since $v_n(y_i c_{-i}) = \infty$ for $n < 0$, the sum also converges
under $|\cdot|_s$ for $0<s<r$. Thus it has a limit $z'
\in \Gamma_{\an,r}$; put $z = z' + \sum_{i=0}^\infty y_i \pi^i$.
Then $y-z = \sum_{i=-\infty}^{-1} y_i (\pi^i - c_{-i})$; since each term
in the sum is divisible by $x$, so is the sum. This verifies (a).
To verify (b), note that $w_s(y_i c_{-i}) \geq w_s(y_i \pi^i)$ for
$i < 0$, so $w_s(z') \geq w_s(y)$, and clearly $w_s(z-z')
\geq w_s(y)$, so $w_s(z) \geq w_s(y)$.
To verify (c), simply note that each term in the sum defining
$z$ satisfies the same condition.
\end{proof}

\subsection{The B\'ezout property for analytic rings}
\label{subsec:anbezout2}

Again, we assume that the valued field
$K$ satisfies the conditions of Proposition~\ref{prop:uniform},
so that $\Gancon = \GKancon$ is defined.
With the factorization results of the previous section in hand, we
now focus on establishing the B\'ezout property for $\Gancon$
(Theorem~\ref{thm:bezout}). We proceed by establishing principality
of successively more general classes of finitely generated ideals,
culminating in the desired result.

\begin{lemma} \label{lem:gcdfinmix}
Let $x$ and $y$ be elements of $\Gamma_{\an,r}$, each with
finitely many slopes, and having no slopes in common.
Then the ideal $(x,y)$ is the unit ideal.
\end{lemma}
\begin{proof}
We induct on the sum of the multiplicities of the slopes
of $x$ and $y$; the case where either $x$ or $y$ has total
multiplicity zero is vacuous, as then $x$ or $y$ is a unit
and so $(x,y)$ is the unit ideal. So we assume that both $x$ and
$y$ have positive total multiplicity.

If $x$ is not pure, then by Lemma~\ref{lem:factor} it factors
as $x_1 x_2$, where $x_1$ is pure and $x_2$ is not a unit. By
the induction hypothesis, the ideals $(x_1,y)$ and $(x_2, y)$
are the unit ideal; in other words, $x_1$ and $x_2$ have multiplicative
inverses modulo $y$. In that case, so does $x = x_1 x_2$,
so $(x,y)$ is the unit ideal.
The same argument applies in case $y$ is not pure. 

It thus remains
to treat the case where $x$ and $y$ are both pure. Let $s$ and $t$
be the slopes of $x$ and $y$, and let $m$ and $n$ be the corresponding
multiplicities. Put $M = m/v_p(\pi)$ and $N = n/v_p(\pi)$.
Without loss of generality, we may assume $s < t$.
By Lemma~\ref{lem:sorite}, we can find units $u$ and $v$ such that
$ux = \sum_{i=-M}^0 x_i \pi^i$ and $vy = \sum_{i=-N}^0 y_i \pi^i$.

Put
\[
X = ux \pi^M x_{-M}^{-1}, \qquad
Y = vy \pi^N y_{-N}^{-1}, \qquad
z = X - Y.
\]
We can read off information about the Newton polygon of $z$ by comparing
$w_r(X)$ with $w_r(Y)$;
see Figure~\ref{fig:newton2} for an illustration. (In both
diagrams, the dashed lines
have slope $-r$.)
If $w_r(X) < w_r(Y)$ (left side
of Figure~\ref{fig:newton2}),
then the highest vertex of the
lower convex hull of the set of points $(v_l(z), l)$ occurs
at $(v_m(X), m)$ and the lowest vertex
has positive $y$-coordinate. Moreover, the
slope of the first segment
of the lower convex hull is at least $-s$.
Thus the sum of all multiplicities
of $z$ is strictly less than $m$, and $y$ and $z$ have no
common slopes, so the induction hypothesis
implies that $(x,y) = (y,z)$ is the unit ideal.

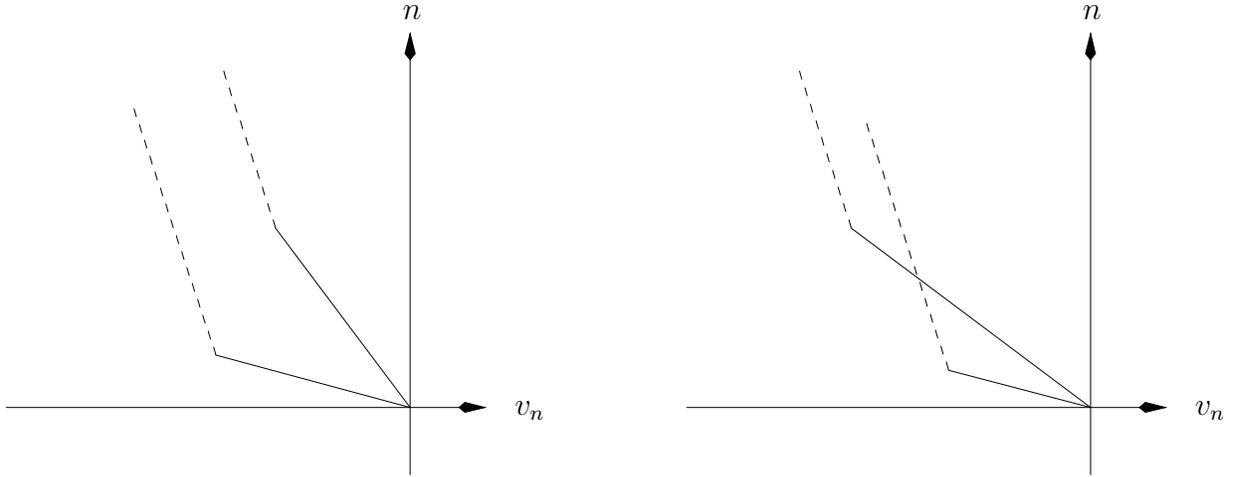
\begin{figure}[ht]
  \centering
\setlength{\unitlength}{0.00087489in}
\begingroup\makeatletter\ifx\SetFigFont\undefined
\def\x#1#2#3#4#5#6#7\relax{\def\x{#1#2#3#4#5#6}}%
\expandafter\x\fmtname xxxxxx\relax \def\y{splain}%
\ifx\x\y   
\gdef\SetFigFont#1#2#3{%
  \ifnum #1<17\tiny\else \ifnum #1<20\small\else
  \ifnum #1<24\normalsize\else \ifnum #1<29\large\else
  \ifnum #1<34\Large\else \ifnum #1<41\LARGE\else
     \huge\fi\fi\fi\fi\fi\fi
  \csname #3\endcsname}%
\else
\gdef\SetFigFont#1#2#3{\begingroup
  \count@#1\relax \ifnum 25<\count@\count@25\fi
  \def\x{\endgroup\@setsize\SetFigFont{#2pt}}%
  \expandafter\x
    \csname \romannumeral\the\count@ pt\expandafter\endcsname
    \csname @\romannumeral\the\count@ pt\endcsname
  \csname #3\endcsname}%
\fi
\fi\endgroup
{\renewcommand{\dashlinestretch}{30}
\begin{picture}(7167,2937)(0,-10)
\path(6537,417)(5682,642)
\path(6537,417)(5097,1497)
\path(4107,417)(6987,417)
\blacken\path(6867.000,387.000)(6987.000,417.000)(6867.000,447.000)(6831.000,417.000)(6867.000,387.000)
\path(6537,12)(6537,2667)
\blacken\path(6567.000,2547.000)(6537.000,2667.000)(6507.000,2547.000)(6537.000,2511.000)(6567.000,2547.000)
\path(2442,417)(1272,732)
\path(2442,417)(1632,1497)
\dashline{60.000}(1632,1497)(1317,2442)
\dashline{60.000}(1272,732)(777,2217)
\path(12,417)(2892,417)
\blacken\path(2772.000,387.000)(2892.000,417.000)(2772.000,447.000)(2736.000,417.000)(2772.000,387.000)
\path(2442,12)(2442,2667)
\blacken\path(2472.000,2547.000)(2442.000,2667.000)(2412.000,2547.000)(2442.000,2511.000)(2472.000,2547.000)
\dashline{60.000}(5682,642)(5187,2127)
\dashline{60.000}(5097,1497)(4782,2442)
\put(2397,2757){\makebox(0,0)[lb]{\smash{{{\SetFigFont{12}{14.4}{rm}$n$}}}}}
\put(3072,372){\makebox(0,0)[lb]{\smash{{{\SetFigFont{12}{14.4}{rm}$v_n$}}}}}
\put(6492,2757){\makebox(0,0)[lb]{\smash{{{\SetFigFont{12}{14.4}{rm}$n$}}}}}
\put(7167,372){\makebox(0,0)[lb]{\smash{{{\SetFigFont{12}{14.4}{rm}$v_n$}}}}}
\end{picture}
}

  \caption{The Newton polygons of $X = ux \pi^{M} x_{-M}^{-1}$ and
$Y = vy \pi^{N} y_{-N}^{-1}$}
  \label{fig:newton2}
\end{figure}

If $w_r(X) \geq w_r(Y)$
(right side of Figure~\ref{fig:newton2}),
then the highest vertex of the
lower convex hull of the set of points $(v_l(z), l)$ occurs
at $(v_n(Y), n)$ and the lowest vertex
has positive $y$-coordinate. Moreover,
$(v_m(X), m)$ is also a vertex of the lower
convex hull, and the line joining it to
$(v_n(Y), n)$ is a support line of the lower
convex hull. Thus the segment joining the two points is a segment
of the lower convex hull, of
slope less than $-t$; the remainder
of the lower convex hull consists of segments of slope at
least $-s$, of total multiplicity less than $m$.
By Lemma~\ref{lem:factor}, $z$ factors as $z_1 z_2$, where $z_1$ is
pure of some slope greater than $t$, and $z_2$ has all slopes
less than or equal to $s$ and total multiplicity less than $m$.
By the induction hypothesis, $(x, z_1)$ and $(y, z_2)$ both
equal the unit ideal. But $(y,z_1) = (x,z_1)$ since $z_1$
divides $z = ux \pi^M x_{-M}^{-1} - vy \pi^N y_{-N}^{-1}$,
so $(y, z_1z_2) = (y,z) = (x,y)$ is also
equal to the unit ideal.

We conclude that the induction goes through for all $x$ and
$y$. This completes the proof.
\end{proof}

\begin{lemma} \label{lem:gcdpure}
Let $x$ and $y$ be elements of $\Gamma_{\an,r}$ with $x,y$ pure of
the same slope $s$.
Then $(x,y)$ is either the unit ideal or is generated by a
pure element of slope $s$.
\end{lemma}
\begin{proof}
Let $m$ and $n$ be the multiplicities of $s$ as a slope of $x$ and $y$;
we induct on $m+n$. Put $M = m/v_p(\pi)$ and $N = n/v_p(\pi)$.
By Lemma~\ref{lem:sorite}, we may choose units $u,v$ so that
$ux$ and $vy$ admit semiunit decompositions
$ux = \sum_{i=-M}^0 x_i \pi^i$ and $vy = \sum_{i=-N}^0 y_i \pi^i$.
Put
\[
X = ux \pi^M x_{-M}^{-1}, \qquad
Y = vy \pi^N y_{-N}^{-1}, \qquad
z = X - Y.
\]
By symmetry, we may assume $m \leq n$ without loss of generality.

First suppose $m < n$. Then the highest vertex of the lower
convex hull of the set of points $(v_l(z), l)$ occurs at
$(v_n(Y), n)$, and the lowest vertex has 
positive $y$-coordinate. Thus the total multiplicity of $z$ is
strictly less than $n$. By Lemma~\ref{lem:factor}, we may write
$z = z_1 z_2$, where $z_1$ has no slopes equal to $s$, and
$z_2$ is either a unit or is pure of slope $s$. By
Lemma~\ref{lem:gcdfinmix}, $(x, z_1)$ is the unit ideal, so
$(x, y) = (x,z) = (x,z_2)$, which is principal by the
induction hypothesis.

Next suppose $m = n$.
We define a sequence $\{z_l\}$ as follows, starting
with $z_0 = z = X - Y$.
Given $z_l$, let $\sum_i u_i \pi^i$ be a semiunit decomposition of
$z_l$, and set
\[
z_{l+1} = z_l - \sum u_i \pi^i Y
\]
taking the sum over all $i$ for which $w_r(u_i \pi^i) > 0$.
Clearly $z_l - z_0$ is divisible by $y$ for each $l$.

Suppose $z_l = z_{l+1}$ for some $l$. If $\sum_i u_i \pi^i$ is
the chosen semiunit decomposition of $z_l$, then
$w_r(u_i \pi^i) \leq 0$ for each $i$.
On the other hand, we have $w_r(z_l - z_0) > w_r(Y)$ by induction
on $l$.

We claim that the multiplicity of $s$ as a slope of $z_l$ is less than $n$.
Suppose that this is not the case; then there exist indices
$i<j$ with $w_r(u_i \pi^i) - w_r(u_j \pi^j) = -w_r(Y)$.
As shown above, we have $w_r(u_i \pi^i) \leq 0$; since
$w_r(z_l-z_0) > w_r(Y)$ and 
$w_r(z_0) \geq w_r(Y)$,
we also have $w_r(u_j \pi^j) \geq w_r(Y)$.
Since these inequalities imply $w_r(u_i \pi^i) - w_r(u_j \pi^j) \leq -w_r(Y)$,
they must all be equalities; that is,
$w_r(u_i \pi^i) = 0$, $w_r(u_j \pi^j) = w_r(Y)$, and
$w_r(z_0) = w_r(Y)$. This means that the highest vertex
of the lower convex hull of the set of points $(v_c(z_l), c)$
occurs at $(v_n(Y), n)$, and there is another vertex with
$y$-coordinate at most $0$. This is a contradiction.

Still under the hypothesis $z_l = z_{l+1}$, we conclude that
the multiplicity of $s$ as a slope of $z_l$ is less than $n$.
By Lemma~\ref{lem:factor}, we can write $z_l = a_1 a_2$,
where $a_1$ has no slopes equal to $s$, and $a_2$ is pure
of slope $s$ and multiplicity less than $n$. By
Lemma~\ref{lem:gcdfinmix}, $(y, a_1)$ is the unit ideal;
by the induction hypothesis, $(y, a_2)$ is principal. Thus
$(x,y) = (y, z_0) = (y, z_l) = (y, a_1a_2) = (y, a_2)$
is principal, as desired.

Finally, suppose that $z_l \neq z_{l+1}$ for all $l$.
We claim that $z_l \equiv 0 \pmod{\pi^{l+1}}$ for each $l$.
This holds for $l=0$ by the choice of $z_0$. Given the congruence
for $z_l$, let $\sum_i u_i \pi^i$ be the chosen semiunit
decomposition of $z_l$. By hypothesis, $u_i = 0$ for $i \leq l$.
If $u_{l+1} \neq 0$, then we must have $w_r(u_{l+1} \pi^{l+1}) >0$,
or else $w_r(u_i\pi^i) \leq 0$ for all $i$ (recall that in a semiunit
decomposition, whenever $i<j$ and $u_i,u_j \neq 0$ we have
$w_r(u_i \pi^i) \geq w_r(u_j \pi^j)$). Thus
$z_{l+1} \equiv z_l - u_{l+1} \pi^{l+1} \equiv 0 \pmod{\pi^{l+2}}$.
Thus the congruence holds by induction.

Since the $z_l$ converge to zero $\pi$-adically and $w_r(z_l)$
is bounded below, the $z_l$ also
converge to zero in the Fr\'echet topology.
This means the series $\sum_{l=0}^\infty (z_l - z_{l+1})$
converges to $z_0$, the series $\sum_{l=0}^\infty (z_l - z_{l+1})
/ y$ converges to some limit $a$, and we have $y a = z_0$.
Therefore $z_0$ is divisible by $y$, as is $x$, and the ideal
$(x,y)$ is generated by $y$.

This completes the induction in all cases, whence the desired
result.
\end{proof}
\begin{cor} \label{cor:gcdpure}
For $x,y \in \Gamma_{\an,r}$ with $x$ pure of slope $s$,
the ideal $(x,y)$ is principal.
\end{cor}
\begin{proof}
By Lemma~\ref{lem:rollup}, there exists $z \in \Gamma_{\an,r}$
such that $y-z$ is divisible by $x$ and $v_n(z) = \infty$ for
$n < 0$. Thus $z$ has only finitely many slopes. By
Lemma~\ref{lem:factor}, we can factor $z$ as $z_1 z_2$,
where $z_1$ is pure of slope $s$ and $z_2$ has no slopes equal
to $s$. Then $(x,z_2)$ is the unit ideal, so $(x,y)
= (x,z) = (x,z_1)$, which is principal by Lemma~\ref{lem:gcdpure}.
\end{proof}

\begin{lemma}[Principal parts theorem] \label{lem:crt}
Let $s_n$ be a decreasing sequence of positive rationals 
with limit $0$,
and suppose
$x_n \in \Gamma_{\an,r}$ is pure of slope $s_n$ for all $n$.
Then for any sequence $y_n$ of elements of $\Gamma_{\an,r}$,
there exists $y \in \Gamma_{\an,r}$ such that $y \equiv y_n \pmod{x_n}$
for all $n$.
\end{lemma}
\begin{proof}
As in the proof of Lemma~\ref{lem:factor}, we can replace each $x_n$
with itself times a unit, in such a way that $\prod_n x_n$ converges.
Put $x = \prod_n x_n$ and $u_n = x/x_n$.
By Lemma~\ref{lem:gcdpure}, $x_n$ is coprime to each of
$x_1, \dots, x_{n-1}$. By Corollary~\ref{cor:gcdpure}, the
ideal $(x_n, \prod_{i>n} x_i)$
is principal, but if it were not the unit ideal,
any generator would both be pure of slope $s_n$ and have
all slopes less than $s_n$. Thus $x_n$ is coprime to $\prod_{i>n} x_i$,
hence also to $u_n$.

We construct a sequence $\{z_n\}_{n=1}^\infty$ such that
$u_n z_n \equiv y_n \pmod{x_n}$ and $\sum u_n z_n$ converges
for the Fr\'echet topology;
then we may set $y = \sum u_n z_n$ and be done.
For the moment, fix $n$ and
choose $v_n$ with $u_n v_n \equiv y_n \pmod{x_n}$.

For $s > s_n$, we have $|1-x_n|_s < 1$, so the sequence
$c_m = -1 -(1-x_n) - \cdots - (1-x_n)^m$ is Cauchy for the norm
$|\cdot|_s$, and $|1 + c_m x_n|_s = |1-x_n|^{m+1}_s \to 0$ under $|\cdot|_s$.
In particular, for any $\epsilon > 0$,
there exists $m$ such that
$|1 + c_m x_n|_s < \epsilon$
for $s_{n-1} \leq s \leq r$.

Now choose $\epsilon_n >0$ such that 
$\epsilon_n |u_nv_n|_s < 1/n$ for all $s \geq s_{n-1}$
(with $n$ still fixed),
choose $m$ as above, and put $z_n = v_n(1 + c_m x_n)$.
Then $u_n z_n \equiv u_n v_n \equiv y_n \pmod{x_n}$.
Moreover, for any $s>0$, we have $s \geq s_{n-1}$ for 
sufficiently large $n$ since the $s_n$ tend to zero.
Thus for $n$ sufficiently large,
\begin{align*}
|u_n z_n|_s &= |u_n v_n(1 + c_m x_n)|_s \\
&< \epsilon_n |u_n v_n|_s  < 1/n.
\end{align*}
Hence $\sum_n u_n z_n$ converges with respect to 
$|\cdot|_s$ for $0 <s \leq r$, and its limit $y$ has the
desired property.
\end{proof}

At long last,
we are ready to prove the generalization of Lazard's result,
that $\Gamma_{\an,r}$ is a B\'ezout ring.
\begin{proof}[Proof of Theorem~\ref{thm:bezout}]
By induction on the number of generators of the ideal,
it suffices to prove that if $x,y \in \Gamma_{\an,r}$
are nonzero,
then the ideal $(x,y)$ is principal. 

Pick a slope factorization $\prod_j y_j$ of $y$. By
Corollary~\ref{cor:gcdpure},
we can choose a generator $d_j$ of $(x, y_j)$ for
each $j$, such that $d_j$ is either $1$ or is pure of the same slope
as $y_j$.
As in the proof of Lemma~\ref{lem:factor}, we can choose the $d_j$ so that
$\prod_j d_j$ converges.
Since the $d_j$ are pairwise coprime
by Lemma~\ref{lem:gcdfinmix},
$x$ is divisible by the product of any finite subset of the
$d_j$, and hence by $\prod_j d_j$. 

Choose
$a_j$ and $b_j$ such that $a_jx + b_jy_j = d_j$, and apply 
Lemma~\ref{lem:crt} to find $z$ such that
$z \equiv a_j \prod_{k \neq j} d_j \pmod{y_j}$ for each $j$. 
Then $zx-\prod_j d_j$ is divisible by
each $y_j$, so it is divisible by $y$, and so
$\prod_j d_j$ generates the ideal $(x,y)$. Thus $(x,y)$ is principal
and the proof is complete.
\end{proof}
\begin{cor} \label{cor:bezout}
For $K$ a finite extension of
$k((t))$, the ring $\GK_r[\fp]$ is a B\'ezout ring.
\end{cor}
\begin{proof}
For $x,y \in \GK_r[\fp]$, Theorem~\ref{thm:bezout} implies
that the ideal $(x,y)$ becomes principal in $\GK_{\an,r}$.
Let $d$ be a generator; then $d$ must have finite total multiplicity,
and so belongs to $\GK_{r}[\fp]$.

Put $x'= x/d$ and $y' = y/d$, so that $(x', y')$ becomes the
unit ideal in $\GK_{\an,r}$.
By Lemma~\ref{lem:factor},
$x'$ factors in $\GK_{\an,r}$ as $a_1 \dots a_l$, where 
each $a_i$ is pure. Since each of those factors has
finite total multiplicity, each lies in $\GK_r[\fp]$.

Since $(x', y')$ is the unit ideal in $\GK_{\an,r}$, so is
$(a_i, y')$ for each $i$. That is, there exist $b_i$ and $c_i$
in $\GK_{\an,r}$
such that $a_i b_i + c_i y' = 1$. Since $a_i$ is pure, 
Lemma~\ref{lem:rollup} implies that $c_i \equiv d_i \pmod{a_i}$
for some $d_i$ with finite total multiplicity, which thus belongs to
$\GK_{r}[\fp]$. Now $d_i y' \equiv 1 \pmod{a_i}$, and
$e_i = (d_i y' - 1)/a_i$ has finite total multiplicity, so itself
lies in $\GK_r[\fp]$. We now have the relation
$a_i e_i + d_i y' =1$ within $\GK_r[\fp]$, so $(a_i, y')$
is the unit ideal in $\GK_r[\fp]$. Since this is true for each
$i$, $(x', y')$ is also the unit ideal and so
$(x,y) = (d)$.

We conclude that any ideal generated by two elements is principal.
By induction, this implies that $\GK_r[\fp]$ has the B\'ezout
property.
\end{proof}
One presumably has the same result if $K$ is perfect, but it does
not follow formally from Theorem~\ref{thm:bezout}, since
$\GK_r$ is not Fr\'echet complete in $\GK_{\an,r}$. That is,
an element of $\GK_{\an,r}$ of finite total multiplicity need not
lie in $\GK_r$. So one must repeat the arguments used to prove 
Theorem~\ref{thm:bezout} working withing $\GK_r[\fp]$; as we have
no use for the result, we leave this to the reader.

\section{The special Newton polygon}
\label{sec:special}

In this chapter, we construct a Newton polygon for $\sigma$-modules
over $\Gancon$, the ``special Newton polygon''. More precisely,
we give a slope filtration over $\Galgancon$
that, in case the $\sigma$-module
is quasi-unipotent, is precisely the filtration that makes it quasi-unipotent.
The special Newton polygon is a numerical invariant of this filtration.

Throughout this chapter, we assume
$K$ is a valued field satisfying the condition of
Proposition~\ref{prop:uniform}.
The choice of $K$ will only be relevant once or twice, as most of
the time
we will be working with $\Galgancon = \Gamma^{k((t))^{\alg}}_{\an,\con}$.
When this is the case,
we will also assume $k$ is algebraically closed and that $\pi^\sigma = \pi$.

We will use without further comment the facts that every
element of $\Galgancon$ has a strong semiunit decomposition
(Proposition~\ref{prop:semidec3}) and that $\Gancon$ and $\Galgancon$
are B\'ezout rings (Theorem~\ref{thm:bezout}). In particular, any
$\sigma$-module over $\Gancon$ or $\Galgancon$ is free by
Proposition~\ref{prop:free}, so admits a basis.

\subsection{Properties of eigenvectors}

Recall that we call a nonzero element $\bv$ of a $\sigma$-module
$M$ an \emph{eigenvector} if there exists $\lambda \in \calO[\fp]$
such that $F\bv = \lambda \bv$. Also recall that if $\bv$ an
eigenvector, the \emph{slope} of $\bv$ is defined to be
$v_p(\lambda)$. (Beware: this differs from the notion of slope
used in Section~\ref{subsec:anbezout}.)
Our method of constructing the special Newton polygon of a
$\sigma$-module
over $\Galgancon$ is to exhibit a basis of eigenvectors after
enlarging $\calO$ suitably. Before proceeding, it behooves us
to catalog some basic properties of eigenvectors of $\sigma$-modules
over $\Galgancon$. Some of these assertions will also hold
more generally over $\Gancon$ (i.e., for arbitrary $K$), so we distinguish between
$\Gancon$ and $\Galgancon$ in the statements below.

For $M$ a $\sigma$-module over $\Gancon$, we say $\bv \in M$ is
\emph{primitive} if $\bv$ extends to a basis of $M$. By Lemma~\ref{lem:det1},
if $\be_1, \dots, \be_n$ is a basis of $M$ and $\bv = \sum c_i \be_i$,
then $\bv$ is primitive if and only if the $c_i$ generate the unit
ideal.
\begin{lemma} \label{lem:prim}
Let $M$ be a $\sigma$-module over $\Galgancon$. Then every eigenvector of $M$
is a multiple of a primitive eigenvector.  
\end{lemma}
\begin{proof}
Suppose $F\bv = \lambda \bv$. Choose a basis $\be_1, \dots, \be_n$,
put $\bv = \sum_i c_i \be_i$, and let $I$ be the ideal generated by
the $c_i$. Then $I$ is invariant under $\sigma$ and $\sigma^{-1}$. 
By Theorem~\ref{thm:bezout},
$I$ is principal; if $r$ is a generator of $I$, then 
$r^\sigma$ is also a generator.
Put $r^\sigma = cr$, with $c$ a unit,
and write $c = \mu d$, with $\mu \in \calO[\fp]$,
$v_0(d) < \infty$ and $v_n(d) = \infty$ for $n<0$.
By Proposition~\ref{prop:rank1},
 there exists a unit $s \in \Galgancon$ such that
$s^\sigma = ds$; then $(r/s)^\sigma = \mu (r/s)$. Therefore
$\sum_i s(c_i/r)\be_i$ is a primitive eigenvector of $M$ of which
$\bv$ is a multiple, as desired.
\end{proof}
A sort of converse to the previous statement is the following.
\begin{prop} \label{prop:raise0}
For $M$ a $\sigma$-module over $\Galgancon$, if $M$ contains
an eigenvector of eigenvalue $\lambda \in \calO[\fp]$,
then it contains an eigenvector of eigenvalue $\lambda \mu$
for any $\mu \in \calO$.
\end{prop}
\begin{proof}
Let $\bv \in M$ be a nonzero eigenvector with $F\bv = \lambda \bv$.
If $\mu$ is a unit, there exists a unit $c \in \calO$ such that
$c^\sigma = \mu c$. If $\mu$ is not a unit, then by
Proposition~\ref{prop:ansigeq}(b) there exists a nonzero $c
\in \Galgancon$ such that $c^\sigma = \mu c$. In either case,
we have $F(c\bv) = c^\sigma \lambda \bv = \lambda \mu (c\bv)$.
\end{proof}

\begin{prop} \label{prop:split}
Let $0 \to M_1 \to M \to M_2 \to 0$ be an exact sequence of $\sigma$-modules
over $\Galgancon$. Assume $M_1$ and $M_2$ have bases $\bv_1, \dots, \bv_m$
and $\bw_1, \dots, \bw_n$ of eigenvectors such that the slope of $\bv_i$
is less than or equal to the slope of $\bw_j$ for $1 \leq i \leq m$
and $1 \leq j \leq n$. Then the exact sequence splits over $\Galgancon$.
\end{prop}
\begin{proof}
Choose a basis $\bv_1, \dots, \bv_m, \bx_1, \dots, \bx_n$ of $M$ such
that $\bx_j$ projects onto $\bw_j$ in $M_2$ for $j=1, \dots, n$.
Suppose $F\bv_i = \lambda_i \bv_i$ for some $\lambda_i \in \calO[\fp]$.
Then $F\bx_j = \mu_j \bx_j + \sum_{i=1}^m A_{ij} \bv_i$ for some
$\mu_j \in \calO[\fp]$ and $A_{ij} \in \Galgancon$.
If $\by_j = \bx_j + \sum_{i=1}^m c_{ij} \bv_i$, then
\[
F\by_j = \mu_j \by_j + \sum_{i=1}^m (\lambda_i c_{ij}^\sigma - \mu_j 
c_{ij} + A_{ij}) \bv_i.
\]
By Proposition~\ref{prop:ansigeq}(a) and (b),
we can choose $c_{ij} \in \Galgancon$ for each $i,j$ so that
$\lambda_i c_{ij}^\sigma - \mu_j c_{ij} + A_{ij} = 0$.
For this choice, $\bv_1, \dots,
\bv_m, \by_1, \dots, \by_n$ 
form a basis of eigenvectors, so the exact sequence splits
as desired.
\end{proof}

\begin{prop} \label{prop:minslope}
Let $M$ be a $\sigma$-module over $\Gancon$ with a basis
$\bw_1, \dots, \bw_n$ such that $F\bw_i =
\mu_i \bw_i + \sum_{j < i} A_{ij} \bw_j$ for some $\mu_i \in \calO[\fp]$
and $A_{ij} \in \Gancon$. Then any eigenvector of $M$ has slope at least
$\min_i \{v_p(\mu_i)\}$.
\end{prop}
\begin{proof}
Let $\bv$ be any eigenvector of $M$, with $F\bv = \lambda \bv$. 
Write $\bv = \sum_i b_i \bw_i$ for some $b_i \in \Gancon$. Suppose that
$v_p(\lambda) < v_p(\mu_i)$ for all $i$. Then
\[
  \sum_i \lambda b_i \bw_i = F\bv = \sum_i b_i^\sigma \mu_i \bw_i + \sum_i b_i^\sigma \sum_{j<i} A_{ij} \bw_j.
\]
Comparing the coefficients of $\bw_n$ yields $\lambda b_n = \mu_n b_n^\sigma$,
which implies $b_n = 0$ by Proposition~\ref{prop:ansigeq}(c). Then comparing
the coefficients of $\bw_{n-1}$ yields $\lambda b_{n-1} = \mu_{n-1}
b_{n-1}^\sigma$, so $b_{n-1}=0$. Continuing in this fashion, we deduce
$b_1 = \cdots = b_n = 0$, contradiction. Thus $v_p(\lambda) \geq
v_p(\mu_i)$ for some $i$, as desired.  
\end{proof}

Recall that a sequence $(a_1, \dots, a_n)$ of real numbers
is said to \emph{majorize}
another sequence $(b_1, \dots, b_n)$ if $a_1 + \cdots +a_n = 
b_1 +\cdots + b_n$ and for $i=1, \dots, n-1$, the sum of
the $i$ smallest of $a_1, \dots, a_n$ is less than or equal to the sum of the
$i$ smallest of $b_1, \dots, b_n$. Note that two sequences majorize each
other if and only if they are equal up to permutation.

\begin{prop} \label{prop:basis}
Let $M$ be a $\sigma$-module over $\Gancon$ with a basis
$\bv_1, \dots, \bv_n$ of eigenvectors, with
$F\bv_i = \lambda_i \bv_i$ for $\lambda_i \in \calO[\fp]$.
Let $\bw_1, \dots, \bw_n$ be a basis of $M$ such that $F\bw_i =
\mu_i \bw_i + \sum_{j < i} A_{ij} \bw_j$ for some $\mu_i \in \calO[\fp]$
and $A_{ij} \in \Gancon$. Then the sequence $v_p(\mu_1), \dots,
v_p(\mu_n)$ majorizes the sequence $v_p(\lambda_1), \dots, v_p(\lambda_n)$.
\end{prop}
\begin{proof}
Assume without loss of generality that $v_p(\lambda_1) \geq \cdots
\geq v_p(\lambda_n)$.
Note that $v_p(\mu_1)+\cdots+v_p(\mu_n) = v_p(\lambda_1) + \cdots +
v_p(\lambda_n)$ since both are equal to the slopes of primitive
eigenvectors of $\wedge^n M$.
Note also that $\wedge^i M$ satisfies the conditions of 
Proposition~\ref{prop:minslope}
for all $i$, using the exterior products of the $\bw_j$
as the basis and the corresponding products of the $\mu_j$ as the
diagonal matrix entries. (More precisely, view the exterior products as being
partially ordered by sum of indices; any total ordering of the products
refining this partial order satisfies the conditions of the proposition.)
Since $\bv_{n-i+1} \wedge \cdots \wedge \bv_n$ is an eigenvector of
$\wedge^i M$ of slope $v_p(\lambda_{n-i+1}) + \cdots + v_p(\lambda_n)$,
by Proposition~\ref{prop:minslope}
this slope is greater than or equal to the smallest valuation of an
$i$-term product of the $\mu_j$, i.e., the sum of the $i$ smallest of
$v_p(\mu_1), \dots, v_p(\mu_n)$. This is precisely the
desired majorization.
\end{proof}
\begin{cor} \label{cor:newuniq}
Let $M$ be a $\sigma$-module over $\Gancon$.
If $\bv_1, \dots, \bv_n$ and $\bw_1, \dots, \bw_n$ are bases of $M$
such that $F\bv_i = \lambda_i \bv_i$ and $F\bw_i = \mu_i \bw_i$ for some
$\lambda_i, \mu_i \in \calO[\fp]$,
then the sequences $v_p(\lambda_1), \dots,
v_p(\lambda_n)$ and $v_p(\mu_1), \dots, v_p(\mu_n)$
are permutations of each other.
\end{cor}

Finally, we observe that the existence of an eigenvector of a specified
slope does not depend on what ring of scalars $\calO$ is used, so long as the
value group of $\calO$ contains the desired slope.
\begin{prop} \label{prop:valgp}
Let $M$ be a $\sigma$-module over $\Galgancon$. Suppose $\lambda 
\in \calO[\fp]$
occurs as the eigenvalue of an eigenvector of $M \otimes_{\calO} \calO'$
for some finite extension $\calO'$ of $\calO$. Then $\lambda$ occurs as
the eigenvalue of an eigenvector of $M$.  
\end{prop}
\begin{proof}
Since $k$ here is algebraically closed, we can choose
a basis $\mu_1, \dots, \mu_m$ of $\calO'$ over $\calO$ consisting
of elements fixed by $\sigma$. (Namely, let $\pi'$ be a uniformizer
of $\calO'$ fixed by $\sigma$, and take $\mu_i = (\pi')^{i-1}$.)
Given an eigenvector $\bv$ over $\calO'[\fp]$ with $F\bv = \lambda \bv$, we can
write $\bv = \mu_1 \bw_1 + \cdots + \mu_m \bw_m$ for a unique
choice of $\bw_1, \dots, \bw_m \in M$. Now
\[
0 = F\bv - \lambda \bv = \mu_1 (F\bw_1 - \lambda \bw_1) +
\cdots + \mu_m (F\bw_m - \bw_m).
\]
Since the representation $0 = \mu_1(0) + \cdots + \mu_m(0)$ is unique,
we must have $F\bw_i = \lambda \bw_i$ for $i=1, \dots, m$. Since 
$\bv$ is nonzero, at least one of the $\bw_i$ must be nonzero, and it
provides the desired eigenvector within $M$.
\end{proof}

\subsection{Existence of eigenvectors}

In this section, we prove that every $\sigma$-module over $\Galgancon$
has an eigenvector. 

\begin{prop} \label{prop:exist}
For every $\sigma$-module $M$ over $\Galgancon$,
there exist $\lambda \in \calO_0$ and $\bv \in M$, both nonzero,
such that $F\bv = \lambda \bv$.
\end{prop}
Note that once this assertion is established for a single $\lambda$,
it holds for all $\lambda \in \calO$ of sufficiently high valuation by
Proposition~\ref{prop:ansigeq}(b).

\begin{proof}
Let $v$ denote the valuation on $k((t))^{\alg}$ normalized so that 
$v(t) = 1$.
Let $\be_1, \dots, \be_n$ be a basis for $M$, and
suppose $F\be_i = \sum_j A_{ij} \be_j$.
Choose $r>0$ so that the entries of $A_{ij}$ all lie in $\Galg_{\an,r}$,
and let $c$ be an integer less than $\min\{w_r(A), w_r((A^{-1})^{\sigma^{-1}})\}$. For $0 < s \leq r$, we define the valuations $w_s$ on $M$ in terms of
the basis $\be_1, \dots, \be_n$. That is, $w_s(\sum_i c_i \be_i)
= \min_i \{w_s(c_i)\}$.

Notice that for $\lambda \in \calO_0$ and $u$ a strong semiunit,
\begin{align*}
v_0(u) \geq \frac{-c + v_p(\lambda)}{(q-1)r} 
&\iff
v_0(u)r \leq -v_p(\lambda)+v_0(u)qr+c \\
&\implies w_r(u \be_i) < w_r(\lambda^{-1} F(u \be_i)) \\
v_0(u) \leq \frac{qc+qv_p(\lambda)}{(q-1)r}
&\iff
v_0(u)r \leq v_p(\lambda)+v_0(u)r/q+c \\
 &\implies w_r(u \be_i) < w_r(\lambda F^{-1}(u \be_i)).
\end{align*}
Choose $\lambda \in \calO_0$ of large enough valuation so that
$-c + v_p(\lambda) < qc + qv_p(\lambda)$, and let $d$ be a rational
number such that $d(q-1)r \in (-c+v_p(\lambda), qc+qv_p(\lambda))$.

Define functions $a,b,f: M \to M$ as follows.
Given $\bw \in M$, write
$\bw = \sum_{i=1}^n z_i \be_i$, let $z_i = \sum_m \pi^m u_{i,m}$
be a strong semiunit decomposition for each $i$,
let $x_i$ be the sum of $\pi^m u_{i,m}$
over all $m$ such that $v_0(u_{i,m}) < d$, and put $y_i = z_i - x_i$.
Put $a(\bw) = \sum_{i=1}^n x_i \be_i$,
$b(\bw) = \sum_{i=1}^n y_i \be_i$, and
\[
f(\bw) = \lambda^{-1} b(\bw) - F^{-1} a(\bw).
\]
(Note: the definitions of $a,b,f$ depend on the choices of semiunit
decompositions above, but this does not cause any trouble.)
From the inequalities tabulated above, we have
\[
w_r(\lambda F^{-1} a(\bw)) \geq w_r(a(\bw)) + \epsilon,
\qquad
w_r(\lambda^{-1} F b(\bw)) \geq w_r(b(\bw)) + \epsilon
\]
for some $\epsilon > 0$.
Therefore
\begin{align*}
w_r(f(\bw)) &= w_r(\lambda^{-1} b(\bw) - F^{-1} a(\bw)) \\
&\geq w_r(\lambda^{-1} \bw) \\
w_r(F(f(\bw)) - \lambda f(\bw) + \bw) &=
w_r(F \lambda^{-1} b(\bw) - a(\bw) - b(\bw) - \lambda F^{-1} a(\bw) + \bw) \\
&= w_r(\lambda^{-1} F b(\bw) + \lambda F^{-1} a(\bw)) \\
&\geq w_r(\bw) + \epsilon
\end{align*}
for all nonzero $\bw \in M$.

Now define sequences $\{\bv_l\}_{l=0}^\infty$ and $\{\bw_l\}_{l=0}^\infty$
as follows. First choose $T \in k((t))^{\alg}$ of valuation $d$, and set
\[
\bv_0 = \lambda^{-1} [T] \be_1 + [T^{1/q}] F^{-1}\be_1,
\]
where the brackets again denote Teichm\"uller lifts.
Then define $\bv_l$ and $\bw_l$ recursively by the formulas
\[
\bw_l = F \bv_l - \lambda \bv_l, \qquad
\bv_{l+1} = \bv_l + f(\bw_l).
\]
For each $l$, $\bv_l$ is defined over $\Galg_{\an,rq}$
and $\bw_l$ is defined over $\Galg_{\an,r}$.
By the final remark of the previous paragraph, we have
\[
w_r(\bv_{l+1} - \bv_l) = w_r(f(\bw_l)) \geq w_r(\lambda^{-1} \bw_l)
\]
and
\begin{align*}
w_r(\bw_l) &= w_r(F\bv_l - \lambda \bv_l) \\
&= w_r(F \bv_{l-1} + F f(\bw_{l-1}) - \lambda \bv_{l-1} - 
\lambda f(\bw_{l-1})) \\
&= w_r(F f(\bw_{l-1}) - \lambda f(\bw_{l-1}) + \bw_{l-1}) \\
&\geq w_r(\bw_{l-1}) + \epsilon.
\end{align*}
Thus $w_r(\bw_l)$ is a
strictly increasing function of $l$
that tends to $\infty$, and $w_r(\bv_{l+1} - \bv_l)$ also tends to
$\infty$ with $l$.

We claim that in the Fr\'echet topology, $\bw_l$ converges to $0$ and so
$\bv_l$ converges to a limit $\bv$, from which it follows
that $F\bv - \lambda \bv = \lim_{l \to \infty} \bw_l = 0$.
We first show that $w_s(\lambda F^{-1} a(\bw_l)) \to \infty$
as $l \to \infty$ for $0 < s \leq qr$.

Let $a(\bw_l) = \sum_{i,m} \pi^m a_{l,i,m} \be_i$
be a strong semiunit decomposition, in which we must have
$v_0(a_{l,i,m}) < d$ whenever $a_{l,i,m} \neq 0$. Then
\begin{align*}
w_s(\lambda F^{-1}a(\bw_l)) &\geq w_s(\lambda (A^{-1})^{\sigma^{-1}})
+ w_s(a(\bw_l)^{\sigma^{-1}}) \\
&= w_s(\lambda (A^{-1})^{\sigma^{-1}})
+ w_{s/q}(a(\bw_l)) \\
&= w_s(\lambda (A^{-1})^{\sigma^{-1}}) +
\min_{i,m} \{ m v_p(\pi) + (s/q) v_0(a_{l,i,m}) \} \\
&\geq w_s(\lambda (A^{-1})^{\sigma^{-1}}) +
\min_{i,m} \{ m v_p(\pi) + r v_0(a_{l,i,m})\} + \min_{i,m} \{ (-r+s/q) v_0(a_{l,i,m})
\} \\
& > w_s(\lambda (A^{-1})^{\sigma^{-1}}) +
w_r(a(\bw_l)) - (r-s/q)d.
\end{align*}
In particular, $w_s(\lambda F^{-1} a(\bw_l)) \to \infty$ as $l \to \infty$.

We next show that $w_s(\lambda^{-1} F b(\bw_l)) \to \infty$
as $l \to \infty$ for $0 < s \leq r$.
Let $b(\bw_l) = \sum_{i,m} \pi^m b_{l,i,m} \be_i$
be a strong semiunit decomposition, necessarily with
$v_0(b_{l,i,m}) \geq d$ whenever $b_{l,i,m} \neq 0$. Then
\begin{align*}
w_s(\lambda^{-1} F b(\bw_l)) &\geq
w_s(\lambda^{-1} A) + w_s(b(\bw_l)^\sigma) \\
&= w_s(\lambda^{-1} A) + w_{sq}(b(\bw_l)) \\ 
&= w_s(\lambda^{-1} A) + \min_{i,m} \{mv_p(\pi) + sq v_0(b_{l,i,m})\}. 
\end{align*}
Choose $e>0$ large enough so that $s(q-1)e + w_s(\lambda^{-1} A) > 0$.
If $v_0(b_{l,i,m}) < e$, then
\begin{align*}
m v_p(\pi) + sq v_0(b_{l,i,m}) &= m v_p(\pi) + r v_0(b_{l,i,m}) + (sq-r) v_0(b_{l,i,m}) \\
&\geq w_r(b(\bw_l)) +h,
\end{align*}
where $h=(sq-r)d$ if $sq-r \geq 0$ and $h=(sq-r)e$ if $sq-r < 0$.
If $v_0(b_{l,i,m}) \geq e$, then
\begin{align*}
m v_p(\pi) + sq v_0(b_{l,i,m}) &= m v_p(\pi)
+ s v_0(b_{l,i,m}) + s(q-1) v_0(b_{l,i,m}) \\
&\geq w_s(b(\bw_l)) + s(q-1)e.
\end{align*}

Suppose $\liminf_{l \to \infty} w_s(b(\bw_l)) < L$
for some $L < \infty$.
For $l$ sufficiently large, we have $w_s(\lambda F^{-1} a(\bw_l))
\geq L$ and $w_r(b(\bw_l)) \geq L -h - w_s(\lambda^{-1} A)$;
 by the previous
paragraph, this implies
\begin{align*}
w_s(b(\bw_{l+1})) &\geq w_s(\bw_{l+1}) \\
&= w_s(F\bv_{l+1} - \lambda \bv_{l+1}) \\
&= w_s(F\bv_l + Ff(\bw_l) - \lambda \bv_l - \lambda f(\bw_l)) \\
&= w_s(\bw_l + F f(\bw_l) - \lambda f(\bw_l)) \\
&= w_s(\lambda^{-1} F b(\bw_l) + \lambda F^{-1} a(\bw_l)) \\
&\geq \min\{w_s(\lambda^{-1}A) + w_s(b(\bw_l)) + s(q-1)e, L\}.
\end{align*}
We first deduce from this inequality that $w_s(b(\bw_l))$ is bounded
below: pick any $l$, choose $C<L$ such that $w_s(b(\bw_l))>C$,
then note that $w_s(b(\bw_{l+1})) \geq \min\{L, C+w_s(\lambda^{-1} A) +
s(q-1)e\} > C$.
If we put $M = \liminf w_s(b(\bw_l))$, we thus have
$-\infty < M < L$.
However, in the inequality above, the limit inferior of the left side is $M$,
while the limit inferior of the smaller right
side is $\min\{L, M+w_s(\lambda^{-1} A) + s(q-1)e\} > M$.
This contradiction shows
that no $L$ can exist as above, and so
$w_s(\lambda^{-1}  F b(\bw_l)) \to \infty$.

From $w_s(\lambda F^{-1} a(\bw_l)) \to \infty$ for $0 < s \leq qr$,
and $w_s(\lambda^{-1} F b(\bw_l)) \to \infty$ for $0 < s \leq r$,
we conclude that $w_s(a(\bw_l)) \to \infty$ and 
$w_s(b(\bw_l)) \to \infty$ for $0 < s \leq r$. Thus
$\bw_l$ converges to 0 in the Fr\'echet topology,
and $\bv_l$ converges to
a limit $\bv$ satisfying $F\bv = \lambda \bv$.

Finally, we check that $\bv \neq 0$. First note that
$w_r(\lambda^{-1} [T] \be_1) = dr - v_p(\lambda)$,
while
\begin{align*}
w_r(\bv_0 - \lambda^{-1}[T] \be_1) &= w_r([T^{1/q}]F^{-1} \be_1) \\
&\geq dr/q + c \\
&> dr - v_p(\lambda)
\end{align*}
by our choice of $d$.
Therefore $w_r(\bv_0) = dr - v_p(\lambda)$.
On the other hand,
\begin{align*}
w_r(\lambda^{-1} \bw_0) &=
 w_r(\lambda^{-1} F \bv_0 - \bv_0) \\
&= w_r(\lambda^{-2} F [T] \be_1 + \lambda^{-1} [T] \be_1
- \lambda^{-1} [T] \be_1 - [T^{1/q}] F^{-1} \be_1) \\
&= w_r(\lambda^{-2} F [T] \be_1 - [T^{1/q}] F^{-1} \be_1) \\
&\geq \min\{ r d q + c - 2 v_p(\lambda), r d/q + c\}.
\end{align*}
We have just checked that the second term in braces is greater than
$dr - v_p(\lambda) = w_r(\bv_0)$. As for the first term,
\[
rdq + c - 2 v_p(\lambda) - (dr - v_p(\lambda))
= dr(q-1) + c - v_p(\lambda)
\]
is positive, again by the choice of $d$.
Therefore $w_r(\lambda^{-1} \bw_0) > w_r(\bv_0)$. 

Since we showed
earlier that $w_r(\bw_l)$ is a strictly increasing
function of $l$, we have $w_r(\lambda^{-1} \bw_l) \geq w_r(\lambda^{-1}
\bw_0)$ for 
$l \geq 0$. We also showed earlier that $w_r(\bv_{l+1} - 
\bv_l) \geq w_r(\lambda^{-1} \bw_l)$ for $l \geq 0$. Thus
$w_r(\bv_{l+1} - \bv_l) \geq w_r(\lambda^{-1} \bw_0)$ for each $l$, and so
$w_r(\bv_l - \bv_0) \geq w_r(\lambda^{-1} \bw_0)$.
It follows that $w_r(\bv - \bv_0) \geq w_r(\lambda^{-1} \bw_0) > w_r(\bv_0)$; in particular,
$\bv \neq 0$, so $\lambda$
and $\bv$ satisfy the desired conditions.
\end{proof}
\begin{cor} \label{cor:specfilt}
Every $\sigma$-module $M$ over $\Galgancon$ admits a basis $\bv_1,
\dots, \bv_n$ such that $\bv_i$ is an eigenvector in $M/\Span(\bv_1,
\dots, \bv_{i-1})$ for $i=1, \dots, n$.
\end{cor}
\begin{proof}
By the proposition and Lemma~\ref{lem:prim}, every $\sigma$-module
over $\Galgancon$ contains a primitive eigenvector. The corollary
now follows by induction on the rank of $M$.
\end{proof}
\begin{cor} \label{cor:bdbelow}
The set of slopes of eigenvectors of $M$, over all finite extensions
of $\calO$, is bounded below.
\end{cor}
\begin{proof}
Combine the previous corollary with Proposition~\ref{prop:minslope}.
\end{proof}

\subsection{Raising the Newton polygon}

In the previous section, we produced within any $\sigma$-module
over $\Galgancon$ a basis on which $F$ acts by a triangular matrix.
By Proposition~\ref{prop:basis}, if there is a basis of eigenvectors,
the valuations
of the diagonal entries of this matrix majorize the slopes of 
the eigenvectors.
Thus to produce a basis of eigenvectors,
we need to ``raise the Newton polygon'', i.e., find eigenvectors
whose eigenvalues have smaller slopes than the ones we started with.
In this section, we carry this process out by direct computation
in an important special case; the general process, using this case
in some basic steps, will follow in the next section.

By a \emph{Puiseux polynomial} over a field $K$, we shall mean a
formal expression of the form
\[
P(z) = \sum_{i \in I} c_i z^i
\]
where $I$ is a finite set of nonnegative rationals
and $c_i \in K$ for each $i \in I$.
If $K$ has a valuation $v_K$, we define the Newton polygon of a
Puiseux polynomial, by analogy with the definition for an ordinary
polynomial, as the lower convex hull of the set of points
$(-i, v_K(c_i))$. In fact, for some integer $n$,
$P(z^n)$ is an ordinary polynomial; by comparing the Newton polygons
of $P(z)$ and $P(z^n)$, and using the usual theory of Newton polygons
of polynomials
over fields complete with respect to a valuation, 
we obtain the following result.
\begin{lemma} \label{lem:puis}
Let $P(z)$ be a Puiseux polynomial over the $t$-adic completion
of $k((t))^{\alg}$.
Then $P$ has a root of valuation $l$ if and only if the 
Newton polygon of $P$ has a segment of slope $l$.
\end{lemma}

For $x \in \Galgancon$ a strong
semiunit, we refer to $v_0(x)$ as the \emph{valuation} of $x$.
\begin{lemma} \label{lem:eq}
Let $n$ be a positive integer, and let
$x = \sum_{i=0}^{n} u_i \pi^i$ for some strong 
semiunits $u_i \in \Galgancon$ of negative (or infinite)
valuation, not all zero. Then the system of equations
\begin{equation} \label{eq:system}
a^\sigma = \pi a, \qquad
\pi b^{\sigma^n} = b - a x
\end{equation}
has a solution with $a, b \in \Galgancon$ not both zero.
\end{lemma}
\begin{proof}
For $i \in \{0, \dots, n\}$ for which $u_i \neq 0$, 
$l \in \ZZ$ and $m \in \RR^+$, put
\[
f(i, l, m) = (v_0(u_i) + m q^{-l})q^{-n(i+l)},
\]
Note that for fixed $i$ and $m$, $f(i,l,m)$ approaches 0 from below
as $l \to +\infty$, and tends to $+\infty$ as $l \to -\infty$.
Thus the minimum $h(m) = \min_{i,l} \{f(i,l,m)\}$ is well-defined.
Observe that the map $h: \RR^+ \to \RR$ is continuous and
piecewise linear with everywhere positive slope, and
$h(qm) = q^{-n} h(m)$ because $f(i,l+1,qm) = q^{-n} f(i,l,m)$.
Since $f(i,l,m)$ takes negative values
for fixed $i,l$ and small $m$, $h(m) < 0$ for some $m$,
implying $h(q^j m)<0$ for all $j \in \ZZ$, so $h$ takes only
negative values. We conclude that $h$ is a continuous
increasing bijection of $\RR^+$ onto $\RR^-$.

Pick $t \in \RR^+$ at which $h$ changes slope,
let $S$ be the finite set of ordered pairs $(i,l)$ for which
$f(i,l,t) < q^{-n} h(t)$, and let $T$ be the
set of ordered pairs $(i,l)$ for
which $f(i,l,t) < 0$; then $T$ is infinite (and contains $S$),
but the values of $l$ for pairs $(i,l) \in T$ are bounded below.
For each pair $(i,l) \in T$, put $s(i,l) = \lfloor \log_{q^n} 
(h(t)/f(i,l,t))
\rfloor$. This function has the following properties:
\begin{enumerate}
\item[(a)] $s(i,l) \geq 0$ for all $(i,l) \in T$;
\item[(b)] $f(i,l,t) q^{ns(i,l)} \in [h(t), q^{-n}h(t))$ for all
$(i,l) \in T$;
\item[(c)] $(i,l) \in S$ if and only if $(i,l) \in T$ and
$s(i,l) = 0$;
\item[(d)] for any $e>0$, there are only finitely many pairs
$(i,l) \in T$ such that $s(i,l) \leq e$.
\end{enumerate}

For $c \in \RR$, let $U_c$ be the set of $z \in \Galgancon$ such that
$v_m(z) = \infty$ for $m < 0$ and $v_m(z) \geq c$ for $m \geq 0$.
Then the function
\[
r(z) = \sum_{(i,l) \in T} \pi^{s(i,l)} u_i^{\sigma^{-ni-nl+n s(i,l)}} 
z^{\sigma^{-ni-(n+1)l+n s(i,l)}}
\]
is well-defined (by (d) above, the series is $\pi$-adically convergent)
and carries $U_t$ into $U_{h(t)}$ because for $z \in U_t$ and $m \geq 0$,
\begin{align*}
v_m\left(u_i^{\sigma^{-ni-nl+n s(i,l)}}
z^{\sigma^{-ni-(n+1)l+n s(i,l)}} \right)
&\geq q^{-ni-nl+ns(i,l)} v_0(u_i) + \min_{j} 
  \{ q^{-ni-(n+1)l+n(s,i,l)} v_j(z) \} \\
&\geq q^{ns(i,l)} q^{-n(i+l)} (v_0(u_i) + t q^{-l}) \\
&= q^{ns(i,l)} f(i,l,t) \\
&\geq h(t).
\end{align*}
The reduction of $r(z)$ modulo $\pi$ is
congruent to a finite sum over pairs $(i,l) \in S$, so it is a
Puiseux polynomial in the reduction of $z$.
Since $s(i,l) = 0$ for all $(i,l) \in S$ and the values
$-ni-(n+1)l$ are all distinct (because $i$ only runs over
$\{0, \dots, n\}$),
we get a distinct monomial modulo $\pi$ for each pair $(i,l) \in S$.

We now consider the Newton polygon of the Puiseux polynomial given
by the reduction of $r(z)-w$, for $w \in U_{h(t)}$.
It is the convex hull of the points
$(-q^{-ni-(n+1)l}, v_0(u_i) q^{-ni-nl})$ for each $(i,l) \in S$,
together with $(0, v_0(w))$.
The line $y = tx + h(t)$
either passes through or lies below the point corresponding to
$(i,l)$, 
depending on whether $f(i,l,t)$ is equal to or strictly
greater than $h(t)$.
Moreover, $(0, v_0(w))$ lies on or above the
line because $v_0(w) \geq h(t)$.
Since $h$ changes slope at $t$, there must be at
least two points on the line; therefore the Newton polygon has
a segment of slope $t$. By Lemma~\ref{lem:puis}, the Puiseux polynomial has a
root
of valuation $t$. In other words,
there exists $z \in U_t$ with $v_0(z) = t$
such that $r(z) \equiv w \pmod{\pi}$.

As a consequence of the above reasoning, we see that
the image of $U_t$ is dense in $U_{f(t)}$ with respect to the
$\pi$-adic topology. Since $U_t$ is complete,
$U_t$ must surject onto $U_{f(t)}$. Moreover, we can take $w=0$ and obtain
$z_0 \in U_t$ with $v_0(z_0) = t$ such that $r(z_0) \equiv 0 \pmod{\pi}$;
in particular, $z_0$ is nonzero modulo $\pi$.
We may then obtain $z_1 \in U_t$ such that $r(z_1) = r(z_0)/\pi$.
Put $z = z_0 - \pi z_1$; then $z \not\equiv 0 \pmod{\pi}$ and so
is nonzero, but $r(z) = 0$.

Now set $a = \sum_{l=-\infty}^\infty \pi^{l} z^{\sigma^{-l}}$;
the sum converges in $\Galgancon$
because for $s>0$, $w_s(\pi^l z^{\sigma^{-l}})
\geq l v_p(\pi) + r q^{-l} t$ and the latter tends to $\infty$ as
$l \to \pm \infty$ (because $t>0$).
Then
\begin{align*}
ax &= \sum_{i=0}^{n} \sum_{l=-\infty}^\infty \pi^{i+l} u_i z^{\sigma^{-l}} \\
&=  \sum_{(i,l) \in T} \pi^{i+l} u_i z^{\sigma^{-l}} +
\sum_{(i,l) \notin T} \pi^{i+l} u_i z^{\sigma^{-l}}.
\end{align*}
Let $A$ and $B$ denote the two sums in the last line; then $v_m(B)
\geq 0$
for all $m$, so by Proposition~\ref{prop:ansigeq}(d) (with $\sigma$
replaced by $\sigma^n$),
$B$ can be written as $\pi b_1^{\sigma^n} - b_1$ for some $b_1 \in \Galgancon$.
On the other hand, we claim that $A$ can be rewritten
as $r(z) + \pi b_2^{\sigma^n} - b_2$ for
\begin{align*}
b_2 &= \sum_{(i,l) \in T} \sum_{j=1}^{i+l-s(i,l)} 
\pi^{i+l-j} u_i^{\sigma^{-nj}} z^{\sigma^{-l-nj}} \\
&= \sum_{(i,l) \in T} \sum_{k=0}^{i+l-s(i,l)-1} 
\pi^{k+s(i,l)} \left(
u_i^{\sigma^{-ni-nl+ns(i,l)}} z^{\sigma^{-ni-(n+1)l+ns(i,l)}} \right)^{\sigma^{nk}}
\end{align*}
(via the substitution $k = i+l-s(i,l) - j$);
we must check that this series converges $\pi$-adically and that its limit
is overconvergent. Note that as $l \to +\infty$ for $i$ fixed, 
$f(i,l,m)$ is asymptotic to $v_0(u_i) q^{-n(i+l)}$.
Therefore
$i+l-s(i,l)$ is bounded, so the possible values of $k$ are uniformly
bounded over all pairs $(i,l) \in T$. This implies on one hand that
the series converges $\pi$-adically (since $l$ is bounded below over
pairs $(i,l) \in T$ and $s(i,l) \to \infty$ as
$l \to \infty$), and on the other hand that
$v_m(b_2)$ is bounded below uniformly in $m$ (since the quantity
in parentheses in the second sum belongs to $U_{h(t)}$),
so $b_2 \in \Galgancon$.

Having shown that the series defining $b_2$ converges, we can
now verify that $b_2 - \pi b_2^{\sigma^n} = r(z) - A$:
the quantity on the left is the sum over pairs $(i,l) \in T$ of
a sum over $k$ which telescopes, leaving the term
$k=0$ minus the term $k=i+l-s(i,l)$, or
\[
\pi^{s(i,l)}
u_i^{\sigma^{-ni-nl+ns(i,l)}} z^{\sigma^{-ni-(n+1)l+ns(i,l)}} 
- \pi^{i+l}
u_i z^{\sigma^{-l}},
\]
which when summed over pairs $(i,l) \in T$ yields
$r(z) - A$.

Since $r(z) = 0$ by construction, we have
$\pi b^{\sigma^n} = b - a x$ for $b = -(b_1+b_2)$. Thus $(a,b)$ constitute
a solution of (\ref{eq:system}), as desired.
\end{proof}

We apply the previous construction to study the system of equations
\begin{equation} \label{eq:system2}
a^\sigma = \pi a, \qquad \pi b^{\sigma^n} = b - ac,
\end{equation}
where $c \in \Galgancon$ is given.
Notice that replacing $c$ by $c + \pi^{n+1} y^{\sigma^n} - y$
does not alter whether (\ref{eq:system2}) has a solution: for any $a$
such that $a^\sigma = \pi a$, if $\pi b^{\sigma^n} = b - ac$, then
\[
\pi (b - a y)^{\sigma^n} = (b - a y) 
- a(c + \pi^{n+1} y^{\sigma^n} - y).
\]

We begin by analyzing \eqref{eq:system2} in a restricted case.
\begin{lemma} \label{lem:solve1}
For any positive integer $n$ and any
$c \in \Galgancon$ such that $v_m(c) \geq -1$ for all $m$ and
$v_m(c) = \infty$ for some $m$,
there exist $a,b \in \Galgancon$ not both zero satisfying
\eqref{eq:system2}.
\end{lemma}
\begin{proof}
By multiplying $c$ by a power of $\pi$, we may reduce to the case
where $v_m(c) = \infty$ for $m < 0$.
Define the sequences $c_0, c_1, \dots$ and $d_0, d_1, \dots$ as far
as is possible by the following iteration. First put $c_0 = c$ and
$d_0 = 0$. Given $c_i$, if $v_0(c_i) < -1/q^n$, stop. Otherwise,
let $d_i$ be a strong semiunit congruent to $c_i$ modulo $\pi$ and
put $c_{i+1} = (c_i + \pi^{n+1} d_i^{\sigma^n} - d_i)/\pi$. Note that
$v_m(c_i) \geq -1$ and $v_m(d_i) \geq -1/q^n$ for all $m \geq 0$ and all $i$.

If the iteration never terminates, then we have
$c + \pi^{n+1} d^{\sigma^n} - d = 0$ for $d = \sum_{i=0}^\infty d_i \pi^i$.
In this case, apply Proposition~\ref{prop:ansigeq}(b) to produce $a$
nonzero such that $a^\sigma = \pi a$ and set
$b = a d$ to obtain a solution to \eqref{eq:system2}.

If the iteration terminates at $c_l$, 
set $d = \sum_{i=0}^{l-1} d_i \pi^i$, so that
$\pi^l c_l = c + \pi^{n+1} d^{\sigma^n} - d$.
Let $\sum_{j=0}^\infty u_j \pi^j$ be a strong semiunit
decomposition of $c_l$, necessarily having $v_0(u_0) < -1/q^n$.
Put $e = \sum_{j=n+1}^\infty u_j^{\sigma^{-n}} \pi^{j-n-1}$ and set
$x = c_l - \pi^{n+1} e^{\sigma^n} + e$.
Then
\[
x = \sum_{j=1}^{n} \pi^j u_j + \left( u_0 + \sum_{j=n+1}^\infty
u_j^{\sigma^{-n}} \pi^{j-n-1} \right),
\]
and the quantity in parentheses is a strong semiunit of the same valuation
as $u_0$, since $v_0(u_0) < -1/q^n \leq v_0(u_j^{\sigma^{-n}})$ 
for all $j$.
Thus $x$ satisfies the condition of Lemma~\ref{lem:eq}, so there exist
$a', b' \in \Galgancon$ not both zero so that
\[
(a')^\sigma = \pi a', \qquad \pi (b')^{\sigma^n} = b' - a' x.
\]
We obtain a solution of \eqref{eq:system2} by setting
$a = a'$, $b = a'd - \pi^l a' e + \pi^l b'$.
\end{proof}

We now analyze \eqref{eq:system2} in general
by reducing to the special case treated above.
\begin{lemma} \label{lem:solve2}
For any positive integer $n$ and any
$c \in \Galgancon$,
there exist $a,b \in \Galgancon$ not both zero such that
\eqref{eq:system2} holds.
\end{lemma}
\begin{proof}
Let $\sum_i u_i \pi^i$ be a strong semiunit decomposition of $c$,
and let $N$ be the smallest integer such that $v_{0}(u_N) < 0$,
or $\infty$ if there is no such integer.
By Proposition~\ref{prop:ansigeq}(d), there exists $y \in \Galgancon$
such that $\pi^{n+1} y^{\sigma^n} - y + \sum_{i =-\infty}^{N-1} u_i \pi^i
= 0$.

If $N = \infty$, then in fact $\pi^{n+1} y^{\sigma^n} - y + c = 0$,
so we obtain a solution of \eqref{eq:system2} by choosing $a$ nonzero
with $a^\sigma = \pi a$ via Proposition~\ref{prop:ansigeq}(b),
and setting $b = a y$. So suppose hereafter that $N < \infty$.

For each $i \geq N$ for which $u_i \neq 0$, set
$t_i = \lceil \log_{q^n} (-v_0(u_i)) \rceil$, so that
$-1 \leq v_0(u_i^{\sigma^{-nt_i}}) < -1/q^{n}$ for all such $i$.
Then the sum
\[
z = \sum_{i=N}^\infty \sum_{j=1}^{t_i} u_i^{\sigma^{-nj}}
\pi^{i - (n+1)j}
\]
is $\pi$-adically convergent: $-v_0(u_i)$ grows at most
linearly in $i$, so $t_i$ grows at most logarithmically and
$i - (n+1)t_i \to \infty$ as $i \to \infty$.
Moreover, $t_i$ is bounded below, so 
$v_0(u_i^{\sigma^{-nj}} \pi^{i - (n+1)j})$ is as well; thus
the sum $z$ is in $\Galgancon$.

Put $c' = c + \pi^{n+1}(y-z)^{\sigma^n} - (y-z)$; then
\begin{align*}
c' &=
\sum_{i=N}^\infty \left( u_i \pi^i + \sum_{j=1}^{t_i}
u_i^{\sigma^{-nj}} \pi^{i-(n+1)j}
- \sum_{j=1}^{t_i} u_i^{\sigma^{-n(j-1)}} \pi^{i-(n+1)(j-1)}
\right) \\
&= \sum_{i=N}^\infty u_i^{\sigma^{-nt_i}} \pi^{i-(n+1)t_i},
\end{align*}
so that $v_m(c') \geq -1$ for all $m$. By Lemma~\ref{lem:solve1},
there exist $a', b' \in \Galgancon$ not both zero such that
\[
(a')^\sigma = \pi a', \qquad \pi (b')^{\sigma^n} = b' - a'c';
\]
we obtain a solution of \eqref{eq:system2} by setting
$a = a', b = b' + a'(y-z)$.
\end{proof}

We now prove our basic result on raising the Newton polygon, i.e., reducing
the slope of an eigenvector. 
\begin{prop} \label{prop:direct}
Let $m$ and $n$ be positive integers, and let
$M$ be a $\sigma$-module over $\Galgancon$ admitting a basis
$\bv_1, \dots, \bv_n, \bw$ such that for some $c_i \in \Galgancon$,
\begin{align*}
F\bv_i &= \bv_{i+1} \qquad (i=1, \dots, n-1) \\
F\bv_n &= \pi \bv_1 \\
F\bw &= \pi^{-m} \bw + c_1 \bv_1 + \cdots + c_n \bv_n. 
\end{align*}
Then there exists $\by \in M$ such that $F\by = \by$.
\end{prop}
This will ultimately be a special case
of our main results; what
makes this case directly tractable is that if $\Span(\bv_1,
\dots, \bv_n)$
does not admit an $F$-stable complement in $M$ (i.e.,
is not a direct summand of $M$ in the category of $\sigma$-modules),
then the map $\by \mapsto
F\by -
\by$ is actually surjective, as predicted by the expected behavior
of the special Newton polygon.
\begin{proof}
Suppose $\by = d \bw + b_1 \bv_1+\cdots + b_n\bv_n$
satisfies $F\by = \by$, or in other words
\[
d \bw + \sum_{i=1}^n b_i \bv_i =
\pi^{-m} d^\sigma \bw + \sum_{i=1}^n d^\sigma c_i \bv_i
+ \sum_{i=1}^{n-1} b_i^\sigma \bv_{i+1} + \pi b_n^\sigma \bv_1.
\]
Comparing coefficients in this equation, we have
$b_i^\sigma = b_{i+1} - d^\sigma c_{i+1}$ for $i=1, \dots, n-1$,
as well as
$\pi b_n^\sigma = b_1 - d^\sigma c_1$ and
$d^\sigma = \pi^m d$.
If we use the first $n$ relations to
eliminate $b_2, \dots, b_n$, we get 
\begin{align*}
b_1^{\sigma^n} &= b_2^{\sigma^{n-1}} - 
d^{\sigma^n} c_2^{\sigma^{n-1}} \\
&= b_3^{\sigma^{n-2}} - d^{\sigma^{n-1}} c_3^{\sigma^{n-2}}
- d^{\sigma^n} c_2^{\sigma^{n-1}} \\
&\qquad \vdots \\
&= b_n^\sigma - d^{\sigma^2} c_{n}^\sigma - 
\cdots - d^{\sigma^n} c_2^{\sigma^{n-1}} \\
&= \pi^{-1} b_1 - d(\pi^{m-1} c_1 + \pi^{2m} c_{n}^\sigma + 
\pi^{3m} c_{n-1}^{\sigma^2} + \cdots + \pi^{nm} c_2^{\sigma^{n-1}}).
\end{align*}
Let $c'$ be the quantity in parentheses in the last line. We have shown
that if $F\by = \by$ has a nonzero solution, then the system
of equations
\begin{equation} \label{eq:system3}
d^\sigma = \pi^m d, \qquad
\pi b_1^{\sigma^n} =
b_1 - \pi c'd
\end{equation}
has a solution with $b_1, d$ not both zero.
Conversely, from any nonzero solution of \eqref{eq:system3}
we may construct a nonzero $\by \in M$ such that $F\by = \by$, by
using the relations $b_i^\sigma = b_{i+1} - d^\sigma c_{i+1}$
to successively define $b_2, \dots, b_n$.

By Proposition~\ref{prop:ansigeq}(b), we can find $e \in \Galgancon$
nonzero such that $e^\sigma = \pi^{m-1} e$; we will construct
a solution of (\ref{eq:system3}) with $d = ae$
for some $a$ such that $a^\sigma = \pi a$.
Namely, put $c = \pi c' e$, and apply Lemma~\ref{lem:solve2}
to find $a,b \in \Galgancon$, not both zero, such that
\[
a^\sigma = \pi a, \qquad
\pi b^{\sigma^n} = b - a c.
\]
Then $b_1 = b$ and $d = ae$ constitute a nonzero solution of
\eqref{eq:system2}; as noted above, this implies that there exists
$\by \in M$ nonzero with $F\by = \by$, as desired.
\end{proof}

\subsection{Construction of the special Newton polygon}

We now assemble the results of the previous sections into the following
theorem, the main result of this chapter.
\begin{theorem} \label{thm:decomp}
Let $M$ be a $\sigma$-module over $\Galgancon$. Then $M$
can be expressed as a direct sum of standard $\sigma$-submodules.
\end{theorem}
As the proof of this theorem is somewhat intricate, we break
off parts of the argument into separate lemmas. In these lemmas,
a ``suitable extension'' of $\calO[\fp]$ means one whose value group contains
whatever slope is desired to be the slope of an eigenvector.
By Proposition~\ref{prop:valgp}, proving the existence of an eigenvector
of prescribed slope over a single suitable extension implies the same
over any suitable extension.

\begin{lemma} \label{lem:rank1}
Let $M$ be a $\sigma$-module over $\Galgancon$ of rank $1$,
and suppose $F$ acts on some generator $\bv$ via $F\bv = c \bv$.
Then $M$ contains an eigenvector, and any eigenvector has slope
$v_p(c)$.
\end{lemma}
Note that $v_n(c) = \infty$ for some $n$ by Corollary~\ref{cor:units},
so that $v_p(c)$ makes sense.
\begin{proof}
The existence of an eigenvector of slope $v_p(c)$ follows from
Proposition~\ref{prop:rank1}. The uniqueness of the slope follows
from Corollary~\ref{cor:newuniq}.
\end{proof}
For $M$ of rank $1$,
 we call this unique slope the slope of $M$.
Note that if $0 \to L \to M \to N \to 0$ is an exact sequence of 
$\sigma$-modules and $L,M,N$ have ranks $l,m,n$, respectively,
then the slope of $\wedge^{m} M$ is the sum
of the slopes of $\wedge^l L$ and $\wedge^n N$.
(This assertion will be vastly generalized by Proposition~\ref{prop:operations}
later.)

\begin{lemma} \label{lem:raise1}
Let $M$ be a $\sigma$-module over $\Galgancon$ of rank $2$,
and let $d$ be the slope of $\wedge^2 M$. Then $M$ contains an eigenvector
of slope $d/2$ over a suitable extension of $\calO[\fp]$.
\end{lemma}
\begin{proof}
We may assume without loss of generality
that $d/2$ belongs to the value group of $\calO[\fp]$. Let $e$ be the
smallest integer such that $M$ contains an eigenvector of slope $e v_p(\pi)$.
(There is such an integer by Proposition~\ref{prop:exist},
and there is a smallest one by Corollary~\ref{cor:bdbelow}.)
By twisting, we may reduce to the case where $e = 1$. 

Put $m = 1-(d/v_p(\pi))$ and suppose by way of contradiction that $m > 0$. 
Choose a eigenvector $\bv$ with $F\bv = \pi \bv$, 
which is necessarily primitive by Lemma~\ref{lem:prim}; 
then by Lemma~\ref{lem:rank1} applied to
$M / \Span(\bv)$, we can
find $\bw$ such that $\bv, \bw$ form a basis of $M$ and
$F\bw = \pi^{-m} \bw + c\bv$ for some $c \in \Galgancon$. 
Now by Proposition~\ref{prop:direct}, $M$ contains an eigenvector
$\bv_1$ with $F\bv_1 = \bv_1$, contradicting the definition of $e$.

Hence $m \leq 0$, which implies $d \geq v_p(\pi)$. Since $d/2$ is also
a multiple of $v_p(\pi)$, we must have $d/2 \geq v_p(\pi)$;
by Proposition~\ref{prop:raise0}, $M$ contains an eigenvector of slope
$d/2$.
\end{proof}

\begin{lemma} \label{lem:raise2}
Let $M$ be a $\sigma$-module over $\Galgancon$ of rank $n$,
and let $d$ be the slope of $\wedge^n M$. Then $M$ contains eigenvectors
of all slopes greater than $d/n$ over suitable extensions of $\calO[\fp]$.
\end{lemma}
\begin{proof}
We proceed by induction on $n$. The case $n=1$ follows from
Lemma~\ref{lem:rank1}, and the case $n=2$ follows from
Lemma~\ref{lem:raise1}.
Suppose $n>2$ and that the lemma has been proved for all smaller
values of $n$.
Let $s$ be the greatest lower bound of the set of rational numbers that
occur as slopes of eigenvectors of $M$ (over suitable extensions of
$\calO[\fp]$). Again, the set is nonempty by Proposition~\ref{prop:exist}
and is bounded below by Corollary~\ref{cor:bdbelow}.

For each $\epsilon>0$ such that $s+\epsilon \in \QQ$, 
over a suitable extension of $\calO[\fp]$ 
there exist an eigenvector $\bv$ of $M$ of
slope $s+\epsilon$ and (by the induction hypothesis) an eigenvector
$\bw$ of $M/\Span(\bv)$ of slope at most
$s' = (d-s-\epsilon)/(n-1) + \epsilon$.
The preimage of $\Span(\bw)$ in $M$ has rank 2, so is covered
by the induction hypothesis; it thus contains, for any
$\delta>0$, an eigenvector of slope at most
\[
\frac{s+\epsilon}{2} + \frac{d-s+(n-2)\epsilon}{2(n-1)} + \delta
\]
over a suitable extension of $\calO[\fp]$.
Such an eigenvector is also an eigenvector of $M$, so its slope
is at least $s$. 
Letting $\epsilon$ and $\delta$ go to 0 in the resulting inequality yields
\[
\frac{s}{2} + \frac{d-s}{2(n-1)} \geq s,
\]
which simplifies to $s \leq d/n$, as desired.
\end{proof}

\begin{lemma} \label{lem:goodeigen}
Let $M$ be a $\sigma$-module over $\Galgancon$ of rank $n$,
and let $d$ be the slope of $\wedge^n M$. Then $M$ contains an eigenvector
of slope $d/n$ over a suitable extension of $\calO[\fp]$.
\end{lemma}
\begin{proof}
We proceed by induction on $n$; again, the case $n=1$ follows from
Lemma~\ref{lem:rank1} and the case $n=2$ follows from
Lemma~\ref{lem:raise1}.
Without loss of generality, we may assume the value group of $\calO$
contains $d/n$, and then that $d=0$.

By Lemma~\ref{lem:raise2},
there exists an eigenvector $\bv$ of $M$ of slope $v_p(\pi)/(n-1)$
over $\calO[\pi^{1/(n-1)}]$; we may as well assume $F\bv = \pi^{1/(n-1)} \bv$.
Let $N$ be the saturated span
of $\bv$ and its conjugates over $\calO[\fp]$; let $m$ be the rank
of $N$ and $s$ the slope of $\wedge^m N$.
Then $m \leq n-1$ and $s \leq m v_p(\pi)/(n-1)$. If $m<n-1$, then
$0 < mv_p(\pi)/(n-1) < v_p(\pi)$, so $s \leq 0$ and
the induction hypothesis implies that $N$
contains an eigenvector of slope $0$. The same argument
applies if $m=n-1$ and $s < v_p(\pi)$.

Suppose instead that $m=n-1$ and $s = v_p(\pi)$.
Write $\bv = \bv_1 + \pi^{-1/(n-1)} \bv_2 + \cdots
+ \pi^{-(n-2)/(n-1)}
\bv_{n-1}$ with each $\bv_i$ defined over $\Galgancon$
(with no extension of $\calO[\fp]$); then $\bv_1, \dots, \bv_{n-1}$
are linearly independent in $N$, and we have
$F\bv_i = \bv_{i+1}$ for $i=1, \dots, n-1$ and $F\bv_{n-1} = 
\pi \bv_1$. In particular, the $\bv_i$ must be a basis of $N$
or else $\wedge^{n-1} N$ would have slope less than $s$.
The slope of $M/N$ is $-v_p(\pi)$, so by Lemma~\ref{lem:rank1},
we can choose $\bw \in M$ such that $F\bw \equiv \pi^{-1} \bw
\pmod{N}$.
Proposition~\ref{prop:direct} then implies that
$M$ contains an eigenvector of slope $0$, as desired.
\end{proof}

\begin{proof}[Proof of Theorem~\ref{thm:decomp}]
We proceed by induction on the rank of $M$. If $\rank M = 1$,
then $M$ is standard by Lemma~\ref{lem:rank1}. Suppose 
$\rank M = n> 1$, and that the proposition has been established
for all $\sigma$-modules of rank less than $n$.
For any rational number $c$,
define the \emph{$\calO$-index} of $c$ as the smallest integer $m$ such that
$mc$ lies in the value group of $\calO[\fp]$. 
The set of rational numbers of $\calO$-index less than or equal to $n$
which occur as slopes of eigenvectors of $M$ is discrete (obvious),
nonempty (by Proposition~\ref{prop:exist}), and bounded below
(by Corollary~\ref{cor:bdbelow}), so has a smallest element $r$.

Let $d$ be the slope of $\wedge^n M$.
By Lemma~\ref{lem:goodeigen}, we have $r \leq d/n$. Let $s$ be the
$\calO$-index of $r$, and let $\lambda$ be an element of a degree $s$
extension $\calO'[\fp]$ of $\calO[\fp]$ such that $v_p(\lambda) = r$ and
$\lambda^s \in \calO[\fp]$.
Choose an eigenvector $\bv$ over $\calO'[\fp]$ with $F\bv = \lambda \bv$,
and write $\bv = \sum_{i=0}^{s-1}
\lambda^{-i} \bw_i$ for $\bw_i \in M$, 
so that $F\bw_i = \bw_{i+1}$ for $i=0,\dots,s-2$
and $F\bw_{s-1} = \lambda^s \bw_0$.
Put $N = \Span(\bw_0, \dots, \bw_{s-1})$ and
$m=\rank N$; then $s \geq m$,
and the slope of $\wedge^m N$ is at most $mr$,
since $N$ is the saturated span of eigenvectors of slope $r$.

If $m=n$, then also $s=n$ and $\bw_0 \wedge \cdots \wedge \bw_{n-1}$
is an eigenvector of $\wedge^n M$ of slope $rn$. Thus $rn \geq d$;
since $r \leq d/n$ as shown earlier, we conclude $r=d/n$, $\bw_0,
\dots, \bw_{n-1}$ form a basis of $M$, and $M$ is standard, completing
the proof in this case. Thus we assume $m<n$ hereafter.

Given that $m<n$,
we may apply the induction hypothesis to $N$, deducing in particular
that its smallest slope is at most $r$ and has $\calO$-index not greater
than $m$. This yields a contradiction unless that slope is $r$, which is only
possible if the slope of $\wedge^m N$ is $mr$. In turn, $mr$ belongs to
the value group of $\calO[\fp]$ only if $m=s$. Thus $m=s$, and since
$\bw_0 \wedge \cdots \wedge \bw_{s-1}$ is an eigenvector of $N$ of slope $rs$,
$\bw_0, \dots, \bw_{s-1}$ form a basis of $N$, and $N$ is standard.

Apply the induction hypothesis to $M/N$ to express it as a sum
$P_1 \oplus \cdots \oplus P_l$ of standard $\sigma$-submodules;
note that the $\calO$-index of the slope of $P_i$ divides the rank
of $P_i$, and so is at most $n$.
If $l=1$, then the slope of $P_1$ cannot be less than $r$ (else the slope of
$\wedge^n M$ would be less than $d$), so by Proposition~\ref{prop:split},
$M$ can be split as
a direct sum of $N$ with a standard $\sigma$-module. 
If $l>1$, let $M_i$ be the preimage of $P_i$ under the projection $M \to M/N$;
again the slope
of each $P_i$ cannot be less than $r$, else the induction hypothesis would
imply that $M_i$ contains an eigenvector of slope less than $r$ and
$\calO$-index not exceeding $n$, contradiction. Thus by
Proposition~\ref{prop:split} again,
each $M_i$ can be split as a direct sum $N \oplus N_i$ of $\sigma$-submodules,
and we may decompose $M$ as $N \oplus N_1 \oplus \cdots \oplus N_l$.
This completes the induction in all cases.
\end{proof}

By Corollary~\ref{cor:newuniq},
the multiset union of the slopes of the standard summands of
a $\sigma$-module $M$ over $\Galgancon$ (each summand contributing
its slope as many times as its rank)
does not depend on the decomposition. Thus we define the
\emph{special Newton polygon} of $M$ as the polygon with
vertices $(i,y_i)$ ($i=0, \dots, n$), where $y_0=0$ and $y_i-y_{i-1}$
is the $i$-th smallest slope of $M$ (counting multiplicity).
We extend this definition to $\sigma$-modules over $\Gancon$ by
base extending to $\Galgancon$.

\section{The generic Newton polygon}
\label{sec:generic}

In this chapter, we recall the construction of the generic Newton polygon
associated to a $\sigma$-module over $\Gamma$.
The construction uses a classification result, the Dieudonn\'e-Manin
classification, for $\sigma$-modules over a complete discrete valuation
ring with algebraically closed residue field. This classification does not
descend very well, so we describe some weaker versions of the classification
that can be accomplished under less restrictive conditions.
These weaker versions either appear in or are inspired directly by
\cite{bib:dej1}.

\subsection{Properties of eigenvectors}

Throughout this section, let $R$ be a discrete valuation ring 
which is unramified over $\calO$.
Again, we call an element $\bv$ of a $\sigma$-module $M$ over $R$ or $R[\fp]$
an \emph{eigenvector}
if there exists $\lambda \in \calO$ or $\calO[\fp]$, respectively,
 such that $F\bv = \lambda \bv$,
and refer to $v_p(\lambda)$ as the \emph{slope} of $\bv$. 
We call an eigenvector \emph{primitive} if it
forms part of a basis of $M$, but this definition is not
very useful: every eigenvector is a $\calO$-multiple of a primitive
eigenvector of the same slope. In fact, in contrast with the situation
over $\Galgancon$, the slopes of eigenvectors over $R$ are ``rigid''.
\begin{prop} \label{prop:genslope1}
Let $M$ be a $\sigma$-module over $R[\fp]$, with $k$
algebraically closed. Suppose
$M$ admits a basis $\bv_1, \dots, \bv_n$ of eigenvectors. Then
any eigenvector $\bw$ is an $\calO[\fp]$-linear combination of
those $\bv_i$ of the same slope. In particular, any eigenvector has
the same slope as one of the $\bv_i$.
\end{prop}
\begin{proof}
Suppose $F\bv_i = \lambda_i \bv_i$ for some $\lambda_i \in \calO[\fp]$,
and write $\bw = \sum_i c_i \bv_i$ with $c_i \in R[\fp]$.
If $F\bw = \mu \bw$ for $\mu \in \calO[\fp]$, 
then equating the coefficients of $\bv_i$ yields
$\lambda_i c_i^\sigma = \mu c_i$. If $v_p(\lambda_i) \neq v_p(\mu)$, this forces
$c_i = 0$; if $v_p(\lambda_i) = v_p(\mu)$, it forces $c_i \in \calO[\fp]$.
This proves the claim.
\end{proof}

By imitating the proof of Proposition~\ref{prop:basis} using 
Proposition~\ref{prop:genslope1} in lieu of Proposition~\ref{prop:minslope},
we obtain the
following analogue of Corollary~\ref{cor:newuniq}.
\begin{prop} \label{prop:genslopes}
Let $M$ be a $\sigma$-module over $R[\fp]$. Suppose $\bv_1, \dots, \bv_n$
and $\bw_1, \dots, \bw_n$ are bases of eigenvectors with $F\bv_i =
\lambda_i \bv_i$ and $F\bw_i = \mu_i \bw_i$, for some $\lambda_i, \mu_i
\in \calO[\fp]$. Then the sequences $v_p(\lambda_1), \dots, v_p(\lambda_n)$
and $v_p(\mu_1), \dots, v_p(\mu_n)$ are permutations of each other.
\end{prop}

In case $M$ has a full set of eigenvectors of one slope, we have the following
decomposition result.
\begin{prop} \label{prop:stand}
Suppose $k$ is algebraically closed, and
let $M$ be a $\sigma$-module over $R$
spanned by eigenvectors of a single slope
over $R \otimes_{\calO} \calO'$, for some finite extension $\calO'$
of $\calO$.
Then $M$ is isogenous to the direct sum of standard $\sigma$-modules
of that slope.
\end{prop}
\begin{proof}
Let $s$ be the common slope, and let $m$ be the smallest positive
integer such that $ms$ is a multiple of $v_p(\pi)$. Since $k$
is algebraically closed, there exists $\lambda \in \calO'$
such that $\lambda^m \in \calO$. Let $\calO''$ be the integral
closure of $\calO$ in $\calO[\fp](\lambda)$.

Note that $M$ is spanned over $R \otimes_{\calO} \calO'$
by eigenvectors $\bv$ with $F\bv = \lambda \bv$ because $k$
is algebraically closed: if $F\bw = \mu \bw$ for some
$\mu$ with $v_p(\mu) = v_p(\lambda)$, we can find $c \in \calO'$
nonzero such that
$c^\sigma = (\lambda/\mu)c$ and obtain a new eigenvector
$\bv = c\bw$ with $F\bv = \lambda \bv$.
We next verify that $M$ is also spanned over $R \otimes_{\calO} \calO''$
by eigenvectors $\bv$ with $F\bv = \lambda \bv$.
Let $\mu_1, \dots, \mu_n$ be a basis of $\calO'$ over $\calO''$
consisting of elements fixed by $\sigma$ (possible because $k$
is algebraically closed).
If $\bv$ is an eigenvector over $R \otimes_{\calO} \calO'$
with $F\bv = \lambda \bv$, we can write $\bv = \sum_i \mu_i \bw_i$
for some $\bw_i$ over $R \otimes_{\calO} \calO''$, and
we must have $F\bw_i = \lambda \bw_i$ for each $i$. Thus $\bv$
is in the span of the $\bw_i$, so the span of eigenvectors of eigenvalue
$\lambda$ over $R \otimes_{\calO} \calO''$ has full rank over 
$R \otimes_{\calO} \calO'$, and thus has full rank over
$R \otimes_{\calO} \calO''$.

Finally, we establish that $M$ is isogenous to a direct sum
of standard $\sigma$-modules.
Let $\bv$ be an eigenvector of eigenvalue $\lambda$
over $R \otimes_{\calO} \calO''$;
we can write $\bv = \sum_{i=0}^{m-1} \bw_i \lambda^{-i}$ for some
$\bw_i \in M$. Then $F\bw_i = \bw_{i+1}$ for $i=0,\dots, m-2$
and $F\bw_{m-1} = \lambda^m \bw_0$, so the span
of $\bw_0, \dots, \bw_{m-1}$ is standard. (Notice that $\bw_0,
\dots, \bw_{m-1}$ must be linearly independent, otherwise any
eigenvalue of their span would lie in an extension of $\calO$
of degree strictly less than $m$, contrary to 
Proposition~\ref{prop:genslopes}.) Let $M_1$ be the standard
submodule just produced. Next, choose an eigenvector
of eigenvalue $\lambda$
linearly independent from $M_1$, and produce another standard
submodule $M_2$. Then choose an eigenvector linearly independent
from $M_1 \oplus M_2$, and so on until $M$ is exhausted.
\end{proof}

\subsection{The Dieudonn\'e-Manin classification}
\label{subsec:dieumanin}

Again, let $R$ be a discrete valuation ring unramified over $\calO$.
\begin{lemma} \label{lem:cycvec}
Suppose that $R$ is complete with algebraically closed residue field.
Given elements $a_0, \dots, a_{n-1}$ of $R$ with $a_0$ nonzero, let $M$ be
the $\sigma$-module with basis $\bv_1, \dots, \bv_n$ such that
\begin{align*}
  F\bv_i &= \bv_{i+1} \qquad (i=1, \dots, n-1) \\
F \bv_n &= a_0 \bv_1 + \cdots + a_{n-1} \bv_n.
\end{align*}
Suppose $s$ belongs to the value group of $R$.
Then the maximum number of linearly independent eigenvectors of slope $s$
in $M$
is less than or equal to the multiplicity $m$ of $s$ as a slope of
the Newton polygon
of the polynomial $x^n + a_{n-1} x^{n-1} + \cdots + a_0$ over $R$.
Moreover, if $m>0$, then $M$ admits an eigenvector of slope $s$.
\end{lemma}
\begin{proof}
Let $l = \min_j \{-js + v_p(a_{n-j}) \}$ (setting $a_n = 1$ for consistency);
then
there exists an index $i$ such that $l = -js + v_p(a_{n-j})$ for
$j=i$, $j=i+m$, and possibly for some values of $j \in \{i+1, \dots, i+m-1\}$,
but not for any other values.

Let $\lambda$ be an element of valuation $s$ fixed by $\sigma$.
Suppose $\bw = \sum_j c_j \bv_j$ satisfies $F\bw = \lambda \bw$.
Then $\lambda c_1 = a_0 c_n^\sigma$ and
$\lambda c_j = a_{j-1} c_n^\sigma + c_{j-1}^\sigma$ for $j=2, \dots, n$.
Solving for $c_n$ yields 
\begin{align*}
c_n &= \lambda^{-1} a_{n-1} c_n^\sigma + \lambda^{-1} c_{n-1}^\sigma \\
&= \lambda^{-1} a_{n-1} c_n^\sigma + \lambda^{-2} 
a_{n-2}^\sigma c_n^{\sigma^2} + \lambda^{-2} c_{n-2}^{\sigma^2}\\
&\qquad \vdots \\
&= \lambda^{-1} a_{n-1} c_n^\sigma +
 \lambda^{-2} a_{n-2}^\sigma
c_n^{\sigma^2} + \cdots +
\lambda^{-n+1} a_1^{\sigma^{n-2}} c_{n}^{\sigma^{n-1}} + 
\lambda^{-n+1} c_1^{\sigma^{n-1}} \\
&= \lambda^{-1} a_{n-1} c_n^\sigma +
 \lambda^{-2} a_{n-2}^\sigma c_n^{\sigma^2} + \cdots 
+ \lambda^{-n+1} a_1^{\sigma^{n-2}} c_{n}^{\sigma^{n-1}} 
+ \lambda^{-n} a_0^{\sigma^{n-1}} c_n^{\sigma^{n}}.
\end{align*}
In other words, $f(c_n) = 0$, where
\[
f(x) = 
-x + \frac{a_{n-1}}{\lambda} x^\sigma + \frac{a_{n-2}^\sigma}{\lambda^2}
x^{\sigma^2} + \cdots
+ \frac{a_0^{\sigma^{n-1}}}{\lambda^n} x^{\sigma^n}.
\]
The coefficients of $f$ of minimal valuation
are on $x^{\sigma^i}$, $x^{\sigma^{i+m}}$, and possibly some in between.

Now suppose $\bw_1, \dots, \bw_{m+1}$ are linearly independent eigenvectors
of $M$ with $F\bw_h = \lambda \bw_h$ for $h=1, \dots, m+1$.
Write $\bw_h = \sum_j c_{hj} \bv_j$. Then $c_{1n}, \dots, c_{(m+1)n}$
are linearly independent over $\calO_0$: if there were a relation
$\sum_h d_h c_{hn} = 0$ with $d_h \in \calO_0$ not all zero, we would
have 
\[
\lambda \sum_h d_h c_{hj} = \left( \sum_h d_h c_{h(j-1)} \right)^\sigma + 
a_{j-1} \left( \sum_h d_h c_{hn} \right)^\sigma \qquad (j=2,\dots, n)
\]
and successively deduce $\sum_h d_h c_{hj} = 0$ for $j=n-1, \dots, 1$. That
would mean $\sum_h d_h \bw_h = 0$, but the $\bw_h$ are linearly independent.

By replacing the $\bw_h$ with suitable $\calO_0$-linear combinations,
we can ensure that the $c_{hn}$ are in $R$
and their reductions modulo $\pi$ are linearly independent over $\FF_q$.
Now on one hand, the reduction of $(\lambda^{i}/a_{n-i}^{\sigma^{i-1}}) f(x)$
modulo $\pi$ is a polynomial in $x$ of 
the form $b_{i+m} x^{q^{i+m}} + \cdots + b_i x^{q^i}$,
which has only $q^m$ distinct roots in $R/\pi R$. On the other hand,
the $\FF_q$-linear combinations of the reductions of the $c_{hn}$ yield
$q^{m+1}$ distinct roots in $R/\pi R$, contradiction.

We conclude that the multiplicity of $s$ as a slope
of $M$ is at most $m$; 
this establishes the first assertion. To establish the second,
note that if $m>0$, then there exists $c_n \neq 0$ such that
$f(c_n) = 0$ by
Proposition~\ref{prop:hensel};
letting $c_n$ be this root, one can then solve for $c_{n-1}, \dots, c_1$
and produce an eigenvector $\bv$ with $F\bv = \lambda \bv$.
\end{proof}

Using this lemma, we can establish the Dieudonn\'e-Manin classification
theorem (for which see also Katz \cite{bib:katz}). We first state it
not quite in the standard form. Note: a ``basis up to isogeny'' means
a maximal linearly independent set.
\begin{prop} \label{prop:geneigs}
Suppose $R$ is complete with algebraically closed residue field. Then
every $\sigma$-module $M$ over $R$ has a basis up to isogeny 
of eigenvectors of
nonnegative slopes over
$R \otimes_{\calO} \calO'$ for some finite extension $\calO'$ of $\calO$
(depending on $M$).
\end{prop}
\begin{proof}
We proceed by induction on $n = \rank M$. Let $\bv$ be any nonzero element of
$M$, and let $m$ be the smallest integer such that $\bv, F\bv, \dots, F^m \bv$
are linearly dependent. Then $N = \Span(\bv, F\bv, \dots, F^{m-1} \bv)$
is a $\sigma$-submodule of $M$, and Lemma~\ref{lem:cycvec} implies that it
has a primitive
eigenvector $\bv_1$ of nonnegative slope
over $R \otimes_{\calO} \calO'$ for some $\calO'$
(since the corresponding polynomial has a root of nonnegative valuation
there). By the induction hypothesis, we can choose
$\bw_2, \dots, \bw_n$ over $R \otimes_{\calO} \calO''$ for some $\calO''$,
whose images in $M/\Span(\bv_1)$ form a basis up to isogeny
of eigenvectors
of nonnegative slopes.
We then have $F\bv_1 = \lambda_1 \bv_1$, where we may take $\lambda_1$
fixed by $\sigma$, and $F\bw_i = \lambda_i \bw_i +
c_i \bv_1$ for some $\lambda_i \in \calO$ and $c_i \in R$.
Apply
Proposition~\ref{prop:hensel} to find
$a_i \in R$ such that $\lambda_1 
c_i + \lambda_1 a_i^\sigma - \lambda_i a_i = 0$,
and set $\bv_i = \lambda_1 \bw_i + a_i \bv_1$ for $i=2, \dots, n$; then
$F\bv_i = \lambda_i \bv_i$, so $\bv_1, \dots, \bv_n$ form a basis up
to isogeny
of eigenvectors of nonnegative slope over $R \otimes_{\calO} \calO''$,
as desired.
\end{proof}
From this statement we deduce the
 Dieudonn\'e-Manin classification theorem in its more standard form.
\begin{theorem}[Dieudonn\'e-Manin] \label{thm:dm}
Suppose $R$ is complete with algebraically closed residue field. Then
every $\sigma$-module over $R$ is canonically isogenous to the direct
sum of $\sigma$-modules, each with a single slope, with all of these slopes
distinct. Moreover, every
$\sigma$-module of a single slope is isogenous to a direct sum of
standard $\sigma$-modules of that slope.
\end{theorem}
\begin{proof}
Let $M$ be a $\sigma$-module over $R$.
For each slope $s$ that occurs in a basis up to isogeny
of eigenvectors produced by
Proposition~\ref{prop:geneigs} over $R \otimes_{\calO} \calO'$,
let $M_s$ be the span of all eigenvectors of
$M$ over $R \otimes_{\calO} \calO'$ of slope $s$. Then $M_s$ is invariant under
$\Gal(\calO'/\calO)$, so by Galois descent, $M_s$ descends to a 
$\sigma$-submodule of $M$. Moreover, $M_s$ is isogenous to a direct
sum of standard $\sigma$-modules of slope $s$ by Proposition~\ref{prop:stand}.
This proves the desired result.
\end{proof}

Given a $\sigma$-module $M$ over a discrete valuation ring $R$ unramified over
$\calO$, we can embed $R$ into a complete discrete valuation ring
over which $M$ has a basis up to isogeny
of eigenvectors by Proposition~\ref{prop:geneigs}.
(First complete $R$, then take its maximal unramified extension, then
complete again, then tensor with a suitable $\calO'$ over $\calO$.)
By Proposition~\ref{prop:genslopes}, the slopes and multiplicities
do not depend on the choice of the basis.
Define the \emph{generic slopes} of $M$
as the slopes of the eigenvectors in the basis,
and the \emph{generic Newton polygon}
of $M$ as the polygon with vertices $(i, y_i)$ for $i = 0, \dots, \rank M$,
where $y_0 = 0$ and $y_i - y_{i-1}$ is the $i$-th smallest
generic slope of $M$ (counting multiplicity).
If $M$ has all slopes equal to 0, we say $M$ is \emph{unit-root}.

With the definition of the generic Newton polygon in hand, we can refine
the conclusion of Lemma~\ref{lem:cycvec} as follows.
\begin{prop} \label{prop:cycvec}
Given elements $a_0, \dots, a_{n-1}$ of $R$ with $a_0$ nonzero, let $M$ be
the $\sigma$-module with basis $\bv_1, \dots, \bv_n$ such that
\begin{align*}
  F\bv_i &= \bv_{i+1} \qquad (i=1, \dots, n-1) \\
F \bv_n &= a_0 \bv_1 + \cdots + a_{n-1} \bv_n.
\end{align*}
Then the generic Newton polygon of $M$ coincides with the
the Newton polygon
of the polynomial $x^n + a_{n-1} x^{n-1} + \cdots + a_0$ over $R$.
\end{prop}
\begin{proof}
The two Newton polygons have the same length $n$, and every number
occurs at least as often as a slope of the polynomial as it occurs as
a slope of $M$ by Lemma~\ref{lem:cycvec}. Thus all multiplicities must coincide.
\end{proof}
For our purposes, the principal consequence of this fact is the following.
\begin{prop} \label{prop:nodenom}
Let $M$ be a $\sigma$-module over $R[\fp]$ with all slopes nonnegative.
Then $M$ is isomorphic to a $\sigma$-module defined over $R$.
\end{prop}
\begin{proof}
We proceed by induction on $n = \rank M$. Let $\bv \in M$ be nonzero, and let
$m$ be the smallest integer such that $\bv, F\bv, \dots, F^m\bv$
are linearly dependent. Then 
$F^m \bv = a_0 \bv + \cdots + a_{m-1} F^{m-1}\bv$ for some
$a_0, \dots, a_{m-1} \in R[\fp]$;
by Proposition~\ref{prop:cycvec}, the $a_i$ belong to $R$.
Let $N = \Span(\bv, F\bv, \dots, F^{m-1}\bv)$; by the induction hypothesis,
$M/N$ is isomorphic to a $\sigma$-module defined over $R$. So
we can choose $\bw_1, \dots, \bw_{n-m}$ that form a basis of $M$ together
with $\bv, F\bv, \dots, F^{m-1} \bv$, such that for $i=1, \dots, n-m$,
$F\bw_i$ equals an $R[\fp]$-linear combination of the $F^j\bv$ plus an
$R$-linear combination of the $\bw_j$. For $\lambda$ sufficiently
divisible by $\pi$, the basis $\lambda \bv, \lambda F\bv, 
\dots, \lambda F^{m-1}\bv, \bw_1,
\dots, \bw_{n-m}$ has the property that the image of each basis vector
under Frobenius is an $R$-linear combination of basis vectors. This
gives the desired isomorphism.
\end{proof}


We close the section with another method for reading off the generic
Newton polygon of a $\sigma$-module, 
inspired by an observation of Buzzard and Calegari \cite[Lemma~5]{bib:bc}.
(We suspect it may date back earlier, possibly to Manin.)
\begin{prop} \label{prop:genspec}
Let $M$ be a $\sigma$-module over a discrete valuation ring $R$.
Suppose $M$ has
a basis on which $F$ acts by the matrix $A$, where $AD^{-1}$ is congruent
to the identity matrix modulo $\pi$ for some diagonal matrix $D$ over $\calO$.
Then the slopes of the generic Newton
polygon of $M$ equal the valuations of the diagonal entries of $D$.
\end{prop}
\begin{proof}
Without loss of generality we may assume $R$ is complete with algebraically
closed residue field.
We produce a sequence of matrices $\{U_l\}_{l=1}^\infty$ such that $U_1 = I$,
$U_{l+1} \equiv U_l \pmod{\pi^{l}}$ and $U_l^{-1} AU_l^\sigma D^{-1} \equiv
I\pmod{\pi^l}$;
the $\pi$-adic limit $U$ of the $U_l$ will satisfy $AU^\sigma = UD$,
proving the proposition. The conditions for $l=1$ are satisfied by
the assumption that $AD^{-1} \equiv I \pmod{\pi}$.

Suppose $U_l$ has been defined. Put $V = U_l^{-1} AU_l^\sigma D^{-1} -
I$. Define a matrix $W$ whose entry $W_{ij}$, for each $i$ and $j$,
is a solution of the
equation $W_{ij} - D_{ii} W_{ij}^\sigma D_{jj}^{-1} = V_{ij}$
with $\min \{v_p(W_{ij}), v_p(D_{ii} W_{ij}^\sigma D_{jj}^{-1})\}
= v_p(V_{ij})$ (such a solution exists by Proposition~\ref{prop:hensel}).
 Then
$W$ and $DW^\sigma D^{-1}$ are both congruent to 0 modulo $\pi^l$.
Put $U_{l+1} = U_l (I+W)$; then
\begin{align*}
U_{l+1}^{-1} A U_{l+1}^\sigma D^{-1} &=
(I+W)^{-1}U_l^{-1} A U_l^\sigma (I + W)^\sigma D^{-1} \\
&= (I+W)^{-1} U_l^{-1} A U_l^\sigma D^{-1} (I + D W^\sigma D^{-1}) \\
&= (I+W)^{-1}(I + V)(I + DW^\sigma D^{-1}) \\
&\equiv I-W+V+DW^\sigma D^{-1} = I \pmod{\pi^{l+1}}.  
\end{align*}
Thus the conditions for $U_{l+1}$ are satisfied, and the proposition
follows.
\end{proof}

\subsection{Slope filtrations}

The Dieudonn\'e-Manin classification holds over $\GK$ only if $K$
is algebraically closed, and even then does not descend to $\GKcon$ in general.
In this section, we exhibit two partial versions of the classification
that hold with weaker conditions on the coefficient ring. One (the
descending filtration) is due to de~Jong \cite[Proposition~5.8]{bib:dej1};
for symmetry, we present independent proofs of both results.

The following filtration result applies for
any $K$ but does not descend to $\Gcon$.
\begin{prop}[Ascending generic filtration] \label{prop:ascfilt}
Let $K$ be a valued field. Then
any $\sigma$-module $M$ over $\Gamma = \Gamma^K$ admits a unique
filtration
$M_0 = 0 \subset M_1 \subset \cdots \subset M_m = M$
by $\sigma$-submodules such that
\begin{enumerate}
\item 
for $i=1, \dots, m$, $M_{i-1}$ is saturated in $M_i$ and
$M_i/M_{i-1}$ has all generic slopes equal to $s_i$, and
\item
$s_1 < \cdots < s_m$.
\end{enumerate}
Moreover, if $K$ is separably closed and $k$ is algebraically closed, 
each $M_i/M_{i-1}$ is isogenous
to a direct sum of standard $\sigma$-modules.
\end{prop}
Warning: this proof uses the object $\Gsep$ even though this has only
so far been defined for $k$ perfect. Thus we must
give an \emph{ad hoc} definition here. For any finite separable extension
$L$ over $K$, Lemma~\ref{lem:sep1} produces a finite extension of
$\GK$ with residue field $L$, and Lemma~\ref{lem:sep2} allows us to identify
that extension as a subring of $\Galg$. We define $\Gsep$ as the completed
union of these subrings; note that $\Gperf \cap \Gsep = \Gamma$.

\begin{proof}
By the Dieudonn\'e-Manin classification 
(Theorem~\ref{thm:dm}), $M$ is canonically isogenous
to a direct sum of $\sigma$-submodules, each of a different single slope.
By Corollary~\ref{cor:galois}, these submodules descend to $\Gperf$;
let $M_1$ be the submodule of minimum slope. It suffices to show
that $M_1$ is defined over $\Gamma$, as an induction on rank will
then yield the general result. Moreover, it is enough to establish this
when $M_1$ has rank 1: if $M_1$ has rank $d$, then the lowest slope
submodule of $\wedge^d M$ is the rank one submodule $\wedge^d M_1$,
and if $\wedge^d M_1$ is defined over $\Gamma$, then so is $M_1$.

So suppose that $M_1$ has rank 1; this implies that the lowest slope of $M$
belongs to the value group of $\calO$. By applying an isogeny and then
twisting, we may reduce to the case where the lowest slope is 0.
Let $\be_1, \dots, \be_n$ be a basis of $M$ and let $A$ be the matrix
such that $F\be_j = \sum_{jl} A_{jl} \be_l$.

Let $\bv$ be an eigenvector of $M$ over $\Galg$ with $F\bv = \bv$.
We will show that $\bv$ is 
congruent to an element of $M \otimes_{\Gamma} \Gsep$
modulo $\pi^m$ for each $m$, by induction on $m$.
The case $m=0$ is vacuous, so assume the result is known for some $m$,
that is, $\bv = \bw + \pi^m \bx$ with $\bw \in M \otimes_{\Gamma}
\Gsep$ and $\bx \in
M \otimes_{\Gamma} \Galg$. Then $0 = F\bv - \bv = (F\bw - \bw) + \pi^m
(F\bx - \bx)$; that is, $F\bx - \bx$ belongs to $M \otimes_{\Gamma}
\Gsep$. Write
$\bx = \sum_j c_j \be_j$ and $F\bx - \bx = \sum d_j \be_j$,
and let $s$ be the smallest nonnegative integer
such that the reduction of $c_j$ modulo $\pi$ lies in $K^{1/q^s}$ for all $j$.
Then $d_j = - c_j + \sum_l A_{jl} c_l^\sigma$; if $s > 0$, then writing
$c_j = - d_j + \sum_l A_{jl} c_l^\sigma$ shows that the reduction of
$c_j$ lies in $K^{1/q^{s-1}}$ for all $j$, contradiction. Thus $s=0$,
and $\bx$ is congruent modulo $\pi$ to an element of $M \otimes_{\Gamma}
 \Gamma^{\sep}$, completing the induction.

We conclude that $\bv \in M \otimes_{\Gamma} \Gsep$. Thus $M_1$ is defined
both over $\Gperf$ and over $\Gsep$, so it is in fact defined over
$\Gperf \cap \Gsep = \Gamma$, as desired. This proves the desired result,
except for the final assertion. In case $K$ is separably closed, one can
repeat the above argument over a suitable finite extension of $\calO$
to show that each $M_i/M_{i-1}$ is spanned by eigenvectors, then
apply Proposition~\ref{prop:stand}.
\end{proof}

The following filtration result applies over $\Gcon$, not just over
$\Gamma$, but requires that $K$ be perfect.
\begin{prop}[Descending generic filtration] \label{prop:descfilt}
Let $K$ be a perfect valued field over $k$. Then
any $\sigma$-module $M$ over $\Gcon = \GKcon$ admits a unique filtration
$M_0 = 0 \subset M_1 \subset \cdots \subset M_m = M$
by $\sigma$-submodules such that
\begin{enumerate}
\item 
for $i=1, \dots, m$, $M_{i-1}$ is saturated in $M_i$ and
$M_i/M_{i-1}$ has all generic slopes equal to $s_i$, and
\item
$s_1 > \cdots > s_m$.
\end{enumerate}
Moreover, if $K$ is algebraically closed, each $M_i/M_{i-1}$ is isogenous
to a direct sum of standard $\sigma$-modules.
\end{prop}
\begin{proof}
By the Dieudonn\'e-Manin classification (Theorem~\ref{thm:dm}), 
$M$ is canonically isogenous
to a direct sum of $\sigma$-submodules, each of a different single slope.
By Corollary~\ref{cor:galois}, these submodules descend to $\Gamma$;
let $M_1$ be the submodule of maximum slope. It suffices to show
that $M_1$ is defined over $\Gcon$, as an induction on rank will
then yield the general result. Moreover, it is enough to establish this
when $M_1$ has rank 1: if $M_1$ has rank $d$, then the lowest slope
submodule of $\wedge^d M$ is the rank one submodule $\wedge^d M_1$,
and if $\wedge^d M_1$ is defined over $\Gcon$, then so is $M_1$.

So suppose that $M_1$ has rank 1; this implies that the highest slope of $M$
belongs to the value group of $\calO$. Choose $\lambda \in \calO$ whose
valuation equals that slope.
Let 
$\bv_1, \dots, \bv_n$ be a basis of $M \otimes_{\Gcon} \Galg$,
in which $F\bv_1 = \lambda \bv_1$ and
the remaining $\bv_i$ span the submodules of $M$ of lower slopes.
Choose $\bw_i \in M \otimes_{\Gcon} \Galgcon$ 
sufficiently close $\pi$-adically to $\bv_i$ for $i=1,
\dots, n$ so that the matrix $B$ with $\lambda \bw_i = \sum_j B_{ij} F \bw_j$
has entries in $\Galg$ and 
\[
B_{ij} \equiv \begin{cases} 1 & i = j = 1 \\ 0 & \mbox{otherwise} 
\end{cases}
\pmod{\pi};
\]
this is possible because the congruence holds for $\bw_i = \bv_i$. Then
the $\bw_i$ form a basis of $M \otimes_{\Gcon} \Galgcon$.

Write $\bv_1 = \sum_i c_i \bw_i$, so that
$c_i^\sigma = \sum_j B_{ji} c_j$. 
Since $v_0(B) \geq 0$, we can find $r$ such that $w_r(B) \geq 0$.
We now show that $r v_h(c_i) + h \geq 0$ for all $i$ and $h$, by induction
on $h$. Suppose this holds with $h$ replaced by any smaller value.
Then the equality $c_i^\sigma = \sum_j B_{ji} c_j$ implies
\[
q v_h(c_i) \geq \min_{l,j} \{v_l(B_{ji}) + v_{h-l}(c_j)\}.
\]
Choose $j,l$ for which the minimum is achieved.
If $l=0$, then we must have $i=j=1$,
in which case $v_0(B_{11}) = 0$ and $q v_h(c_1) \geq v_h(c_1)$, whence
$v_h(c_1) \geq 0$ and $r v_h(c_1) + h \geq 0$ as well.
If the minimum occurs for some $l>0$, then
\begin{align*}
r v_h(c_i) + h &\geq r q^{-1} (v_l(B_{ji}) + v_{h-l}(c_j)) + h\\
 &\geq r q^{-1} (v_l(B_{ji}) + v_{h-l}(c_j)) + q^{-1} h\\
&\geq q^{-1} (r v_l(B_{ji}) + l + r v_{h-l}(c_j) + (h-l)) \\
&\geq q^{-1}(0 + 0) = 0
\end{align*}
by the induction hypothesis. Therefore $r v_h(c_i) + h \geq 0$ for all $h$,
so $c_i \in \Galgcon$ for each $i$.

We conclude that $\bv_1 \in M \otimes_{\Gamma} \Galgcon$. Thus $M_1$ is defined
both over $\Gamma$ and over $\Galgcon$, so it is in fact defined over
$\Gamma \cap \Galgcon = \Gcon$, as desired. This proves the desired result,
except for the final assertion. In case $K$ is algebraically closed, one can
repeat the above argument over a suitable finite extension of $\calO$
to show that each $M_i/M_{i-1}$ is spanned by eigenvectors, then
apply Proposition~\ref{prop:stand}.
\end{proof}
Although we will not use the following result explicitly, it is worth pointing out.
\begin{cor} \label{cor:sepcon}
Let $K$ be a valued field, for $k$ algebraically closed.
Then any $\sigma$-module $M$ over $\GKcon$,
all of whose generic slopes are equal, is isogenous over $\Gsepcon$
to a direct sum of standard $\sigma$-modules.
\end{cor}
\begin{proof}
In this case, the ascending and descending filtrations coincide, so both
are defined over $\Gamma^K \cap \Gperfcon = \GKcon$ and the eigenvectors
are defined over $\Gsepcon \otimes_{\calO} \calO'$ for some finite
extension $\calO'$ of $\calO$. Thus the claim follows from
Proposition~\ref{prop:stand}.
\end{proof}


\subsection{Comparison of the Newton polygons}

A $\sigma$-module over $\Gcon$ can be base-extended both to $\Gamma$
and to $\Gancon$; as a result, it admits both a generic and a special
Newton polygon. In this section, we compare these two polygons.
The main results are that the special polygon lies above the generic
polygon, and that when the two coincide, the $\sigma$-module admits a
partial decomposition over $\Gcon$ (reminiscent of 
the Newton-Hodge decomposition of \cite{bib:katz}).

Throughout this section, $K$ is an arbitrary valued field, which we
suppress from the notation.

\begin{prop} \label{prop:operations}
Let $M$ and $N$ be $\sigma$-modules over $\Gcon$.
Let $r_1, \dots, r_m$ and $s_1, \dots, s_n$ be the generic (resp.\ special)
slopes of $M$ and $N$.
\begin{enumerate}
\item The generic (resp.\ special) slopes of $M \oplus N$ are
$r_1, \dots, r_m, s_1, \dots, s_n$.
\item The generic (resp.\ special) slopes of $M \otimes N$ are
$r_i + s_j$ for $i=1, \dots, m$ and $j=1, \dots, n$.
\item The generic (resp.\ special) slopes of $\wedge^l M$ are
$r_{i_1} + \cdots + r_{i_l}$ for $1 \leq i_1 < \cdots < i_l \leq m$.
\item The generic (resp.\ special) slopes of $M^*$ are
$-r_1, \dots, -r_m$.
\end{enumerate}
\end{prop}
\begin{proof}
These results follow
immediately from the definition of the generic (resp.\ special)
Newton slopes as the valuations of the eigenvalues
of a basis of eigenvectors of $M$ over $\Galg$ (resp.\ $\Galgancon$).
\end{proof}

\begin{prop} \label{prop:above}
Let $M$ be a $\sigma$-module over $\Gcon$.
Then the special Newton polygon lies above the generic Newton
polygon, and both have the same endpoint.
\end{prop}
\begin{proof}
The Newton polygons coincide for $M$ of rank~1 because $M$ has an eigenvector
over $\Galgcon$ by Proposition~\ref{prop:rank1}. Thus the Newton polygons
of $\wedge^n M$ coincide for $n = \rank M$; that is, the Newton polygons
of $M$ have the same endpoint. 
By the descending slope filtration (Proposition~\ref{prop:descfilt}),
$M$ admits a basis $\bw_1, \dots, \bw_n$ over $\Galgcon$
such that modulo
$\Span(\bw_1, \dots, \bw_{i-1})$, $\bw_i$ is an eigenvector whose
slope is the $i$-th largest generic slope of $M$.
Let $\bv_1, \dots, \bv_n$ be a basis of eigenvectors of $M$ over
$\Galgancon$; then by Proposition~\ref{prop:basis},
the sequence of valuations of the eigenvalues of the $\bw_i$
majorizes that of the $\bv_i$. In other words, the sequence
of generic slopes majorizes the sequence of special slopes,
whence the comparison of Newton polygons.
\end{proof}

\begin{prop} \label{prop:union}
Let $0 \to M_1 \to M \to M_2 \to 0$
be an exact sequence of $\sigma$-modules over $\Gcon$. Suppose
the least generic
slope of $M_2$ is greater than the greatest generic slope of
$M_1$. Then the special Newton polygon of $M$ is equal to the union of
the special Newton polygons of $M_1$ and $M_2$.
\end{prop}
\begin{proof}
The least generic slope of $M_2$ is less than or equal to its least special
slope, and the greatest generic slope of $M_1$ is greater than or equal to
its greatest special slope, both by Proposition~\ref{prop:above}.
Thus we may apply Proposition~\ref{prop:split} over $\Galgancon$
(after extending $\calO$
suitably) to deduce the desired result.
\end{proof}

It is perhaps
not surprising that when the generic and special Newton polygons
coincide, one gets a slope filtration
that descends further than usual.
\begin{prop} \label{prop:asc2}
Let $M$ be a $\sigma$-module over $\Gcon$ whose generic and 
special Newton
polygons coincide. Then $M$ admits an ascending slope filtration
over $\Gcon$.
\end{prop}
\begin{proof}
We need to show that the ascending slope
filtration of Proposition~\ref{prop:ascfilt} is defined over
$\Gcon$; it is enough to verify this after enlarging $\calO$.
This lets us assume that $k$ is algebraically closed,
and that the value group of
$\calO$ contains all of the slopes of $M$.
By Theorem~\ref{thm:decomp}, we can find
a basis $\bv_1, \dots, \bv_n$ of eigenvectors
of $M$ over $\Galgancon$, with $F\bv_i = \lambda_i \bv_i$ for
$\lambda_i \in \calO_0[\fp]$ such that $v_p(\lambda_1) \geq \cdots \geq 
v_p(\lambda_n)$.
By Proposition~\ref{prop:descfilt} (the descending slope filtration),
we can find a basis up to isogeny $\bw_1, \dots, \bw_n$ of $M$ over $\Galgcon$ such that
$F\bw_i = \lambda_i \bw_i + \sum_{j<i} A_{ij} \bw_j$ for some
$A_{ij} \in \Galgcon$.

Write $\bv_n = \sum_i b_i \bw_i$ with $b_i \in \Galgancon$;
applying $F$ to both sides,
we have $\lambda_n b_i = \lambda_i b_i^\sigma + \sum_{j>i} b_j^\sigma A_{ji}$
for $i=1, \dots, n$. 
By Proposition~\ref{prop:ansigeq}(a) and (c),
we obtain $b_i \in \Galgcon[\fp]$ for $i=n, n-1, \dots, 1$, and so
$\bv_n$ is defined over $\Galgcon[\fp]$. 

By repeating the above reasoning,
we see that the image of $\bv_i$ in $M/\Span(\bv_{i+1}, \dots, \bv_n)$
is defined over $\Galgcon[\fp]$ for $i=n, \dots, 1$. Thus the ascending slope
filtration is defined over $\Galgcon[\fp]$. Since it is also defined over
$\Gamma$ by Proposition~\ref{prop:ascfilt}, it is in fact defined over
$\Gamma \cap \Galgcon[\fp] = \Gcon$, as desired.
\end{proof}

\section{From a slope filtration to quasi-unipotence}
\label{sec:conclude}

In this chapter we construct 
a canonical filtration of a $\sigma$-module over $\Gamma^{k((t))}_{\an,\con}$.
We do this by partially descending
the special slope filtration obtained over $\Galgancon$ in 
Chapter~\ref{sec:special}. 
More specifically, we show that by changing basis over a finite
extension of $\Gancon$, we can make Frobenius act by a matrix with entries
in a finite extension of $\Gcon$, whose generic Newton polygon coincides with
the special Newton polygon, allowing the use of Proposition~\ref{prop:asc2}.
This will yield the desired filtration (Theorem~\ref{thm:filt}), from which we deduce the $p$-adic local monodromy theorem (Theorem~\ref{thm:monodromy})
using the 
quasi-unipotence of unit-root $(\sigma, \nabla)$-modules over $\Gcon$; the 
latter is a theorem of Tsuzuki \cite{bib:tsu1} (for which see also Christol
\cite{bib:ch}).

\subsection{Approximation of matrices}
\label{subsec:approx}

We collect some results that allow us to approximate matrices from a large
ring with matrices from smaller rings. Note: we will need the notions of slopes
and Newton polygons from Section~\ref{subsec:anbezout}.

\begin{lemma} \label{lem:addpi}
Let $K$ be a finite extension of $k((t))$ and suppose
$\GK_r$ contains a unit lifting a uniformizer of $K$. Then
for any $x,y \in \GK_r[\fp]$,
$x$ is coprime to $y + \pi^j$ for all sufficiently large integers $j$.
\end{lemma}
\begin{proof}
Suppose on the contrary that $x$ and $y + \pi^{j}$ fail to be coprime
for $j = j_1, j_2, \dots$.
By Corollary~\ref{cor:bezout}, the ideal
$(x, y + \pi^{j_l})$ in $\GK_{r}[\fp]$ is principal; let $d_l$ be
a generator. Note that $(y + \pi^{j_i}, y+\pi^{j_l})$ contains the
unit $\pi^{j_i} - \pi^{j_l}$ for $i \neq l$, so is the unit ideal;
that means the $d_l$ are pairwise coprime, and $x$ is divisible by
$d_1 \cdots d_l$ for any $l$. But $x$ has only finite total multiplicity
while each $d_l$ has nonzero total multiplicity, contradiction.
Hence $x$ is coprime to $y + \pi^j$ for $j$ sufficiently large, as desired.
\end{proof}


By an \emph{elementary operation} on a matrix over a ring, 
we mean one of the following
operations:
\begin{enumerate}
\item[(a)]
adding a multiple of one row to another;
\item[(b)]
multiplying one row by a unit of the ring;
\item[(c)]
interchanging two rows.
\end{enumerate}
An \emph{elementary matrix} is one obtained from the identity matrix by
a single elementary operation; multiplying a matrix on the right by
an elementary matrix has the same effect as performing the corresponding
elementary operation.

\begin{lemma} \label{lem:approx}
Pick $s$ such that $0 < s < r$, and let
$U$ be a matrix over 
$\Galg_{\an,r}$ such that $w_l(\det(U) - 1) > 0$ for $s \leq l \leq r$.
Then there exists an invertible matrix $V$ over $\GK_{r}[\fp]$,
for some finite extension $K$ of $k((t))$,
such that $w_l(UV-I) > 0$ for $s \leq l \leq r$.
Moreover, if $U$ is defined over $\Gamma^{k((t))}_{\an,r}$
and $t$ lifts to a semiunit in $\Gamma^{k((t))}_r$,
then we may take $K=k((t))$.
\end{lemma}
Although we only will apply this when $U$ is invertible, we need to
formulate the more general statement in order to carry out the induction.
\begin{proof}
We induct on $n$, the case $n=1$ being vacuous.
Let $M_i$ denote the cofactor of $U_{ni}$ in $U$,
so that $\det(U) = \sum_i M_i U_{ni}$; note that
$M_i = (U^{-1})_{in} \det(U)$ in $\Frac(\Galg_{\an,r})$.
Let $d$ be a generator of the ideal $(M_1, \dots, M_n)$ in $\Galg_{\an,r}$. 
Then
$d$ divides
$\det(U)$; by the hypothesis that $w_l(\det(U)-1) >0$ for $s \leq l
\leq r$, the largest slope of $\det(U)$ 
is less than $s$, so the largest slope of $d$ is also less than $s$.
By Lemma~\ref{lem:sorite}, there exists a unit $u \in \Galg_{\an,r}$
such that $w_l(u d-1) > 0$ for $s \leq l \leq r$.

Let
$\alpha_1, \dots, \alpha_n$ be elements of $\Galg_{\an,r}$ such that
$\sum_i \alpha_i M_{i} = ud$. Choose $\beta_1, \dots, \beta_{n-1}$ and
$\beta_n' \in \GL_r[\fp]$, for some finite extension $L$ of $k((t))$, so
that for $s \leq l \leq r$,
\[
w_l(\beta_i - \alpha_i) > 0 \quad
(i=1, \dots, n-1),
\qquad
\mbox{and} \qquad
w_l(\beta'_n - \alpha_n) > 0.
\]
By Lemma~\ref{lem:addpi}, we can find $j$ for which $\beta_n = \beta'_n + 
\pi^j$ has the properties that $w_l(\beta_n - \alpha_n) > 0$ for 
$s \leq l \leq r$ and
$(\beta_1, \dots, \beta_n)$ is the unit ideal in $\GL_r[\fp]$. (Both hold
for $j$ sufficiently large.)

By Corollary~\ref{cor:bezout}, $\GL_r[\fp]$ is a B\'ezout ring.
Thus Lemma~\ref{lem:det1} can be applied to produce a matrix $A$ 
over $\Gamma^{L}_r[\fp]$ of determinant 1 such that $A_{ni} = \beta_i$
for $i=1, \dots, n$.
Put $U' = U A^{-1}$, 
and let $M'_n$ be the cofactor of $U'_{nn}$ in $U'$.
Then
\begin{align*}
M'_n &= ((U')^{-1})_{nn} \det(U') \\
&= (A U^{-1})_{nn} \det(U) \det(A^{-1}) \\
&= \sum_i A_{ni} (U^{-1})_{in} \det(U) \\
&= \sum_i \beta_i M_i,
\end{align*}
so that
\[
M'_n - 1 = ud - 1 + \sum_i (\beta_i - \alpha_i)M_i
\]
and hence $w_l(M'_n - 1) > 0$ for $s \leq l \leq r$.

Apply the induction
hypothesis to the upper left $(n-1) \times (n-1)$ submatrix of $U'$,
let $V'$ be the resulting matrix, and enlarge $L$ if needed so that 
$V'$ has entries in $\GL_{\an,r}$.
Extend $V'$ to an $n \times n$ matrix
by setting $V'_{nn} = 1$ and $V'_{ni} = V'_{in} = 0$ for $i=1,
\dots, n-1$.
Then for $s \leq l \leq r$, $w_l((U'V' - I)_{ij}) > 0$ for $1 \leq i,j
\leq n-1$.
Moreover, $w_l(\det(V') - 1) > 0$, so $w_l(\det(U'V') - 1) > 0$ as well.

We now
exhibit a sequence of elementary operations which can be performed
on $U'V'$ to obtain a new matrix $W$ over $\Galg_{\an,r}$
with $w_l(W-I) > 0$
for $s \leq l \leq r$; it may clarify matters to regard the procedure
as an ``approximate Gaussian elimination''.
First, define a sequence of matrices $\{X^{(h)}\}_{h=0}^\infty$ by
$X^{(0)} = U'V'$ and 
\[
X^{(h+1)}_{ij} = \begin{cases} X^{(h)}_{ij} & i < n \\
X^{(h)}_{nj} - \sum_{m=1}^{n-1} X^{(h)}_{nm} X^{(h)}_{mj} & i = n;
\end{cases}
\]
note that $X^{(h+1)}$ is obtained from $X^{(h)}$ by subtracting $X^{(h)}_{nm}$
times the $m$-th row from the $n$-th row for $m=1, \dots, n-1$.
At each step, $\min_{1 \leq j \leq n-1} \{ w_l(X^{(h)}_{nj})\}$ 
increases by at least $\min_{1 \leq i,j
\leq n-1} \{ w_l((U'V' - I)_{ij}) \}$; thus for $h$ sufficiently large, 
we have
\[
w_l(X^{(h)}_{nj}) > \max\left\{0, \max_{1 \leq i \leq n-1} \{ -w_l(X^{(h)}_{in}) \} \right\}
\qquad (s \leq l \leq r; \, j=1, \dots, n-1).
\]
Pick such an $h$ and set $X = X_h$.
Then $w_l((X - I)_{ij}) > 0$ for $1 \leq i \leq n$ and $1 \leq j \leq n-1$,
$w_l(X_{in} X_{nj}) > 0$ for $1 \leq i,j \leq n-1$,
and $w_l(\det(X)-1) > 0$. These together imply
$w_l(X_{nn} - 1) > 0$. 

Next, define a sequence of matrices $\{W^{(h)}\}_{h=0}^\infty$ by
$W^{(0)} = X$ and 
\[
W^{(h+1)}_{ij} = \begin{cases} 
W^{(h)}_{ij} - W^{(h)}_{in} W^{(h)}_{nj} & i < n \\
W^{(h)}_{ij} & i = n;
\end{cases}
\]
note that $W^{(h+1)}$ is obtained from $W^{(h)}$ by subtracting $W^{(h)}_{in}$
times the $n$-th row from the $i$-th row for $i=1, \dots, n-1$.
At each step, $w_l(X^{(h)}_{in})$ increases by at least 
$w_l(X^{(h)}_{nn} - 1)$; thus for $h$ sufficiently large, we have
\[
w_l(W^{(h)}_{in}) > 0 \qquad (s \leq l \leq r; \, i=1, \dots, n-1).
\]
Pick such an $h$ and set $W = W_h$; then $w_l(W-I) > 0$ for 
$s \leq l \leq r$.

To conclude, note that by construction, $(U'V')^{-1} W$ is a
product of elementary matrices over $\Galg_{\an,r}$ of type (a).
By suitably approximating each elementary matrix by one defined
over $\GK_r[\fp]$ for a suitable finite extension $K$ of $L$,
we get a matrix $X$ such that $w_l(U'V'X-I) > 0$ for $s \leq l
\leq r$.
We may thus take $V = A^{-1}V'X$.
\end{proof}

We will need a refinement of the above result.
\begin{lemma} \label{lem:approx2}
Pick $s$ such that $0 < s < r$, and let
$U$ be a matrix over 
$\Galg_{\an,r}$ such that $w_l(\det(U) - 1) > 0$ for $s \leq l \leq r$.
Then for any $c>0$, there exists a finite extension $K$ of $k((t))$
and an invertible matrix $V$ over $\GK_{r}[\fp]$
such that $w_l(UV-I) \geq c $ for $s \leq l \leq r$.
Moreover, if $U$ is defined over $\Gamma^{k((t))}_{\an,r}$
and $t$ lifts to a semiunit in $\Gamma^{k((t))}_r$,
then we may take $K=k((t))$.
\end{lemma}
\begin{proof}
Put
\[
s' = s (1 + c/v_p(\pi))^{-1}.
\]
Apply Lemma~\ref{lem:approx} to obtain a finite extension $L$
of $k((t))$ and an invertible matrix $V'$ over $\GL_r[\fp]$ (with
$L = k((t))$ in case $U$ is defined over $\Gamma^{k((t))}_{\an,r}$) such that
$w_l(UV'-I) > 0$ for $s' \leq l \leq r$.

Choose semiunit decompositions $\sum_h W_{ijh} \pi^h$ of 
$(UV')_{ij}-I$ for $1 \leq i,j \leq n$.
For $s \leq l \leq r$ and $m<0$ in the value group of
$\calO$, we deduce from $w_{s'}(UV'-I) >0 $ that
\begin{align*}
l v_m(UV' - I) + m &=
(l/s') (s' v_m(U V' - I) + m)
- m (l/s' - 1) \\
&> -m(l/s' - 1) \\
&> v_p(\pi)(s/s' - 1) \\
&= c.
\end{align*}
Define a matrix $X$ by $X_{ij} = \sum_{h\geq 0} W_{ijh} \pi^h$;
then $UV' - I - X = \sum_{h<0} W_{ijh} \pi^h$, so that for $s \leq l \leq r$,
\[
w_l(UV' - I - X) 
= \min_{m < 0} \{l v_m(UV' - I) + m \} \geq c.
\]
By construction, $v_m(X) = \infty$ for $m<0$ and $v_0(X) > 0$. Thus
$I+X$ is invertible over $\Galg_{\an,r}$.
Choose a matrix $W$ over $\GK_r$, for some finite extension $K$
of $L$ (with $K=k((t))$ if $U$ is defined over $\Gamma^{k((t))}_{\an,r}$), 
such that $w_l(W-(I+X)^{-1}) \geq c$ for $s \leq l \leq r$.
Then $W$ is invertible over $\GK_r$, and for $s \leq l \leq r$,
\begin{align*}
w_l(UV'W - I) &= w_l((UV'-I-X)W + (I+X)(W - (I+X)^{-1})) \\
&\geq \min\{ w_l(UV' - I - X) + w_l(W), w_l(I + X) + w_l(W - (I+X)^{-1}) \} \\
&\geq c.
\end{align*}
We may thus take $V = V' W$.
\end{proof}

\subsection{Some matrix factorizations}

Throughout this section, we take $K = k((t))$ and omit it from the notation;
note also the use of the na\"\i ve partial valuations.
Let $\Gamma_{u}$ and $\Gamma_{\an,u}$
denote the subrings of 
$\Gcon$ and $\Gancon$, respectively, consisting of elements $x$ of the
form $\sum_{i=0}^\infty x_i u^i$.

\begin{lemma} \label{lem:fact}
For $r>0$ 
and $c>0$, let $A$ be a matrix over $\Gamma_{\an,r}$
such that $w_r^{\naive}(A-I) \geq c$.
Then there exists a unique pair of 
matrices
$U = I + \sum_{i=1}^\infty U_i u^i$ over $\Gamma_{\an,r}$ 
and 
$V = \sum_{i=0}^\infty V_i u^{-i}$ over $\Gamma_r$
such that $w_r^{\naive}(U-I) > 0$, $w_r^{\naive}(V - I) > 0$,
and $A = UV$. Moreover, these matrices satisfy
$w_r^{\naive}(U-I) \geq c$ and $w_r^{\naive}(V-I) \geq c$.
\end{lemma}
\begin{proof}
Define a sequence of matrices $\{B^{(j)}\}_{j=0}^\infty$ as follows.
Begin by setting $B^{(0)} = I$. Given $B^{(j)}$ for some $j$,
put $A (B^{(j)})^{-1} = \sum_{i=-\infty}^\infty X_i^{(j)} u^i$,
$C^{(j)} = \sum_{i \leq 0} X_i^{(j)} u^i$, $D^{(j)}
= \sum_{i>0} X_i^{(j)} u^i$, and put
$B^{(j+1)} = C^{(j)} B^{(j)}$.

Since $w_r^{\naive}(A-I) \geq c$, we have
$w_r^{\naive}(C^{(0)} - I) \geq c$ and $w_r^{\naive}(D^{(0)}) \geq c$ as well.
Thus $w_r^{\naive}(A (B^{(1)})^{-1} - I) \geq c$, and by induction
one has $w_r^{\naive}(C^{(j)}-I) \geq c$ and $w_r^{\naive}(D^{(j)}) \geq c$
for all $j$. But we can do better,
by showing
by induction that $w_r^{\naive}(C^{(j)} - I) \geq (j+1)c$ and
$w_r^{\naive}(D^{(j+1)} - D^{(j)}) \geq (j+2)c$ for $j \geq 0$.
Given $w_r^{\naive}(C^{(j)} - I) \geq (j+1)c$, we have
\begin{align*}
A(B^{(j+1)})^{-1} - I &= A (B^{(j)})^{-1} (C^{(j)})^{-1} - I \\
&= (C^{(j)} + D^{(j)})(C^{(j)})^{-1} - I \\
&= D^{(j)} (C^{(j)})^{-1} \\
&= D^{(j)} + D^{(j)}((C^{(j)})^{-1} - I).
\end{align*}
Since $D^{(j)}$ has only positive powers of $u$, $C^{(j+1)}$ is equal to
the sum of the terms of $I + D^{(j)}
((C^{(j)})^{-1}-I)$ involving nonpositive powers of $u$.
In particular, 
\[
w_r^{\naive}(C^{(j+1)}-I) \geq w_r^{\naive}(D^{(j)} ((C^{(j)})^{-1}-I)) \geq c + (j+1)c = (j+2)c;
\]
likewise, $D^{(j+1)}-D^{(j)}$ consists of terms from $D^{(j)}
((C^{(j)})^{-1}-I)$, so $w_r^{\naive}(D^{(j+1)}-D^{(j)}) \geq (j+2)c$.
This completes the induction.

Since $C^{(j)}$ converges to $I$,
we see that $B^{(j)}$ converges to
a limit $V$ such that $w_r^{\naive}(V-I) \geq c$.
Under $w_r^{\naive}$, $I + D^{(j)}$ also converges to a limit $U$
such that $w_r^{\naive}(U - I) \geq c$,
and $A(B^{(j)})^{-1} - I - D^{(j)}$ converges to 0.
Therefore $AV^{-1} = U$ has entries in $\Gamma_{\an,r,\naive}$,
and $U$ and $V$ satisfy the desired conditions.

This establishes the existence of the desired factorization. To establish
uniqueness, suppose we have a second decomposition
$A = U'V'$ with $U'-I$ only involving
positive powers of $u$, $V'$ only involving negative powers of $u$,
$w_r^{\naive}(U' - I) > 0$, and $w_r^{\naive}(V' - I) > 0$.
Within the completion of $\Gamma_r[\fp]$ with
respect to $|\cdot|_r$, the matrices $U,V,U',V'$ are invertible
and $(U')^{-1}U = V'V^{-1}$. On the other hand, $(U')^{-1}U - I$
involves only positive powers of $u$, while $V'V^{-1}-I$ involves
no positive powers of $u$. This is only possible if $(U')^{-1}U - I
= V'V^{-1} - I = 0$, which yields $U = U'$ and $V=V'$.
\end{proof}

The following proposition may be of interest outside of its use to prove
the results of this paper.
For example, Berger's proof \cite[Corollaire~0.3]{bib:berger} that any crystalline
representation is of finite height uses a lemma from \cite{bib:methesis}
equivalent to this.
\begin{prop} \label{prop:matfact}
Let $A = \sum_{i=-\infty}^\infty A_i u^i$
be an invertible matrix over $\Gancon$. Then
there exist invertible matrices $U$ over $\Gamma_{\an,u}$ and
$V$ over $\Gcon[\fp]$ such that $A = UV$. Moreover, 
if $w_r^{\naive}(A-I) >0$
for some $r>0$, there is a unique choice of $U$ and $V$
such that 
$U-I$ involves only positive powers of $u$,
$V$ involves no positive powers of $u$,
$w_r^{\naive}(U - I) > 0$ and $w_r^{\naive}(V - I) > 0$; for
these $U$ and $V$,
$\min\{w_r^{\naive}(U - I), w_r^{\naive}(V - I)\} \geq w_r^{\naive}(A-I)$.
\end{prop}
\begin{proof}
 By Lemma~\ref{lem:approx}, there exists an invertible matrix $W$
over $\Gcon[\fp]$ such that $w_r^{\naive}(AW-I) > 0$. Apply
Lemma~\ref{lem:fact} to write $AW = U_1V_1$ for matrices $U_1$
over $\Gamma_{\an,u}$ and $V_1$ over $\Gcon$, and to write
$(AW)^{-T}= U_2 V_2$ for matrices $U_2$ over $\Gamma_{\an,u}$
and $V_2$ over $\Gcon$.
Now $I = (AW)^T (AW)^{-T} = V_1^T U_1^T U_2V_2$, and so
$V^{-T}_1 V_2^{-1} = U_1^T U_2$
has entries in $\Gcon \cap \Gamma_{\an,u} = \Gamma_u$. Moreover,
$U_1^T U_2 - I$ involves only positive powers of $u$, so $U_1^T U_2$
is invertible over $\Gamma_u$ and $U_1$ is invertible over $\Gamma_{\an,u}$.
Our desired factorization
is now $A = UV$ with $U = U_1$ and $V = V_1 W^{-1}$.
If $w_r^{\naive}(A - I) > 0$, we may take
$W=I$ above and deduce the uniqueness from Lemma~\ref{lem:fact}.
\end{proof}

So far we have exhibited factorizations that separate positive and negative
powers of $u$. We use these to give a factorization that separates
a matrix over $\Gancon$ into a matrix over $\Gcon$ times a matrix with
only positive powers of $u$, in such a way that the closer the original matrix
is to being defined over $\Gcon$, the smaller the positive matrix will be.
\begin{prop} \label{prop:matfact2}
Let $A$ be an invertible matrix over $\Gamma_{\an,r}$
such that $w_r^{\naive}(A - I) > 0$. Then there exists a canonical pair of
invertible matrices $U$ over $\Gamma_{\an,u}$ and $V$ over $\Gcon$ such that
$A = UV$, $U-I$ has only positive powers of $u$, $V - I \equiv 0 \pmod{\pi}$,
$w_r^{\naive}(V - I) \geq w_r^{\naive}(A-I)$ and
\[
w_r^{\naive}(U-I) \geq \min_{m \leq 0} \{r v_m^{\naive}(A - I) + m\}.
\]
\end{prop}
Here ``canonical'' does not mean ``unique''. It means that the construction
of $U$ and $V$ depends only on $A$ and not on $r$.
\begin{proof}
Write $A-I = \sum_i A_i u^i$, and let $X$ be the sum of $A_i$ over
all $i$ for which $v_p(A_i) > 0$. Then
\begin{align*}
w_r^{\naive}(A(I+X)^{-1} - I) &\geq w_r^{\naive}(A - I - X) + 
w_r^{\naive}((I+X)^{-1}) \\
&= \min_{v_p(A_i) \leq 0} \{ v_p(A_i) + ri \} \\
&= \min_{m \leq 0} \{r v_m^{\naive}(A - I) + m\}.
\end{align*}
Apply Proposition~\ref{prop:matfact} to factor $A(I+X)^{-1}$ as $BC$, where
\[
\min\{w_r^{\naive}(B-I),w_r^{\naive}(C-I)\}
 \geq 
\min_{m \leq 0} \{r v_m^{\naive}(A-I) + m\},
\]
$B-I$ involves only positive powers of $u$, and
$C$ involves no positive powers of $u$; the desired matrices are
$U = B$ and $V = C(I+X)$.
\end{proof}

\subsection{Descending the special slope filtration}

In this section, we refine the decomposition
given by Theorem~\ref{thm:decomp} in the case of a $\sigma$-module
defined over $\Gamma^{k((t))}_{\an,\con}$, to obtain our main filtration
theorem. 

\begin{lemma} \label{lem:filt1}
For $K$ a valued field and $r>0$ satisfying the conclusion
of Proposition~\ref{prop:uniform},
let $U$ be a matrix over $\GK_{\an,r}$
and $V$ a matrix over $\GK_r$ such that $w_r(V-I) > 0$ and
$v_p(V-I) > 0$. Then
\[
\min_{m \leq 0} \{r v_m(UV - I) + m \}
= \min_{m \leq 0} \{ r v_m(U - I) + m \}.
\]
\end{lemma}
\begin{proof}
In one direction, we have
\begin{align*}
\min_{m \leq 0} \{r v_m(UV - I) +m \} &=
\min_{m \leq 0} \{r v_m((U-I)V + (V-I)) + m\} \\
&= \min_{m \leq 0} \{ r v_m((U-I)V) + m \} \\
&\geq \min_{m \leq 0,l \geq 0} \{r v_l(V) + l + r v_{m-l}(U-I) + (m-l)
\} \\
&\geq \min_{m \leq 0} \{ r v_m(U-I) + m \},
\end{align*}
the last inequality holding because $w_r(V) = 0$.
The reverse direction is implied by the above inequality with $U$ and
$V$ replaced by $UV$ and $V^{-1}$.
\end{proof}

The key calculation is the following proposition. In fact, it should be
possible to give a condition of this form that guarantees that a
$\sigma$-module has a particular special Newton polygon. However, we have
not found such a condition so far.
\begin{prop} \label{prop:filt2}
Let $K$ be a finite extension of $k((t))$ and $r>0$ a number for 
which there exists a semiunit $u$ in $\GK_{qr}$ lifting a uniformizer of $K$. 
Let $A$ be an invertible matrix over $\Gamma_{\an,r}$, and suppose
that there exists a diagonal matrix $D$ over $\calO$ such that
\[
w_r(AD^{-1} - I) > \max_{i,j} \{v_p(D_{ii}) - v_p(D_{jj}) \}.
\]
Then there exists an invertible matrix $U$ over $\Gamma_{\an,qr}$ such that
$w_r(U - I) > 0$, $U-I$ involves only positive powers of $u$,
$U^{-1} A U^\sigma D^{-1}$
 is invertible over $\Gamma_r$ and $v_p(U^{-1}A U^\sigma D^{-1} - I) > 0$.
\end{prop}
\begin{proof}
There is no loss of generality in assuming $K = k((t))$. Then by
Lemma~\ref{lem:compnaive}, for $s \leq qr$ and $x \in \Gamma_{\an,r}$, 
$w_s(x) = w_s^{\naive}(x)$ and
$\min_{m \leq 0} \{ s v_m(x) + m \} =
\min_{m \leq 0} \{ s v_m^{\naive}(x) + m \}$.
This allows us to apply the results of the previous section.

Put $c = \max_{i,j} \{ v_p(D_{ii}) - v_p(D_{jj}) \}$ and 
$d = w_r(AD^{-1}-I)$,
and define sequences $\{A_i\}$, $\{U_i\}$, $\{V_i\}$ for $i=0,1,\dots$
as follows.
Begin with $A_0 = A$. Given $A_i$,
factor $A_i D^{-1}$ as $U_i V_i$ as per Proposition~\ref{prop:matfact2},
and set $A_{i+1} = U_i^{-1} A_i U_i^\sigma$, so that
$A_{i+1} D^{-1} = V_i (D U_i^\sigma D^{-1})$.

Note that the application of Proposition~\ref{prop:matfact2} is only
valid if $w_r(A_i D^{-1} - I) > 0$.
In fact, we will show that
\[
\min_{m \leq 0} \{ r v_m(A_i D^{-1} - I) + m \} \geq d
+ i ((q-1)d - c)
\qquad \mbox{and} \qquad
w_r(A_i D^{-1} - I) \geq d-c > 0
\]
by induction on $i$. Both assertions hold for
$i=0$. Given that they hold for $i$, we
have $w_r(U_i - I) \geq \min_{m \leq 0} \{ r v_m(A_i D^{-1} - I) + m \}$
by Proposition~\ref{prop:matfact2}. On one hand, we have
\begin{align*}
w_r(D U_i^\sigma D^{-1} - I) &\geq w_{qr}(U_i - I) - c\\
&= \min_{m} \{ qr v_m^{\naive}(U_i - I) + m\} - c\\
&\geq \min_m \{ r v_m^{\naive}(U_i - I) + m\} -c\\
&= w_r(U_i - I) - c\\
&\geq \min_{m \leq 0} \{r v_m(A_i D^{-1} - I) + m \} - c\\
&\geq d - c;
\end{align*}
since $w_r(V_i - I) \geq w_r(A_i D^{-1} - I) \geq d-c$, 
we conclude $w_r(A_{i+1} D^{-1} - I) \geq d-c$.
On the other hand,
by Lemma~\ref{lem:filt1}, we have
\begin{align*}
\min_{m \leq 0} \{ r v_m(A_{i+1} D^{-1} - I) + m \} &= 
\min_{m \leq 0} \{ r v_m(V_i (D U_i^\sigma D^{-1}) - I) + m \} \\
&= \min_{m \leq 0} \{ r v_m(D U_i^\sigma D^{-1} - I) + m \} \\
&\geq \min_{m \leq 0} \{ rq v_m(U_i - I) + m\} - c \\
&\geq q \min_{m \leq 0} \{ r v_m(U_i - I) + m \} - c \\
&\geq q \min_{m \leq 0} \{ r v_m(A_i D^{-1} - I) + m \} - c \\
&\geq qd + qi((q-1)d -c) - c \\
&\geq d + (i+1)((q-1)d - c).
\end{align*}
This completes the induction and shows that the sequences are
well-defined. 

We have now shown $\min_{m \leq 0} \{ r v_m(A_{i} D^{-1} - I) + m \}
\to \infty$ as $i \to \infty$.
By Proposition~\ref{prop:matfact2},
this implies $w_r(U_i - I) \to \infty$ as $i \to \infty$, and so
$w_s(U_i - I) \to \infty$ for $s \geq r$ since $U_i-I$ involves
only positive powers of $u$.

We next consider $s \leq r$, for which $w_{s}(V_i - I)
\geq d-c>0$ for all $i$. By Lemma~\ref{lem:filt1},
\begin{align*}
\min_{m \leq 0} \{s v_m(A_{i+1} D^{-1} - I) + m\} &= 
\min_{m \leq 0} \{s v_m(V_i DU_i^\sigma D^{-1} - I) + m\} \\
&= \min_{m \leq 0} \{s v_m(DU_i^\sigma D^{-1} - I) + m\} \\
&\geq w_{sq}(U_i - I) - c.
\end{align*}
For $r/q \leq s \leq r$, we already have
$w_{sq}(U_i - I) - c \to \infty$ as $i \to \infty$, which
yields $\min_{m \leq 0}
\{ s v_m(A_{i+1} D^{-1} - I) + m \} \to \infty$ as $i \to \infty$;
by similar reasoning, $w_s(A_{i+1} D^{-1}-I) \geq d-c$ for large $i$.
By Proposition~\ref{prop:matfact2} (and the fact that the decomposition
therein does not depend on $s$), we deduce
$w_s(U_{i+1} - I) \to \infty$ as $i \to \infty$.
But now we can repeat the same line of reasoning for $r/q^2 \leq s \leq r/q$,
then for $r/q^3 \leq s \leq r/q^2$, and so on. 
Hence $w_s(U_i - I) \to \infty$ for all $s>0$.

We define $U$ as the convergent product $U_0 U_1 \cdots$;
note that $U$ is invertible because the product $\cdots U_1^{-1} U_0^{-1}$
also converges. Moreover,
\[
A_i D^{-1} = (U_0 \cdots U_{i-1})^{-1} A (U_0 \cdots U_{i-1})^\sigma D^{-1}
\]
converges to $U^{-1} A U^\sigma D^{-1}$ as $i \to \infty$. But
for $m \leq 0$, we already have $r v_m(A_i D^{-1} - I) + m \to \infty$ as $i
\to \infty$, so $v_m(U^{-1} A U^\sigma D^{-1} - I) = \infty$. Hence
$U^{-1} A U^\sigma D^{-1}$ and its inverse have entries in $\Gamma_r$ and is
congruent to $I$ modulo $\pi$, as desired.
\end{proof}

This lemma, together with the results of the previous chapters,
allows us to deduce an approximation to our desired result,
but only so far over an unspecified finite extension of $k((t))$.
\begin{prop} \label{prop:sameslope}
Let $M$ be a $\sigma$-module over $\Gancon = \Gamma^{k((t))}_{\an,\con}$
whose special Newton slopes lie in the value group
of $\calO$.
Then there exists a finite extension $K$
of $k((t))$ such that
$M \otimes_{\Gancon} \GKancon$ is isomorphic to
$M_1 \otimes_{\GKcon} \GKancon$
for some $\sigma$-module $M_1$ over $\GKcon[\fp]$ whose generic and special
Newton polygons coincide.
\end{prop}
If $k$ is perfect, we can take $K$ to be separable over $k((t))$, but this
is not necessary for our purposes.
\begin{proof}
Pick a basis of $M$ and let $A$ be the matrix via which
$F$ acts on this basis. By Theorem~\ref{thm:decomp}, 
there exists an invertible matrix $X$ over $\Galgancon$ such that
$A = XDX^{-\sigma}$ for some diagonal matrix $D$ over $\calO$.
Choose $r>0$ such that $A$ is invertible over $\Gamma_r$ and $X$ is invertible
over $\Galg_{\an,rq}$.

Choose $c > \max_{ij} \{v_p(D_{ii}) - v_p(D_{jj})\}$.
By Lemma~\ref{lem:approx2} applied to $X^T$,
there exists a finite extension $K$
of $k((t))$ and an invertible matrix $V$ over $\GK_r[\fp]$ such that
$w_l(VX - I) \geq 2c$ for $r \leq l \leq qr$.
By replacing $K$ by a suitable inseparable extension, we can ensure that
$\GK_{qr}$ contains a semiunit lifting a uniformizer of $K$.

Observe that
\[
(VAV^{-\sigma})D^{-1} = (VX) D (VX)^{-\sigma} D^{-1}.
\]
Since $w_r(VX - I) \geq 2c$ and
\[
w_{r}(D (VX)^{-\sigma} D^{-1} - I) \geq 
w_{qr}(VX-I) - c \geq c,
\]
we have $w_r(VAV^{-\sigma}D^{-1} - I) \geq c$.
By Proposition~\ref{prop:filt2}, there exists an invertible matrix $U$ over 
$\GK_{\an,qr}$
such that $U^{-1}VAV^{-\sigma}U^\sigma D^{-1}$ has entries in
$\GK_r$ and is congruent to $I$ modulo $\pi$.
Put $W = V^{-1}U$; then we
can change basis in $M$ so that $F$ acts on the new basis
via the matrix $W^{-1}AW^\sigma$. Let $M_1$ be
the $\GKcon[\fp]$-span of the basis elements;
by Proposition~\ref{prop:genspec}, the
generic Newton slopes of $M_1$ are the valuations
of the entries of $D$, so they coincide with the special Newton slopes.
Thus $M_1$ is the desired $\sigma$-module.
%
\end{proof}

By descending a little bit more, we now deduce the main result of the
paper,
a slope filtration theorem for $\sigma$-modules over the Robba ring.
\begin{theorem} \label{thm:filt}
Let $M$ be a $\sigma$-module over $\Gancon = \Gamma^{k((t))}_{\an,\con}$.
Then there is a filtration
$0 = M_0 \subset M_1 \subset \cdots \subset M_l = M$
of $M$ by saturated $\sigma$-submodules such that:
\begin{enumerate}
\item[(a)] for $i=1, \dots, l$, the quotient $M_{i}/M_{i-1}$ has a single
special slope $s_i$;
\item[(b)] $s_1 < \cdots < s_{l}$;
\item[(c)] each quotient $M_i/M_{i-1}$ contains an $F$-stable
$\Gcon[\fp]$-submodule $N_i$ of the same rank, which spans $M_{i}/M_{i-1}$
over $\Gancon$, and which has all generic slopes equal to $s_i$.
\end{enumerate}
Moreover, conditions (a) and (b) determine the filtration uniquely,
and the $N_i$ in (c) are also unique.
\end{theorem}
\begin{proof}
Let $s_1$ be the lowest special slope of $M$ and $m$ its multiplicity.
We prove that there exists a saturated $\sigma$-submodule $M_1$ 
of rank $m$ whose special
slopes all equal $s_1$,
that $M_1$ contains a $F$-stable
$\Gcon[\fp]$-submodule $N_1$ of the same rank, which spans $M_1$
over $\Gancon$, and whose generic slopes equal to $s_1$,
and that these properties uniquely characterize $M_1$ and $N_1$.
This implies the desired result by induction on the rank of $M$.
(Once $M_1$ is constructed, apply the induction hypothesis to $M/M_1$.)

We first establish the existence of $M_1$. Let $\calO'$ be a Galois
extension of $\calO$ whose value group contains all of the special 
slopes of $M$. By Proposition~\ref{prop:sameslope},
for some valued field $K$ finite and normal over $k((t))$, $M$ is isomorphic
over $\GKancon \otimes_{\calO} \calO'$
to a $\sigma$-module $M'$ defined over $\GKcon[\fp] \otimes_{\calO} \calO'$ whose
generic and special Newton polygons are equal. By Proposition~\ref{prop:asc2},
$M'$ admits an ascending slope filtration over $\GKcon \otimes_{\calO} \calO'$,
so $M$ admits one over $\GKancon \otimes_{\calO} \calO'$; let $Q_1$ and $P_1$
be the respective first steps of these filtrations.
Then the slope of $P_1$ is $s_1$ with multiplicity $m$.
Moreover, the top exterior power of $P_1$ is defined
both over $\Galgancon$ (because the lowest slope of $\wedge^m M$ is $s_1 m$,
which is in the value group of $\calO$)
and over $\GKancon \otimes_{\calO} \calO'$, and
hence over their intersection $\GKancon$. Thus $P_1$ is defined over
$\GKancon$. 

Let $K_1$ be the maximal purely inseparable subextension of $K/k((t))$
(necessarily a valued field),
and let $M_1$ be the saturated span of the images of
$P_1$ under $\Gal(K/K_1)$; by Corollary~\ref{cor:galois2}, $M_1$ descends
to $\Gamma^{K_1}_{\an,\con}$, and its rank is at least $m$.
Moreover, over $\Galgancon \otimes_{\calO} \calO'$, $M_1$
is spanned by eigenvectors of slope $s_1$,
so the special slopes of $M_1$ are all at most $s_1$
by Proposition~\ref{prop:basis}. 
Thus $M_1$ has the single slope $s_1$ with multiplicity $m$.

We must still check that $M_1$ descends
from $\Gamma^{K_1}_{\an,\con}$ to $\Gancon$. Let $\be_1, \dots, \be_n$ be
a basis of $M$ and let $\bv_1, \dots, \bv_m$ be a basis of $M_1$. Then
we can write $\bv_i = \sum_{j} c_{ij} \be_j$ for some $c_{ij} \in 
\Gamma^{K_1}_{\an,\con}$. Since $K_1/k((t))$ is purely inseparable,
$K_1^{q^d} \subseteq k((t))$ for some integer $d$; for any such $d$,
$F^d \bv_1, \dots, F^d \bv_m$ is a basis of $M_1$ and
$F^d \bv_i = \sum_j c_{ij}^{\sigma^d} F^d \be_j$. Since
each $c_{ij}^{\sigma^d}$ belongs to $\Gancon$, each $F^d\bv_i$ belongs to $M$;
thus $M_1$ descends to $\Gancon$.

We next establish existence of an $F$-stable $\Gcon[\fp]$-submodule
$N_1$ of $M_1$, having the same rank and spanning $M_1$ over $\Gancon$,
and having all generic slopes equal to $s_1$.
Note that $Q_1$, defined above, is an $F$-stable $(\GKcon[\fp] \otimes_{\calO}
\calO')$-submodule 
of $M_1 \otimes_{\Gancon} \GKancon \otimes_{\calO} \calO' = P_1$ 
with the properties desired of $N_1$.
Moreover, $Q_1 \otimes_{\GKcon} \Galgcon$ is equal to
the $(\Galgcon[\fp] \otimes_{\calO} \calO')$-span of the eigenvectors of 
$M$ of slope $s_1$, which is invariant under $\Gal(k((t))^{\alg}
/ k((t))^{\perf}) \times \Gal(\calO'/\calO)$. Thus
$Q_1$ is invariant under $\Gal(K/K_1) \times \Gal(\calO'/\calO)$; by
Galois descent, it descends to $\Gamma^{K_1}_{\con}[\fp]$, and thus
to $\Gcon[\fp]$ (again, by applying Frobenius repeatedly). This yields
the desired $N_1$.

With the existence of $M_1$ and $N_1$ in hand, we check uniqueness.
For $M_1$, note that $M_1 \otimes_{\Gancon} \Galgancon$ is equal
to the $(\Galgancon \otimes_{\calO} \calO')$-span of the eigenvectors
of $M$ of slope $s_1$, because
otherwise some eigenvector of slope $s_1$ would survive quotienting by $M_1$,
contradicting Proposition~\ref{prop:minslope} because the quotient has
all slopes greater than $s_1$.
This description uniquely determines $M_1$. For $N_1$, note
that $N_1 \otimes_{\Gcon} \Galgcon \otimes_{\calO} \calO'$ is
equal to the $(\Galgcon[\fp] \otimes_{\calO} \calO')$-span of the eigenvectors
of $M$ of slope $s_1$, because it contains a basis of eigenvectors of
slope $s_1$ by Proposition~\ref{prop:descfilt}.
This description uniquely determines $N_1$.

Thus $M_1$ and $N_1$ exist and are unique; as noted above, induction on
the rank of $M$ now completes the proof.
\end{proof}
One consequence of this proposition is that if $k$ is perfect,
the lowest slope eigenvectors of
a $\sigma$-module over $\Gancon$ are defined not just over $\Galgancon$,
but over the subring $\Gancon \otimes_{\Gcon} \Gsepcon$. (If $k$ is 
not perfect, then $\Gsepcon$ may not be defined, but we can replace it
with $\Galgcon$ to get a weaker but still nontrivial statement.)

\subsection{The connection to the unit-root case}

In this section, we deduce Theorem~\ref{thm:monodromy} from 
Theorem~\ref{thm:filt}. To exploit the extra data of a connection 
provided by a $(\sigma, \nabla)$-module, we invoke Tsuzuki's
finite monodromy theorem for unit root $F$-crystals 
\cite[Theorem~5.1.1]{bib:tsu1}, as follows.
(Another proof of the theorem appears in \cite{bib:ch}, and yet
another in \cite{bib:methesis}. However, none of these proves the
theorem at quite the level of generality we seek, so we must fiddle a bit
with the statement.)

Recall that a valued field $K/k((t))$ is said to be \emph{nearly separable}
if it is a separable extension of $k^{1/p^m}((t))$ for some nonnegative
integer $m$ (and that not all separable extensions of $k((t))$ are valued
fields).
\begin{prop} \label{prop:tsuzuki}
Let $M$ be a unit-root $(\sigma, \nabla)$-module of rank $n$ over $\Gcon = 
\Gamma^{k((t))}_{\con}$. 
For any finite extension $K$ of $k((t))$, if
there exists a basis of $M \otimes_{\Gcon} \GKcon$
on which $F$ acts via a matrix $A$
with $v_p(A - I) > 1/(p-1)$, then 
the kernel of $\nabla$ on $M \otimes_{\Gcon} \GKcon$ has rank $n$
over $\calO$ and is $F$-stable.
Moreover, such a $K$ can always be chosen which is separable over
$k((t))$ if $k$ is perfect, or nearly
separable if $k$ is imperfect.
\end{prop}
\begin{proof}
The theorem of Tsuzuki \cite[Theorem~5.1.1]{bib:tsu1} establishes
the first assertion for $k$ algebraically closed and $q=p$; in fact,
it produces a basis of eigenvectors in the kernel of $\nabla$. 
The first assertion in general follows from 
this case by a relatively formal argument, given below.
Note that the kernel of $\nabla$ is always $F$-stable, so we do not have
to establish this separately.

We first allow $q = p^f$, still with $k$ algebraically closed.
Let $M$ be a unit-root $(\sigma,
\nabla)$-module over $\Gamma_{\con}$; recall that the Frobenius
structure can be described as a $\Gcon$-linear isomorphism
$F: M \otimes_{\Gcon,\sigma} \Gcon \to M$. For $i=0, \dots, f$, put $M_i = M
\otimes_{\Gcon, \sigma_0^i} \Gcon$, where $\Gcon$
is viewed as a module over itself via $\sigma_0^i$. Then
$M_0 \oplus \cdots \oplus M_{f-1}$ admits the structure of
a unit-root $(\sigma_0, \nabla_0)$-module as follows.
The connection $\nabla_0$ acts factorwise, with the component of $M_i$ being
\[
\nabla \otimes \mathrm{id}_{\Gcon}:
M \otimes_{\Gcon, \sigma_0^i} \Gcon \to (M_i 
\otimes_{\Gcon} \Omega^1_{\Gcon/\calO})
\otimes_{\Gcon, \sigma_0^i} \Gcon.
\]
 The Frobenius map
\[
F_0: (M_0 \oplus \cdots \oplus M_{f-1}) \otimes_{\Gcon,\sigma_0}
\Gcon
\cong M_1 \oplus \cdots \oplus M_f \to
M_0 \oplus \cdots \oplus M_{f-1}
\]
carries $M_i$ to $M_i$ for $i=1, \dots, f-1$ and maps $M_f \cong
M \otimes_{\Gcon, \sigma} M$
to $M_0 \cong M$ via the original $F$.
By Tsuzuki's theorem,
this module admits a basis of eigenvectors in the kernel of
$\nabla$ over $\GKcon$ for some $K$; projecting these eigenvectors
onto the first factor gives a basis of $M$ consisting of
elements in the kernel of $\nabla$. 

We now treat general $k$ by a ``compactness'' argument.
For simplicity of notation, let us assume
$K = k((t))$, and let $\calO'$ be the completion of the maximal unramified
extension of $\calO$. Then Tsuzuki's theorem (plus the above argument
if $q \neq p$) provides a basis
$\bv_1, \dots, \bv_n$ of
the kernel of $\nabla$ over $\Gamma^{k^{\alg}((t))}_{\con}$,
and we must produce a basis of the kernel of $\nabla$ over $\Gcon$.
Let $\be_1, \dots, \be_n$ be a basis of $M$
and put $\bv_i = \sum_{j,l} c_{i,j,l} u^l \be_j$.
Put $d_{j,l} = \min_i \{ v_p(c_{i,j,l})/v_p(\pi)\}$, and whenever
$d_{j,l} < \infty$, write $c_{i,j,l}$ as $\pi^{d_{j,l}} f_{i,j,l}$.

The fact that $\nabla \bv_i = 0$ for $i=1, \dots, n$ can be rewritten
as a set of ``quasilinear'' equations in the $f_{i,j,l}$. That is,
for $h=1,2,\dots$, we have equations of the form
\[
\sum_{i,j,l} g_{h,i,j,l} f_{i,j,l} = 0
\]
for certain $g_{h,i,j,l} \in \calO$, such that for any $h$ and $m$, only
finitely many of the $g_{h,i,j,l}$ are nonzero modulo $\pi^m$. We are
given that these equations have $n$ linearly independent solutions over
$\calO'$, and wish to prove they have $n$ linearly independent solutions
over $\calO$.

For each finite set $S$ of triples $(i,j,l)$,
let $T_{S}(\calO)$ (resp.\ $T_{S}(\calO')$) be the set of functions 
$f: S \to \calO$ (resp.\ $f: S \to \calO'$), mapping a pair $(i,j,l) \in S$
to $f_{i,j,l}$,
which can be extended to a simultaneous solution of any finite subset of
the equations modulo any power of $\pi$.
If we put the $T_{S}$ into an inverse system under inclusion on $S$,
then the restriction
maps are all surjective, and solutions to the complete set of
equations are precisely elements of the inverse limit.
However, each equation modulo each power of $\pi$
involves only finitely many variables, so
$T_{S}$ is defined by linear conditions on the $f_{i,j,l}$. Thus
$T_{S}(\calO') = T_{S}(\calO) \otimes_{\calO} \calO'$. Since
the solutions of the system over $\calO'$ have rank $n$, we have
$\rank_{\calO'} T_{S}(\calO') = n$ for $S$ sufficiently large.
Thus the same holds over $\calO$, which
produces $n$ $\calO$-linearly 
independent elements of the inverse limit, hence
of the kernel of $\nabla$ over $\Gcon$.
This establishes the first assertion of the proposition for general $k$.

Finally, we show that $K$ can be taken to be (nearly) separable over
$k((t))$. By Proposition~\ref{prop:ascfilt}
(where the \textit{ad hoc} definition of $\Gsep$ was given),
$M \otimes_{\Gcon} \Gsep$ 
admits a basis up to isogeny of eigenvectors $\bw_1, \dots, \bw_n$.
By the Dieudonn\'e-Manin classification in the form of
Proposition~\ref{prop:geneigs} (and the fact that the unique
slope is already in the value group), the kernel of $\nabla$
on $M \otimes_{\Gcon} \GKcon$ admits a basis up to isogeny of eigenvectors 
over some unramified extension $\calO'$ of $\calO$; by the proof of 
Proposition~\ref{prop:ascfilt}, the residue field extension of $\calO'$
over $\calO$ is separable. Thus $\calO' \subseteq \Gsep$, and so
each $\bv_i$ in the kernel of $\nabla$
is a $\Gsep[\fp]$-linear combination of the $\bw_i$. Hence
the $\bv_i$ are defined over $\Gsep[\fp] \cap \GKcon$. If $k$ is perfect,
this intersection equals $\Gamma^{K_1}_{\con}$
for $K_1$ the maximal separable subextension of $K$ over $k((t))$. If
$k$ is imperfect, $K_1$ may fail to be a valued field. Instead,
choose an integer $i$ for which
the maximal
purely inseparable subextension of the residue field extension of
$K_1$ over $k((t))$ is contained in $k^{1/p^i}$. Then the compositum
$K_2$ of $K_1$ and $k^{1/p^i}((t))$ is separable and totally ramified
over $k^{1/p^i}((t))$, so is a nearly separable valued field,
and the $\bv_i$ are defined over $\Gamma^{K_2}_{\con}$, as desired.
\end{proof}

Theorem~\ref{thm:monodromy} follows immediately from the following theorem,
which refines the results of Theorem~\ref{thm:filt} in the presence
of a connection, using Tsuzuki's theorem.
\begin{theorem} \label{thm:main2}
Let $M$ be a $(\sigma, \nabla)$-module over $\Gancon = 
\Gamma^{k((t))}_{\an,\con}$. Then the filtration of
Theorem~\ref{thm:filt} satisfies the following additional properties:
\begin{enumerate}
\item[(d)] each $M_i$ is a $(\sigma, \nabla)$-submodule;
\item[(e)] each $N_i$ is $\nabla$-stable;
\item[(f)] there exists a finite nearly separable extension $K/k((t))$ 
(separable
in case $k$ is perfect) such that each $N_i$
is spanned by the kernel of $\nabla$ over $\GKcon[\fp]$;
\item[(g)] if $k$ is algebraically closed, $N_i$
is isomorphic over $\GKcon[\fp]$ to a direct sum of standard
$(\sigma, \nabla)$-modules.
\end{enumerate}
\end{theorem}
\begin{proof}
Again by induction on the rank of $M$, it suffices to prove (d), (e), (f),
(g)
for $i=1$.
For (d) and (e), we may assume without loss of generality (by enlarging
$\calO$, then twisting) that
the special slopes of $M$ belong to the value group of $\calO$ and that
$s_1 = 0$.

By Proposition~\ref{prop:nodenom}, we can
choose a basis for $N_1$ on which $F$ acts by an invertible matrix
$X$ over $\Gcon$. Extend this basis to a basis of $M$; then
$F$ acts on the resulting basis via some block matrix over $\Gancon$
of the form $\begin{pmatrix} X & Y \\ 0 & Z \end{pmatrix}$.
View $\nabla$ as a map from $M$ to itself by identifying $x \in M$
with $x \otimes du \in M \otimes_{\Gancon} 
\Omega^1$; then
$\nabla$ acts on the chosen basis of $M$ by some block matrix
$\begin{pmatrix} P & Q \\ R & S \end{pmatrix}$ over $\Gancon$.
The relation
$\nabla \circ F = (F \otimes d\sigma) \circ \nabla$ translates into
the matrix equation
\[
\begin{pmatrix} P & Q \\ R & S \end{pmatrix}
\begin{pmatrix} X & Y \\ 0 & Z \end{pmatrix}
+
\frac{d}{du} \begin{pmatrix} X & Y \\ 0 & Z \end{pmatrix}
=
\frac{du^\sigma}{du}
\begin{pmatrix} X & Y \\ 0 & Z \end{pmatrix}
\begin{pmatrix} P & Q \\ R & S \end{pmatrix}^\sigma.
\]
The lower left corner of the matrix equation
yields $RX = \frac{du^\sigma}{du}
ZR^\sigma$. We can write $X = U^{-1} U^\sigma$ with
$U$ over $\Galgcon$ by Proposition~\ref{prop:descfilt} (since $M_1$
has all slopes equal to 0)
and $Z = V^{-1} D V^\sigma$ with $V$ over $\Galgancon$ and $D$
a scalar matrix over $\calO$ whose entries have positive valuation
(because $M_1$ is the lowest slope piece of $M$).
We can write $\frac{du^\sigma}{du} = \mu x$ for some $\mu \in \calO$
and $x$ an invertible element of $\Gcon$; since $u^\sigma \equiv u^q
\pmod{\pi}$, we have $|\mu| < 1$. By Proposition~\ref{prop:rank1}, there
exists $y \in \Galgcon$ nonzero such that $y^\sigma = xy$.
Now rewrite the equation $RX = \frac{du^\sigma}{du}
ZR^\sigma$ as 
\[
yVRU^{-1} = \mu D(yVRU^{-1})^\sigma;
\]
by Proposition~\ref{prop:ansigeq}(c)
applied entrywise to this matrix equation, we deduce
$yVRU^{-1} = 0$ and so $R=0$.
In other words, $M_1$ is stable under $\nabla$, and (d) is verified.

We next check that $N_1$ is $\nabla$-stable;
this fact is due to Berger \cite[Lemme~V.14]{bib:berger}, but
our proof is a bit different.
Put $X_1 = \frac{dX}{du}$;
then the top left corner of the matrix equation yields
$PX + X_1 = \frac{du^\sigma}{du} X P^\sigma$, or
\[
y UPU^{-1} + y U X_1 U^{-\sigma} = \mu (y UPU^{-1})^\sigma.
\]
By Proposition~\ref{prop:ansigeq}(c),
each entry of $y UPU^{-1}$ lies in $\Galgcon$,
so the entries of $P$ lie in $\Galgcon \cap \Gancon = \Gcon$. Thus
$N_1$ is stable under $\nabla$, and (e) is verified.

To check (f), we must relax the simplifying assumptions. If they
do happen to hold, then $N_1$
is a unit-root $(\sigma, \nabla)$-module over $\Gcon$, so for 
some finite (nearly) separable extension $K$ of $k((t))$, the kernel of $\nabla$
on $N_1 \otimes_{\Gcon} \GKcon$ has full rank.
Without the simplifying assumptions, we only have that
the kernel of $\nabla$ has full rank in
$N_1 \otimes_{\Gcon} \GKcon \otimes_{\calO} \calO'$ for some finite
extension $\calO'$ of $\calO$. However, decomposing kernel elements
with respect to a basis of $\calO'$ over $\calO$ produces elements of the
kernel of $\nabla$ in $N_1 \otimes_{\Gcon} \GKcon$ which span $M$, so
the kernel has full rank over $N_1 \otimes_{\Gcon} \GKcon$. Thus
(f) is verified.

Finally, suppose $k$ is algebraically closed.
As noted in the proof of Proposition~\ref{prop:tsuzuki}, the kernel
of $\nabla$ is always $F$-stable. By the Dieudonn\'e-Manin classification
(Theorem~\ref{thm:dm}),
it is isogenous as a $\sigma$-module to a direct sum of standard 
$\sigma$-modules. This gives a decomposition of $N_1 \otimes_{\Gcon}
\GKcon$
as a direct sum of standard $(\sigma, \nabla)$-modules. Thus
(g) is verified and the proof is complete.
\end{proof}

\subsection{Logarithmic form of Crew's conjecture}

An alternate formulation of the local monodromy theorem can be given,
that eschews the filtration and instead describes a basis of the original
module given by elements of the kernel of $\nabla$. The tradeoff is that 
these elements are defined not over a Robba ring, but over a ``logarithmic''
extension thereof. As this is the most useful formulation in some
applications, we give it explicitly.

For $r>0$,
the series $\log(1+x) = x - x^2/2 + \cdots$ converges under $|\cdot|_r$
whenever $|x|_r < 1$. Thus if $x \in \Gcon$ satisfies $|x-1|_r<1$,
then $\log(1+x)$ is well-defined and $\log (1+x+y+xy)
= \log(1+x) + \log(1+y)$.

For any valued field $K$ finite over $k((t))$,
we produce the ring $\Gamma^K_{\log,\con}$ (resp.\ $\Gamma^K_{\log,\an,\con}$)
from $\Gamma^K_{\con}$ (resp.\ $\Gamma^K_{\an,\con}$)
by adjoining one variable $l_z$
for each $z \in \Gamma^K_{\con}$ not divisible by
$\pi$, subject to the relations
\begin{align*}
  l_{z^n} &= n l_z \qquad (n \in \ZZ) \\
  l_{z(1+x)} &= l_z + \log (1+x) \qquad (v_0(x) > 0)
\end{align*}
and with the operations $\sigma$ and $\frac{d}{du}$,
for $u \in \GK_{con}$ nonzero, extended as follows:
\begin{align*}
  (l_z)^\sigma &= l_{z^\sigma} \\
  \frac{dl_z}{du} &= \frac{1}{z} \frac{dz}{du}.
\end{align*}
Note that as a ring, $\GK_{\log,\con}$ (resp.\ $\Gamma^K_{\log,\an,\con}$)
is isomorphic to the polynomial ring $\GK_{\con}[l_z]$ (resp.\ $\Gamma^K_{\an,\con}[l_z]$)
for any one $z \in \GKcon$ which lifts a uniformizer of $K$.
We will suggestively write $\log z$ instead of $l_z$.

\begin{theorem}
Let $M$ be a $(\sigma, \nabla)$-module over 
$\Gancon = \Gamma^{k((t))}_{\an,\con}$. 
Then for some finite (nearly) separable extension $K$ of $k((t))$,
$M$ admits a basis over $\GK_{\log,\an,\con}$ 
of elements of the kernel of $\nabla$.
Moreover, if $k$ is algebraically closed, 
$M$ can be decomposed over $\GK_{\log,\an,\con}$ as the direct sum of
standard $(\sigma, \nabla)$-submodules.
\end{theorem}
\begin{proof}
By Theorem~\ref{thm:main2}, there exists a basis $\bv_1, \dots, \bv_n$
of $M$ over $\GKancon$, for some finite nearly separable extension $K$ of $k((t))$,
such that
$\nabla \bv_i \in \Span(\bv_1, \dots, \bv_{i-1}) \otimes 
\Omega^1$.
Choose a lift $u \in \GKcon$ of a uniformizer of $K$, view
$\nabla$ as a map from $M$ to itself by identifying $\bv \in M$ with $\bv \otimes
du$,
and write 
$\nabla \bv_i = \sum_{j<i} A_{ij} \bv_i$ for some $A_{ij} \in \GKancon$.

Define a new basis $\bw_1, \dots, \bw_n$ of $M$ over
$\GK_{\log,\an,\con}$ as follows.
First put $\bw_1 = \bv_1$. Given $\bw_1, \dots, \bw_{i-1}$ with the same
span as $\bv_1, \dots, \bv_{i-1}$ such that $\nabla \bw_j = 0$
for $j=1, \dots, i-1$, put
$\nabla\bv_j = c_{i,1} \bw_1 + \cdots + c_{i,i-1} \bw_{i-1}$ and
write $c_{i,j} = \sum_{l,m} d_{i,j,l,m} u^l (\log u)^m$.
Now recall from calculus that every expression of the form
$u^l (\log u)^m$, with $m$ a nonnegative integer, can be written as the 
derivative with respect to $u$
of a linear combination of such expressions. (If $l = -1$, the expression
is the derivative of a power of $\log u$ times a scalar. Otherwise,
integration by parts can be used to
reduce the power of the logarithm.) Thus there exist
$e_{i,j} \in \GK_{\log,\an,\con}$ such that $\frac{d}{du} e_{i,j} = c_{i,j}$.
 Put
$\bw_i = \bv_i - \sum_{j<i} e_{i,j} \bw_j$; then $\nabla \bw_i = 0$.
This process thus ends with a basis $\bw_1, \dots, \bw_n$
of elements of the kernel of $\nabla$.

As in the proof of Proposition~\ref{prop:tsuzuki}, the kernel of $\nabla$
is $F$-stable.
Thus if $k$ is algebraically closed, we may apply the Dieudonn\'e-Manin
classification (Theorem~\ref{thm:dm}) 
to decompose $M$ over $\GK_{\log,\an,\con}$
as the sum of standard $(\sigma, \nabla)$-modules, as desired.
\end{proof}

\end{document}